\documentclass[final,numbib]{imaiai}

\usepackage{natbib}
\usepackage[utf8]{inputenc}
\usepackage[T1]{fontenc}
\usepackage[english]{babel}
\usepackage[table]{xcolor}
 \definecolor{burgundy}{rgb}{0.5, 0.0, 0.13}
  \definecolor{dark_burgundy}{rgb}{0.6, 0.0, 0.05}
\definecolor{camel}{rgb}{0.76, 0.6, 0.42}
\definecolor{chamoisee}{rgb}{0.63, 0.47, 0.35}
\definecolor{grey1}{RGB}{128,128,128}
\usepackage{url}            
\usepackage{booktabs}       
\usepackage{graphics}
\usepackage{epstopdf}
\usepackage{amsmath,amsthm}
\numberwithin{equation}{section}
\usepackage{tikz}
\usetikzlibrary{positioning}
\usepackage{dsfont}
\usepackage{enumitem}
\usepackage{caption}
\usepackage{booktabs}
\usepackage{manfnt}

\usepackage{stmaryrd}

\usepackage{algorithm} 
\usepackage{algorithmic}  
\usepackage[algo2e,norelsize,ruled,vlined,commentsnumbered,linesnumbered]{algorithm2e}

\makeatletter
\newcommand{\nosemic}{\renewcommand{\@endalgocfline}{\relax}}
\newcommand{\dosemic}{\renewcommand{\@endalgocfline}{\algocf@endline}}
\let\oldnl\nl
\newcommand{\nonl}{\renewcommand{\nl}{\let\nl\oldnl}}
\makeatother

\def\beq{\begin{equation}}
\def\eq{\begin{equation}}
\def\eeq{\end{equation}}
\def\qe{\end{equation}}
\def\beqn{\begin{eqnarray*}}
\def\eeqn{\end{eqnarray*}}
\def\bitem{\begin{itemize}}
\def\eitem{\end{itemize}}
\def\benum{\begin{enumerate}}
\def\eenum{\end{enumerate}}
\def\bmult{\begin{multline*}}
\def\emult{\end{multline*}}
\def\bcenter{\begin{center}}
\def\ecenter{\end{center}}

\DeclareMathOperator*{\argmax}{arg\, max}

\DeclareMathOperator{\rank}{rank}

\def\cA{\mathcal{A}}

\def\cD{\mathcal{D}}

\def\cI{\mathcal{I}}

\def\cP{\mathcal{P}}

\def\cS{\mathcal{S}}
\def\cT{\mathcal{T}}
\def\cU{\mathcal{U}}

\def\bF{\boldsymbol{F}}

\def\bT{\boldsymbol{T}}
\def\bU{\boldsymbol{U}}

\def\bb{\mathbf{b}}

\def\bp{\mathbf{p}}

\def\bz{\mathbf{z}}

\def\bbE{\mathds{E}}
\def\bbF{\mathds{F}}

\def\bbN{\mathds{N}}

\def\bbP{\mathds{P}}

\def\bbR{\mathds{R}}

\def\bb\bU{\mathds{\bU}}

\newcommand{\E}{\operatorname{\mathds{E}}}
\renewcommand{\P}{\operatorname{\mathds{P}}}
\renewcommand{\bar}{\overline}
\renewcommand{\hat}{\widehat}
\renewcommand{\tilde}{\widetilde}

\newcommand\indep{\protect\mathpalette{\protect\independenT}{\perp}}
\def\independenT#1#2{\mathrel{\rlap{$#1#2$}\mkern2mu{#1#2}}}

\def\\bUnif{\text{\bUnif}}

\newcommand{\1}{\mathds{1}}

\begin{document}

\title{Multiple Testing and Variable Selection along\\ the path of the Least Angle Regression}

\shorttitle{GtSt: Generalized t-Spacing test} 
\shortauthorlist{Aza\"is and De Castro} 

\author{{
\sc Jean-Marc Aza\"is}\\[2pt]
Institut de Math\'ematiques de Toulouse\\ Universit\'e Paul Sabatier, 118 route de Narbonne, F-31062 Toulouse, France\\
{jean-marc.azais@univ-toulouse.fr}
\\[6pt]
{\sc and}\\[6pt]
{\sc Yohann De Castro}$^*$,\\[2pt]
Institut Camille Jordan UMR 5208, \'Ecole Centrale Lyon
 \\ 36 Avenue Guy de Collongue, F-69134 Écully, France\\
$^*${\email{Corresponding author: yohann.de-castro@ec-lyon.fr}}}

\maketitle

\begin{abstract}
{We investigate multiple testing and variable selection using the Least Angle Regression (LARS) algorithm in high dimensions under the assumption of Gaussian noise.
LARS is known to produce a piecewise affine solution path with change points referred to as the {\it knots of the LARS path}. The key to our results is an expression in closed form of the exact joint law of a $K$-tuple of knots conditional on the variables selected by LARS, the so-called {\it post-selection} joint law of the LARS knots. Numerical experiments demonstrate the perfect fit of our findings.

This paper makes three  main contributions.
First, we build testing procedures on variables entering the model along the LARS path in the general design case when the noise level can be unknown.
These testing procedures are referred to as the Generalized $t$-Spacing tests~(GtSt) and we prove that they have an exact non-asymptotic level (i.e., the Type I error is exactly controlled).
This extends work of \cite{taylor2014exact} where the spacing test works for consecutive knots and known variance.
Second, we introduce a new exact multiple testing procedure after model selection in the general design case when the noise level may be unknown.
We prove that this testing procedure has exact non-asymptotic level for general design and unknown noise level.
Third, we prove exact control of the false discovery rate under orthogonal design assumption.
Monte Carlo simulations and a real data experiment are provided to illustrate our results in this case.
Of independent interest, we introduce an equivalent formulation of the LARS algorithm based on a recursive function.}
{{Multiple Testing}; {False Discovery Rate}; {High-Dimension}; {Selective Inference}.}
\\
2000 Math Subject Classification: Primary {62E15}, {62F03}, {60G15}, {62H10}, {62H15}; secondary  {60E05}; {60G10}; {62J05}; {94A08}
\end{abstract}


\section{Introduction}
In the past decades, statistical problems have become increasingly high-dimensional, {\it i.e.,} they require estimation of more parameters than the number of available samples/observations. Some examples range from signal processing \citep{Chen_Donoho_Saunders98,Candes_Romberg_Tao06} to genomics \citep{rhee2006genotypic,barber2015controlling}. Some successful techniques of estimation have been developed and a popular approach is based on optimizing a suitable regularized likelihood function. Most models of statistical parameters are well approximated by sparse vectors; and sparsity promoting regularizations, such as the $\ell_1$-norm, are now well recognized to tackle high-dimensional problems. Recent advances have focused on a deeper understanding of the law of the estimates of $\ell_1$-regularization procedures in high-dimension. One goal is to quantify the uncertainty of some linear statistic of the outcomes of sparse regression estimation. Such estimators are non-linear and non-explicit. They are defined as the minimum of some optimization program, or as the outcomes of some greedy method. Most of them estimate some set of relevant parameters, {\it i.e.,} a small number of parameters that may explain the observation. This non-linear framework makes it impossible to characterize the distribution of the estimator. One possibility is to look at some conditional distribution of the estimator and this is the scope of the so-called selective inference, which produces an uncertainty quantification conditional on the set of indices of nonzero estimated parameters, referred to as the selection event. Selective inference aims at building some confidence intervals and some testing procedures on the estimates (see~\cite[Chapter 6]{van2016estimation} and references therein), or controlling the false discovery rate, {\it e.g.},~\cite{barber2015controlling} for instance. 

One of the most popular regularized estimation procedure in high-dimensions is LASSO \citep{Chen_Donoho_Saunders98} and its asymptotic de-biased version referred to as the debiased LASSO. Controlling the FDR (resp., confidence intervals~(CI)) built upon the debiased LASSO procedure has been studied in \cite{javanmard2019false} (resp., \cite{javanmard2014confidence}) which provides an FDR with asymptotic control (resp., the CI with asymptotic control of the confidence level) for designs with some independent sub-Gaussian rows. The LASSO is based on $\ell_1$-norm regularization and one of its offsprings is the sorted-$\ell_1$ regularization, referred to as the SLOPE, which achieves minimax rate of prediction and estimation. Controlling the FDR for SLOPE with the Benjamini-Hochberg~(BH) selection procedure has been achieved in \cite{bogdan2015slope} for orthogonal designs. 

Inference after model selection has been studied in several papers, such as~\cite{fithian2014optimal,taylor2015statistical} (resp.,~\cite{tian2018selective}) for selective inference (resp., for a joint estimate of the noise level). These works give the non-asymptotic law of any linear statistics, {\it i.e.,} any linear combination of the estimates of the parameters, conditional on the selection event. For the first time, this paper provides the non-asymptotic joint law of {\bf several linear statistics} conditional on the selection event. These linear statistics are given by the knots of the LARS procedure. One may note that, conditional on the selection event, the law of three {consecutive} knots has been studied by~\cite{lockhart2014significance} who refer to it as the spacing test (ST)~\citep{hastie2015statistical}. The article~\cite{azais2018power} proved that the spacing test is unbiased and introduce a studentized version of this test. In the same direction, inference after model selection has been studied in several papers, such as~\cite{fithian2014optimal,taylor2015statistical} and respectively~\cite{tian2018selective} for {\it selective inference} and respectively a joint estimate of the noise level. 

In the present paper, our test is based on the conditional joint law of three, {\bf not necessarily consecutive}, knots. In this way, we extend the work from \cite{taylor2014exact} where the spacing test works for consecutive knots. We refer to these new tests as the generalized spacing tests (GSt). Furthermore, the exact formulation of the spacing test of the pioneering work of~\cite{taylor2014exact} requires extra computations of the term denoted by~$M^+$ in \cite[Lemma 5]{taylor2014exact}. They proved that the spacing test is asymptotically equivalent to the conservative spacing test. We remove this restriction and we prove that it suffices to check wether the so-called {\it Irrepresentable Check} Condition holds to get a non-asymptotic equivalence between the Spacing test and the conservative Spacing test. Finally, we theoretically prove that working with non-consecutive knots can render the testing procedure more powerful.


\subsection{Joint law of LARS knots in Post-Selection Inference}
\label{sec:Q1}
In this paper, we consider linear models in high-dimensions where the number of observations~$n$ may be less than the number of predictors~$p$.
We denote by~$Y\in\mathbb R^n$ the response variable and we assume that
\eq
\label{eq:LinModel}
Y=X\beta^0+\eta\sim\mathcal N_n(X\beta^0,\sigma^2 \mathrm{Id}_n),
\qe
where~$\eta\sim\mathcal N_n(0,\sigma^2 \mathrm{Id}_n)$ is a Gaussian noise, the noise level~$\sigma>0$ may be known or may have to be estimated (depending on the context), and~$X\in\mathds R^{n\times p}$ has rank~$r>0$. We consider the LARS and denote by~$(\lambda_k)_{k\geq1}$ the sequence of knots and by~$(\bar \imath_k,\varepsilon_k)_{k\geq1}$ the sequence~of variables~$\bar \imath_k\in[p]$ and signs~$\varepsilon_k\in\{\pm1\}$ that enter the model along the~LARS path.
We encode by 
\[
\widehat\imath_k:=\bar \imath_k+p\,\Big(\frac{1-\varepsilon_k}2\Big)\in[2p],
\]
both  the variables~$\bar \imath_k\in[p]$ and the signs~$\varepsilon_k\in\{\pm1\}$, calling them the `signed variables'\!. Section~\ref{sec:LAR3} recalls LARS (Algorithm~\ref{alg:LAR}) and present equivalent formulations in Algorithm~\ref{alg:LAR2} (using orthogonal projections) and Algorithm~\ref{alg:LAR3} (using a recursion). In particular, Algorithm~\ref{alg:LAR3} consists in three lines, applying the same function recursively, see Section~\ref{sec:larrec}. As far as we know, {\bf Algorithm~\ref{alg:LAR3} is new}. 

For a short moment, consider the simplest linear model, where one observes the target vector $\beta^0=(\beta^0_1,\ldots,\beta^0_p)\in\mathds R^p$, namely there is no noise and the design~$X=\mathrm{Id}_p$ is the identity.
In this case, LASSO and LARS give the same knots~$\lambda_1,\lambda_2,\ldots$  and the estimate of the LASSO is the outcome of the proximal operator of the~$\ell_1$-norm at point~$\beta^0\in\mathds R^p$, see for instance~\cite[Chapter 2]{hastie2015statistical}.
In this simple case, we deduce that the knots are 
\eq
\label{eq_simplest}
\lambda_k=\beta^0_{(k)},
\qe
where we have considered the reordering~$\beta^0_{(1)}\geq\beta^0_{(2)}\geq\ldots$ of the entries of the target.
Obviously, this is no longer true for general designs in high-dimensions with noise, but one may ask:
\begin{itemize}
\item[\bf{[Q1]}]
{\it 
What is the {\bf joint law} of the LARS knots~$\lambda_1,\lambda_2,\ldots$ and how do they relate to the target~$\beta^0$?}
\end{itemize}

{We will answer~$\bf{[Q1]}$ in high-dimensions under the assumption of Gaussian noise in Section~\ref{sec:selection_event} and Section~\ref{subsec:main}.
Working with the so-called `Irrepresentable Check' Condition\footnote{See Section~\ref{sec:Assumption} for a definition and detailed comments on this assumption.}~\eqref{hyp:IrrAlongThePath}, which can be efficiently checked in practice, we are able to provide the joint law of the LAR's knots conditional on the so-called `selection event' defined by
 \[
 \mathcal E:=\big\{\hat \imath_1 =\imath_1, \ldots, \hat \imath_{K}=\imath_{K},\lambda_{K+1}\big\}.
 \]
 This selection event states that the signed variable~$\imath_k$ has been selected by the LARS algorithm at its~$k^{\text{th}}$ step for~$k=1,\ldots,K$.
This is the cornerstone of the paper, showing that the conditional joint distribution of the LARS knots is a mixture of Gaussian order statistics, as presented in the next theorem.
 
 \pagebreak 
 
 \begin{theorem}[Conditional Joint Law of the~LARS knots]
\label{thm:Main1}
Let~$(\lambda_1,\ldots, \lambda_K,\lambda_{K+1})$ be the first knots of the LARS and let~$(\,\hat \imath_1,\ldots, \hat \imath_K)$ be the first variables entering along the~LARS path.
If~$(\,\hat \imath_1,\ldots, \hat \imath_K)$ satisfies~\eqref{hyp:IrrAlongThePath}, then, conditional on the selection event~$\{\,\hat \imath_1,\ldots, \hat \imath_K,\lambda_{K+1}\}$, the vector~$(\lambda_{1},\ldots,\lambda_K)$ obeys a law with the following density~(w.r.t. Lebesgue measure) 
\[
{\mathrm Q}^{-1}_{(\,\hat \imath_1,\ldots,\hat \imath_K,\lambda_{K+1})}
\,
\Big(
\prod_{k=1}^K \varphi_{m_k,v_k^2}(\ell_k)
\Big)
\,
\mathds{1}_{\{\ell_1\geq \ell_{2}\geq\cdots\geq\ell_{K}\geq\lambda_{K+1}\}}\text{ at point }(\ell_1, \ell_{2},\ldots,\ell_{K}),
\]
where~${\mathrm Q}_{(\,\hat \imath_1,\ldots,\hat \imath_K,\lambda_{K+1})}$ is a normalizing constant,~$\varphi_{m_k,v_k^2}$ is the standard Gaussian density with mean~$m_k$ and variance~$v_k^2:=\sigma^2\rho_k^2$, are explicitly given by~\eqref{eq:def_mu_proj} and~\eqref{eq:def_rho}.
\end{theorem}

\noindent The proof of this theorem is given in Section~\ref{proof:Main1}.
Now, let us describe the dependency between~$(m_k,\rho_k^2)$ and~$\bar\mu^0:=X^\top X\beta^0$.
For a design matrix~$X$ with columns $(X_{j})_{j=1}^p$, we denote by\footnote{Recall that the selected variables~$\widehat\imath_k\in[2p]$ are decomposed into~$\widehat\imath_k:=\bar \imath_k+p\,(\frac{1-\varepsilon_k}2)$.} 
\eq
\notag
\{0\}=:H_0\subset
H_1\subset \cdots \subset H_k:=\mathrm{Span}(X_{\bar\imath_1},\ldots,X_{ \bar\imath_k})\subset \cdots \subset H_K.
\qe
By~\eqref{eq:def_mu_proj} and~\eqref{eq:def_rho}, one has 
\[
\forall k\in[K],\quad 
 m_k = c_k \, \varepsilon_k\, \langle X_{\bar\imath_k},P_{k-1}^\perp(X\beta^0)\rangle\quad \text{and}\quad \rho_k^2 =d_k\, \sin \measuredangle (X_{\bar\imath_k}, H_{k-1}),
\]
where~$c_k,d_k>0$ are constants that depends only on~$X_{\bar\imath_1},\ldots,X_{ \bar\imath_{k-1}}$,~$P_{k-1}^\perp$ denotes the orthogonal projection onto the orthogonal of~$H_{k-1}$,~$\varepsilon_k$ is the sign of the~$k^{\text{th}}$ variable entering the LARS path, and~$\measuredangle (X_{\bar\imath_k}, H_{k-1})$ is the angle between~$X_{\bar\imath_k}$ and~$H_{k-1}$.

\subsection{The Generalized $t$-Spacing test (GtSt)}

This paper introduces a class of exact tests built from~$\ell_1$-minimization regression in high-dimensions. {More precisely, we design a testing procedure for a null hypothesis of the form 
\eq
\notag
\mathds H_0\,:\ \text{`}X\beta^0\in H_{a_0}\text{'},
\qe
where~$H_{a_0}:=\mathrm{Span}(X_{\bar\imath_1},\ldots,X_{ \bar\imath_{a_0}})$.
Note that the null~$\mathds H_0$ is equivalent to the hypothesis that all the true positives (i.e., the support of~$\beta^0$) are among the first~$a_0$ variables selected by LARS, namely~$\{\bar\imath_1,\ldots,\bar\imath_{a_0}\}$.} Following the original idea of~\cite{lockhart2014significance}, we study testing procedures of~$\mathds H_0$ based on the knots of the LARS path.
Note that, conditional on the selection event, the law of three consecutive knots has been studied by~\cite{lockhart2014significance}, where it was referred to as the {\it spacing test} (ST) \citep{hastie2015statistical}.
The article~\cite{azais2018power} proved that the spacing test is unbiased, and introduced a Studentized version of this test.
In the same direction, inference after model selection has been studied in several papers, such as~\cite{fithian2014optimal,taylor2015statistical} and respectively~\cite{tian2018selective} for selective inference and respectively a joint estimate of the noise level.
This raises the following questions.
\begin{itemize}
\item[\bf{[Q2]}]
{\it 
Can we provide exact testing procedures based on knots that are not consecutive?}
\item[\bf{[Q3]}]
{\it 
What is the most powerful test among these spacing tests?}
\item[\bf{[Q4]}]
{\it 
Can we provide exact testing procedures when the noise level is not known?}
\end{itemize}

\noindent $\bullet\ \mathbf{{Contribution\ }(i)}:$ First, our test is based on the conditional joint law of three, not necessarily consecutive, knots~$a_0 \leq a < b < c \leq K + 1$.
In this way, we extend the work of \cite{taylor2014exact}, where the spacing test works for consecutive knots.
We present this framework in Section~\ref{subsec:main} and we refer to these new tests are Generalized Spacing tests (GSt). At the first reading of the next theorem, one can set~$\hat m = a_0$ for a fixed value~$0\leq a_0\leq a$. The selection procedure, defining $\hat m$, will be presented in Section~\ref{sec:exact_test_admissible} with the notion of an `admissible procedure'~{\eqref{e:stopping_rule}.

\begin{theorem}
\label{nt:pvalue}
Let~$a, b,$ and $c$ be such that~$0\leq a<b<c\leq K+1$.
Let~$(\lambda_1,\ldots, \lambda_K,\lambda_{K+1})$ be the first knots and let~$(\,\hat \imath_1,\ldots, \hat \imath_K)$ be the first variables entering along the~LARS path.
If~$(\,\hat \imath_1,\ldots, \hat \imath_K)$ satisfies~\eqref{hyp:IrrAlongThePath} and~$\hat m$ is chosen according to a procedure satisfying~{\eqref{e:stopping_rule}}, then under the null hypothesis
\[
\mathds H_0\,:\ \text{‘‘}X\beta^0\in H_{a}\,\text{''},
\] 
and conditional on the selection event~$\big\{\,\hat m\leq a\big\}$, it follows that
\eq
\label{zaza5}
\hat\alpha_{abc}=\hat\alpha_{abc}(\lambda_a,\lambda_b,\lambda_c, \hat \imath_{1}, \ldots,\hat \imath_{c-1}):=1-\frac{\bbF_{abc}(\lambda_b)}{\bbF_{abc}(\lambda_a)}\sim\cU(0,1),
\qe
namely, it is uniformly distributed over~$(0,1)$.
\end{theorem}
\noindent
The proof of Theorem~\ref{nt:pvalue} is presented in Section~\ref{proof:nt:pvalue}. The construction of the~$p$-values~$\hat\alpha_{abc}$ and of~$\bbF_{abc}$ is given in~\eqref{fun:fabc} and Section~\ref{sec:fabc} respectively. We consider the following Generalized Spacing test procedures (GSt):
\eq
\label{zaza6}
\cS_{abc}:=\mathds1_{\{\hat\alpha_{abc}\leq\alpha\}},
\qe
that rejects if the~$p$-value~$\hat\alpha_{abc}$ is less than the level~$\alpha$ of the test. One may remark~that 
\begin{center}
{\bf the~$p$-value~$\hat\alpha_{abc}$ detects abnormally large values of~$\lambda_b$ relatively to the interval~$(\lambda_a,\lambda_c)$}.
\end{center}
\noindent
When the noise variance is unknown, we introduce the Generalized~$t$-Spacing tests (GtSt) whose theoretical guarantees are given in the next theorem. The estimator of the variance $\hat\sigma^2$ is given in Section~\ref{sec:var_estimation}.
\begin{theorem}
\label{t:pvalue}
Let~$a, b,$ and $c$ be such that~$0\leq a<b<c\leq K+1$.
Let~$(\lambda_1,\ldots, \lambda_{K+1})$ be the first knots and let~$(\,\hat \imath_1,\ldots, \hat \imath_K)$ be the first variables entering along the~LARS path.
If~$(\,\hat \imath_1,\ldots, \hat \imath_K)$ satisfies~\eqref{hyp:IrrAlongThePath} and~$\hat m$ is chosen according to a procedure satisfying~{\eqref{e:stopping_rule}} then under the null hypothesis
\[
\mathds H_0\,:\ \text{`}X\beta^0\in H_{a}\,\text{'},
\] 
and conditional on the selection event~$\big\{\hat m \leq a\big\}$, it follows that
\[
 \hat  \beta_{abc}= \hat  \beta_{abc}( {\Lambda_a},{ \Lambda_b},{ \Lambda_c}, \hat \imath_1, \ldots, \hat \imath_{K}):= 1 - 
  \frac{\tilde\bbF_{abc}(\Lambda_b)}
 {\tilde\bbF_{abc}(\Lambda_a)}\sim \cU(0,1),
\]
where $\Lambda_k:=\lambda_k/\hat\sigma$.
\end{theorem}

\noindent
The proof of Theorem~\ref{t:pvalue} is presented in Section~\ref{proof:t:pvalue}.
The construction of the~$p$-values~$\hat\beta_{abc}$, of~$\tilde\bbF_{abc}$ and of the estimation of the noise~$\hat\sigma$ is given in~\eqref{e:talpha}, Section~\ref{fun:tilde_fabc} and Section~\ref{sec:var_estimation} respectively. One may remark that  
\begin{center}
{\bf 
the~$p$-value~$\hat\beta_{abc}$ detects abnormally large values of~$\Lambda_b$ relatively to the interval~$(\Lambda_a,\Lambda_c)$.}
\end{center}

\medskip

{
\noindent
$\bullet\ \mathbf{{Contribution\ }(ii)}:$ Working with three consecutive knots, we recover the spacing test of~\cite{taylor2014exact} and even in this framework, the present  paper improves the current state of knowledge.
We prove that:\smallskip
\newline
\noindent
{\it~$\circ$ Under~\eqref{hyp:IrrAlongThePath}, the Spacing test procedure defined in~\cite[Theorem 1]{taylor2014exact} is exact, and is equal to the so-called `conservative' spacing test defined in \cite[Theorem~2]{taylor2014exact}.}
\smallskip
\newline
The exact formulation of the spacing test of the pioneering work of~\cite{taylor2014exact} requires extra computations of the term denoted by~$M^+$ in \cite[Lemma~5]{taylor2014exact}.
They proved that the spacing test is asymptotically equivalent to the conservative spacing test.
We remove this restriction and we prove that it suffices to check wether the Irrepresentable Check Condition~\eqref{hyp:IrrAlongThePath} holds to get a non-asymptotic equivalence between the spacing test and the conservative spacing test.

\medskip

{
\noindent
$\bullet\ \mathbf{{Contribution\ }(iii)}:$  We theoretically prove that working with non-consecutive knots can allow obtaining higher power for the testing procedure.

 \begin{theorem}
 \label{thm:zazPower}
Assume that the design~$X$ is orthogonal, namely~$X^\top X=\mathrm{Id}_{p}$.
Let~$a_0$ be an integer such that~$0\leq a_0\leq K-1$.
If~$\hat m$ is chosen according to a procedure satisfying~{\eqref{e:stopping_rule}}, then under the null hypothesis
\[
\mathds H_0\,:\ \text{`}X\beta^0\in H_{a_0}\,\text{'},
\] 
and conditional on the selection event~$\big\{\hat m=a_0\big\}$, it follows that the test~$\cS_{a_0,a_0+1,K+1}$ is uniformly more powerful than any of the tests~$ \cS_{a,b,c}$ for~$a_0\leq a <b<c \leq K+1$.
 \end{theorem}

 \noindent
The proof of this result is given in Appendix \ref{proof:zazPower}.
 This shows that the most powerful test among the set of tests~$(\cS_{a,b,c})_{a_0\leq a <b<c \leq K+1}$ is given by
\begin{center}
{\bf  the GSt test~$\cS_{a_0,a_0+1,K+1}$ with the smallest~$a$ and the largest~$c$}.
\end{center}
More precisely, in the proof of Theorem~\ref{thm:zazPower}, it is shown that 
\[
\hat \alpha_{ab(c+1)}\preccurlyeq \hat \alpha_{abc}\text{ and }
\hat \alpha_{a(b-1)c}\preccurlyeq \hat \alpha_{abc}\text{ and }
\hat \alpha_{(a-1)bc}\preccurlyeq \hat \alpha_{abc},
\]
for orthogonal designs, where~$\preccurlyeq$ denotes stochastic ordering.}

\subsection{A new exact testing procedure on false negatives (FN) after support selection}

{
\begin{algorithm}[!t]
\setcounter{AlgoLine}{0}
\captionsetup{labelfont={sc,bf}, labelsep=newline}
\caption{Exact false negative testing after model selection}
  \label{Alg:TestingAStop}
  \SetAlgoLined 
  \KwData{$K$ satisfying~\eqref{e:K}, selection procedure~$\hat m$ satisfying~\eqref{e:stopping_rule}, couple~$(X,Y)$ giving design and response.}
          \BlankLine  
  \KwResult{$p$-value~$\hat \alpha$ on the existence of false negative.}
  {\nonl \tcc{$\1_{\{\hat \alpha\leq\alpha\}}$ is a testing procedure with level exactly~$\alpha$}}
 \BlankLine
  {Compute the LARS path from~$(X,Y)$.}
\BlankLine
  {Check that~$(\,\hat \imath_1,\ldots, \hat \imath_K)$ satisfies~\eqref{hyp:IrrAlongThePath}.
If not {\bf Stop}.}
\BlankLine
 {Compute~$\hat m$, the size of the selected model}.
 \BlankLine      
 {{\bf Return}~$\hat\alpha=\hat\alpha_{\hat m(\hat m+1)(K+1)}$, see~\eqref{zaza5}.\\}
 {\nonl \tcc{When variance is unknown,~$\hat\alpha=\hat\beta_{\hat m\, (\hat m+1)\,(K+1)}$, see~\eqref{e:talpha}.}}
\end{algorithm}
}

One specific task is to estimate the support of the target sparse vector, namely identify the true positives in the context of a multiple testing procedure.
In particular, one may take the support of the LASSO (or SLOPE) solution as an estimate of the support of the solution.
This strategy has been intensively studied in the literature, one may consider~\cite{wainwright2009sharp,bogdan2015slope,van2016estimation,bellec2018slope} and references therein.
Support selection has been studied under the so-called `Irrepresentable Condition' (IC), as presented for instance in~\cite[Page 53]{van2016estimation} and~\cite[Sec.~7.5.1]{Buhlmann_vandeGeer11} and also referred to as the `Mutual Incoherence Condition'~\citep{wainwright2009sharp}.
Under the so-called `Beta-Min Condition', one may prove~\citep{Buhlmann_vandeGeer11,van2016estimation} that the LASSO asymptotically returns the true support.
Following this line of thought, one may ask:

\begin{itemize}
\item[\bf{[Q5]}]
{\it 
Can we provide a false negative testing procedure with a controlled Type I error?}
\end{itemize}

In this article, we build an exact non-asymptotic multiple test for false non-negatives, see Sections~\ref{sec:exact_model_selection} and \ref{sec:FalseNeg}.
The control of the false negatives after model selection in the case  of an unknown noise level is given in Section~\ref{sec:FalseNeg} and the procedure is introduced in Algorithm~\ref{Alg:TestingAStop}.
We assume~\eqref{e:stopping_rule}, which assumes that the model has been selected using an `admissible' procedure, which basically means that the decision to select a model of size~$a$ only depends on the orthogonal projection of the observation~$Y$ onto~$H_a=\mathrm{Span}(X_{\bar\imath_1},\ldots,X_{\bar\imath_a})$.
Assuming further that~\eqref{hyp:IrrAlongThePath} holds, we provide an exact testing method for false negatives.
In order to reach high power, the test statistic is the~$p$-value of three non-consecutive knots of the LARS path.
To compute this~$p$-value, one needs to marginalize the joint law of the knots, leading to a numerical integration whose complexity grows exponentially with the space between the indices of the knots.
We propose to use QMC techniques to compute the statistic, see Appendix~\ref{app:cbc}.

\subsection{False Discovery Rate control for LARS}

Simultaneous controls of confidence intervals independently of the selection procedure have been studied under the concept of {\it post-selection constants} as introduced in~\cite{berk2013valid} and studied for instance in~\cite{bachoc2018post}.
Asymptotic confidence intervals can be build using the {\it de-sparsified LASSO}, the reader may refer to~\cite[Chapter 5]{van2016estimation} and references therein.
We also point a recent study~\citep{javanmard2019false} of the FDR control as the sample size tends to infinity using {\it de-biased LASSO}, which has been implemented in Section~\ref{sec:simulated}.
Asymptotic FDR control has been studied in~\cite{barber2015controlling} and references therein, which has been implemented in Section~\ref{sec:simulated}.
Let us point recent control of the {\it Joint family-wise Error Rate} as in~\cite{blanchard2017post} and references therein.
Following these lines of work, one may ask:
\begin{itemize}
\item[\bf{[Q6]}]
{\it 
Can we provide {multiple} Spacing Tests with a controlled False Discovery Rate (FDR)?}
\end{itemize}
To the best of our knowledge, this paper is the first to study the joint law and an exact control of multiple spacing tests of LARS knots in a non-asymptotic frame, see Sections~\ref{subsec:main} and~\ref{sec:FDR}.
We investigate the consecutive spacings of the knots of the LARS as test statistics and we prove an exact FDR control using a Benjamini–Hochberg procedure~\citep{benjamini1995controlling} in the orthogonal design case, see Theorem~\ref{thm:independentTest} and Section~\ref{sec:FDR}.
Our proof (see Appendix~\ref{proof:yoFDR}) is based on the {\it Weak Positive Regression Dependency} (WPRDS), the reader may consult~\cite{blanchard2008two} or the survey~\cite{roquain2010type}, and {\it Knothe-Rosemblatt transport}, see for instance~\cite[Sec.2.3, Page 67]{santambrogio2015optimal} or~\cite[Page 20]{villani2008optimal}, which is based on conditional quantile transforms.

\subsection{Additional related works on high-dimensional statistics}
Parsimonious models have become ubiquitous tools to tackle high-dimensional representations with a small budget of observations.
Successful applications may be found in signal processing (see for instance the pioneering works of~\cite{Chen_Donoho_Saunders98,Candes_Romberg_Tao06} and references therein) and biology (see for instance~\cite{barber2015controlling} or~\cite[Chapter~1.4]{Buhlmann_vandeGeer11} and references therein).
These applications have shown that there are interesting {\it almost sparse representations} in some well chosen basis.
Nowadays, in many practical situations, this sparsity assumption is recognized as reasonable.

These important successes have put a focus on High-Dimensional Statistics and Compressed Sensing in the past decades, which may be due to the deployment of tractable algorithms with strong theoretical guarantees.
 Among the large panoply of methods, one may consider~$\ell_1$-regularization, which benefits from a remarkable tractability, empirical performance, and theoretical guarantees.
Nowadays, sparse regression techniques based on~$\ell_1$-regularization are a common and powerful tool in high-dimensional settings.
Popular estimators, among which one may point to the LASSO \citep{tibshirani1996regression} and SLOPE~\citep{bogdan2015slope}, are known to achieve a minimax rate of prediction and to satisfy the sharp oracle inequalities under conditions on the design, such as Restricted Eigenvalue~\citep{bickel2009simultaneous,bellec2018slope} or Compatibility~\citep{Buhlmann_vandeGeer11,van2016estimation}.
The sharp oracle inequalities show that the estimation errors, in~$\ell_1$ and~$\ell_2$ norm, of these estimators are optimal, see for instance \cite[Chapter 2.7]{van2016estimation}.

Variable selection has also been investigated, and it has been proven, see for instance \cite[Theorem~7.5]{Buhlmann_vandeGeer11}, that the LASSO selects the true variables ({\it i.e.,} there are no false negatives) under the Compatibility condition and the so-called beta-min condition (which assumes that the true parameters are large enough with respect to some threshold that scales linearly with the regularization parameter $\lambda$ of LASSO).
Under a stronger assumption, referred to as the `irrepresentable condition', one can prove, see for instance \cite[Theorem~7.1]{Buhlmann_vandeGeer11}, that the~$\ell_\infty$-estimation error scales linearly with the regularization parameter $\lambda$ of LASSO.
As the regularization parameter~$\lambda$ tends to zero, when the number of observation goes to infinity and under some assumption on the noise, these results show that LASSO produces a consistent selection of the variables (it asymptotically finds the true support with no errors).   }

\subsection{Outline of the paper}  
\subsubsection{Detailed outline}  

Section~\ref{sec2} introduces the notation (see also Section~\ref{table:mr} for a summary), assumptions~\eqref{hyp:IrrAlongThePath} and \eqref{e:stopping_rule}, and variance estimate $\hat \sigma^2$.
The variance estimate is a key step in our testing procedures: we introduce new variance estimate with properties useful for deriving exact and non-asymptotic post-selection laws, see Section~\ref{sec:var_estimation}.

The main assumption is based on the Irrepresentable Check condition~\eqref{hyp:IrrAlongThePath}, which can be checked in practice, see Section~\ref{sec:Assumption}.
Under~\eqref{hyp:IrrAlongThePath}, we obtain a new characterization of the selection event in Proposition~\ref{prop:reccur} of Section~\ref{sec:frozen_knots}.

Section~\ref{sec3} gives the main results:  Section \ref{subsec:main} describes the joint distribution of the LARS knots as a mixture  of Gaussian order statistics and the GST and GtST tests.
The power in the orthogonal case is considered in Section \ref{subsec:ortho}.
The control of the  false negatives  in a post selection  inference with estimation of the variance  is presented in Section \ref{sec:FalseNeg} (when the variance is known, this procedure is studied in Section~\ref{sec:exact_model_selection}).
A procedure  to control the FDR in the orthogonal case is presented  in Section~\ref{sec:FDR}.

Illustrations of our method, both on simulated data and on real data, are presented in Section~\ref{sec:num}.
A Zenodo repository of the code used in all our experiments can be found at \cite{yohann_de_castro_2021_5079768}.


\pagebreak

\subsubsection{Dependency diagram}  
The outline can be depicted by the following {\it dependency diagram}:

\vspace*{0.5cm}

\hspace*{-1cm}
\begin{tikzpicture}[>=stealth,every node/.style={shape=rectangle,draw,rounded corners},
greennode/.style={rectangle, draw=burgundy!60, fill=burgundy!5, very thick, minimum size=5mm},
rednode/.style={rectangle, draw=camel!60, fill=camel!5, very thick, minimum size=5mm},
]
    \node[greennode] (c20) {Conditional hypotheses $H_a$ in Sec.~\ref{sec:hypotheses}};
    \node[rednode] (c21) [right=2.43cm of c20]{ Proposition~\ref{prop:reccur}, {\bf selection event} (Sec. \ref{sec:frozen_knots})};
    \node[greennode] (c1) [below =0.6cm of c20]{Irrepresentable Check \eqref{hyp:IrrAlongThePath} in Sec. \ref{sec:Assumption}};
    \node[greennode] (c2) [below =0.6cm of c1]{Variance estimate $\hat \sigma^2$ in Sec. \ref{sec:var_estimation}};
    \node[greennode] (c3) [below =0.6cm of c2]{Admissible procedures \eqref{e:stopping_rule} in Sec. \ref{sec:exact_test_admissible}};
    \node[greennode] (c5) [below =1.7cm of c3]{Orthogonal design, $X^\top X=\mathrm{Id}_p$};
        \node[greennode] (c11) [below =0.6cm of c5]{{\bf Assumptions \& tools}};
    \node[rednode] (c6) [right=2.01cm of c1]{Theorem~\ref{thm:Main1}, {\bf conditional joint law} (Sec. \ref{subsec:main})};
    \node[rednode] (c7) [below right =0.6cm and 1.82 cm of c3]{Theorem~\ref{thm:zazPower}, {\bf power/Type II error} (Sec. \ref{subsec:ortho})};
    \node[rednode] (c8) [right=2.79cm of c2]{Theorem~\ref{nt:pvalue}, {\bf GSt \& FN testing} (Sec. \ref{sec:exact_model_selection})};
    \node[rednode] (c9) [right=2.1cm of c3]{Theorem~\ref{t:pvalue}, {\bf GtSt \& FN testing} (Sec. \ref{sec:Student})};
    \node[rednode] (c10) [right=3.15cm of c5]{Theorem~\ref{thm:independentTest}, {\bf FDR control} (Sec. \ref{sec:FDR})};
    \node[rednode] (c12) [right=6cm of c11]{{\bf Results}};

    \draw[->] (c1) to[out=0,in=-180] (c6);
    \draw[->] (c1) to[out=0,in=-180] (c21);
    \draw[->] (c1) to[out=0,in=-180] (c8);
    \draw[->] (c1) to[out=0,in=-180] (c9);
    \draw[->] (c3) to[out=0,in=-180] (c8);
    \draw[->] (c3) to[out=0,in=-180] (c9);
    \draw[->] (c2) to[out=0,in=-180] (c8);
    \draw[->] (c20) to[out=0,in=-180] (c8);
    \draw[->] (c5) to[out=0,in=-180] (c7);
    \draw[->] (c5) to[out=0,in=-180] (c10);
    \draw[->] (c2) to[out=0,in=-180] (c9);
    \draw[->] (c20) to[out=-2,in=-182] (c9);
    \draw[->] (c6) -- (c8);
    \draw[->] (c21) -- (c6);
    \draw[->] (c20) -- (c21);
    \draw[->] (c6) to[out=0,in=3] (c9);

\end{tikzpicture}

\vspace*{0.3cm}

\subsubsection{Notation and commands}  
\label{table:mr}
\vspace*{0.3cm}

\begin{center}
	{\small
	\begin{tabular}{ll}
		\multicolumn{2}{c}{\textit{General notation}} \\
		\midrule
		$[a]$ 
		& the set of integers~$\{1,...,a\}$ \\
		{\ } & {\ }\\
		$ Y = X \beta^0 + \eta~$ 
		&  Linear Model~\eqref{eq:LinModel},~$X$  is~$n\times p$ design matrix with rank~$r$ \\
		$\sigma^2$ 
		& the variance of the errors~$\eta$ \\
		{\ } & {\ }\\
		$K$ 
		&  the number of knots~$\lambda_1,\ldots,\lambda_K$ that are considered, see~\eqref{e:K} \\
		$n_1$,$n_2$ 
		& number of d.o.f.
used for constructing~$  \hat \sigma$ \\
		{\ } & {\ }\\
		$\varphi_{m_k, v^2_k}~$  &standard Gaussian density with mean~$m_k$ and variance~$v_k^2:=\sigma^2\rho_k^2$\\ 
		$\tilde\varphi$
		& multivariate~$t$-distribution with~$\nu=n_2$ degrees of freedom, mean~$m=(m_1,\ldots,m_K)$\\
		& and variance-covariance matrix~$\mathrm{Diag}(\rho_1,\ldots,\rho_K)$\\
		$m_k,v_k^2$ 
		& conditional mean, see~\eqref{eq:def_mu_proj}, and conditional variance~$v^2_k=\sigma^2\rho_k^2$, see~\eqref{eq:def_rho}\\
		{\ } & {\ }\\
		$ \hat \alpha_{abc}$ 
		& the~$p$-value  of the generalized spacing test (GSt), see~\eqref{zaza5}\\
		$ \cS_{abc}$ 
		& $\mathds1_{\{\hat\alpha_{abc}\leq\alpha\}}$, the generalized spacing test (GSt) see~\eqref{zaza6} \\
		{\ } & {\ }\\
		$\Lambda_k$
		& $t$-knots defined by~\eqref{def_Lambda}\\
		$\hat\beta_{abc}$
		& the~$p$-value of the generalized~$t$-spacing test (GtSt), see~\eqref{e:talpha}\\
		$ \cT_{abc}$ 
		& $\mathds1_{\{\hat\beta_{abc}\leq\alpha\}}$, the {\it generalized~$t$-spacing test} (GtSt), see~\eqref{zaza66} 
	\end{tabular}}
\end{center}

\begin{center}
	{\small
	\begin{tabular}{ll}
		\multicolumn{2}{c}{Technical notation} \\
		\midrule
		$ \hat \imath_k~$
		&  a way of coding both indices and signs, see~\eqref{eq:signed_indices} \\
		$ \bar \imath_k ; \varepsilon_k$
		& the indices and the signs of  the variables that enter in the LARS path \\
		$j_1,\ldots,j_k$;~$s_1,\ldots,s_k~$ 
		& a generic  value of the sequences above \\
		$i_1,\ldots,i_k$ 
		& a generic  value of the sequence~$ \hat \imath_k~$\\
		{\ } & {\ }\\
		$Z$ 
		& the vector of correlations, obtained by symmetry from~$\bar{Z}$ defined by~\eqref{eq:defBarZ}\\
		$R$ 
		& the variance-covariance matrix of~$Z$, see~\eqref{eq:white} \\
		$M_{i_1,\ldots, i_{\ell}}$
		&  sub-matrix of~$R$ indexed by~$\{i_1,\ldots, i_{\ell}\}$, see~\eqref{e:M}\\
		{\ } & {\ }\\
		$  \bar S^k$ 
		& $\{ \bar \imath_1,\ldots, \bar \imath_k$\}, a possible selected support~\eqref{eq:sequence_supports} \\
		$S_0$
		& the true support\\
		$\hat {S}$
		& the chosen set of variables :~$\bar S ^{\hat m}$, see \eqref{eq:size_model}\\
		$\hat m$  
		& the chosen size \\
		\eqref{e:stopping_rule}
		& stopping rule, see Section \ref{s:ap}\\
		$ H_k$
		&$ \mathrm{Span}(X_{\bar \imath_1},\ldots,X_{\bar \imath_k})$\\
		$P_k (P_k^\perp)$ 
		& Orthogonal projection on (the orthogonal of)~$H_k$\\
		{\ } & {\ }\\
		\eqref{hyp:IrrAlongThePath} 
		&  Irrepresentable Check, see~\eqref{hyp:IrrAlongThePath}\\
		$\theta_j(i_1, \ldots,i_k)~$  
		& expectation of $Z_j$ conditional on $Z_{i_1}=\cdots=Z_{i_k}=1$, see~\eqref{eq:thetai}, and~$\theta^{\ell}:=\theta(\,\hat \imath _1, \ldots, \hat \imath _\ell)$\\
		{\ } & {\ }\\
		$Z^{(i_1,\ldots, i_{k})}_j$
		& frozen residual, see~\eqref{e:frozenZ}\\
		$ \Pi_{ i_1,\ldots, i_{k}} (Z_j)$ 
		&  regression of $Z_j$ on $(Z_{i_1},\ldots,Z_{i_k})$, see~\eqref{e:Pi}\\
		~$\lambda^f_{k} :=   Z^{i_1,\ldots i_{k-1}}_{i_{k}}$ 
		  & the~$k^{\text{th}}$ frozen knot, see~\eqref{e:frozen}\\
		~$m_k^f,\sigma\rho_k^f$ 
		 & mean~\eqref{e:frozen_m} and standard deviation~\eqref{e:frozen_rho} of~$\lambda^f_{k}$\\
		 {\ } & {\ }\\
		$\bbF_{abc}(t)$
		& up to some numerical constant, the CDF of $\lambda_b\ |\  \lambda_a,\lambda_c$, see~\eqref{fun:fabc} \\
		$ F_i;\mathcal P_{ij}$  
		&$F_i:=\Phi_i(\lambda_i) :=  \Phi(\lambda_i/(\sigma\rho_i))$ and~$\mathcal P_{ij}$ is given by~\eqref{eq:pij}\\
		$\tilde\bbF_{abc}(t)$
		& up to some numerical constant, the CDF of $\Lambda_b\ |\  \Lambda_a,\Lambda_c$, see~\eqref{fun:tilde_fabc}\\
		$\bT_k$
		& up to some numerical constant, the CDF of centered~$t$-Student distribution, see~\eqref{e:Student}
	\end{tabular}}
\end{center}

\medskip


\section{Assumptions, Variance Estimation and Admissible Procedures}
\label{sec2}

\subsection{Signed variables of LARS}
\label{sec:signednotation}
We give some notation that will be useful.
We denote by~$(\,\hat \imath_1,\ldots,\hat \imath_k)\in[2p]^k$ the `signed}' variables that enter the model along the~LARS path with the convention that
\eq
\label{eq:signed_indices}
\hat \imath_k:=\bar \imath_k+p\,\Big(\frac{1-\varepsilon_k}2\Big),
\qe
so that~$\hat \imath_k\in[2p]$ is a useful way of encoding both the variable~$\bar \imath_k\in[p]$ and its sign~$\varepsilon_k=\pm1$ as used in Algorithm~\ref{alg:LAR3}.
We denote by~$\bar Z:=X^\top Y$ the {\it correlation} vector such that~$\bar Z_k$ is the scalar product between the~$k^{\text{th}}$ predictor and the response variable, and we denote by~$\sigma^2\bar R$ its variance-covariance matrix.
For the sake of presentation, we may consider the~$2p$-vector
\eq
\label{eq:defBarZ}
Z:=(\bar Z,-\bar Z)=(X^\top Y ,-X^\top Y ),
\qe
whose mean is given by
\eq
\label{e:mu0}
\mu^0:=(\bar R\beta^0,-\bar R\beta^0)=(X^\top X\beta^0,-X^\top X\beta^0)=(\bar\mu^0,-\bar\mu^0),
\qe
and whose variance-covariance matrix is~$\sigma^2 R$ with 
\eq
\label{eq:white}
R
=\left[\begin{array}{cc} \bar R & -\bar R \\-\bar R & \bar R\end{array}\right]
=\left[\begin{array}{cc} X^\top X & -X^\top X \\-X^\top X & X^\top X\end{array}\right].
\qe

\noindent
 We also denote by
\begin{enumerate}
\item[$\circ$]~$\hat \imath_1,\ldots,  \hat \imath_k$,  the first~$k$ signed variables entering the LARS,
\item[$\circ$]~$ i_1,\ldots,i_k$, a generic value of the sequence above,
\item[$\circ$]~$ \bar \imath_1 ,\ldots, \bar \imath_k~$, the first~$k$ variables entering the LARS,
\item[$\circ$]~$j_1, \ldots,j_k$, a generic value of the sequence above,
\item[$\circ$]~$\varepsilon_1, \ldots,  \varepsilon_k$, the first~$k$ signs of the coefficients of the variables entering in the LARS,
\item[$\circ$]~$s_1,\ldots,s_k$, a generic value of the sequence above.
\end{enumerate}
\noindent
 The quantities above are related by~\eqref{eq:signed_indices} and
 \eq
\label{e:findices}
i_k:=j_k+p\,\Big(\frac{1-s_k}2\Big).
\qe

\subsection{Models, conditional hypotheses, and the notation~$K$}
\label{sec:hypotheses}
We are interested in selecting the true support~$S^0$ of~$\beta^0$, where the support is defined by
\[
S^0:=\big\{k\in[p]\ :\ \beta^0_k\neq0\big\}.
\]
To estimate this support, we will consider the models that appear along the LARS path: the selected model~$\hat S$ would be chosen from the family of nested models
\eq
\label{eq:sequence_supports}
\underbrace{\{\bar \imath_1\}}
_{\bar S^{1}}
\subset
\underbrace{\{\bar  \imath_1,\bar  \imath_2\}}
_{\bar S^{2}}
\subset\cdots\subset
\underbrace{\{\bar  \imath_1,\bar  \imath_2,\ldots,\bar \imath_a\}}
_{\bar S^{a}}
\subset\cdots\subset
\underbrace{\{\bar  \imath_1,\bar  \imath_2,\ldots,\bar  \imath_K\}}
_{\bar S^{K}},
\qe
where~$K$ denotes the maximal model size.
We denote by $\hat m$ the size of the selected model~$\hat S$, and then
\eq
\label{eq:size_model}
\hat S=\bar S ^{\,\hat m}.
\qe
Respectively, denote
\eq
\label{eq:sequence_subspaces}
\{0\}=:H_0\subset
H_1\subset \cdots \subset H_a:=\mathrm{Span}(X_{\bar \imath_1},\ldots,X_{\bar \imath_a})\subset \cdots \subset H_K,
\qe
the corresponding family of nested subspaces of~$\mathds R^n$.
Once the model has been selected, we will construct tests based on the~$K+1$ first knots of the LARS.

\begin{remark}
\label{rem:selective_testing}
The testing procedures under consideration are not standard since the~$(H_{a})_{a=1}^{K}$ are random subspaces.
We are interested in the framework of selective testing, namely, testing procedures conditional on the selection event~$\{\hat m=a,\hat \imath_1 =i_1, \ldots, \hat \imath_{K}=i_{K}\}$, for some fixed~$a\in[K-1]$.
Conditional on the event, note that~$H_{a}$ is fixed.
By convention, we may consider the case~$a=0$, that is, testing the global null hypothesis.
\end{remark}

\newpage

\noindent
Throughout this paper, we assume that~
\eq
\label{e:K}
K\text{ is fixed and such that }1\leq K <\min(n,r)\text{ where }r=\rank(X).
\qe
In practice,~$K$ can be considerably much smaller than~$n$.
Our analysis is conditional on~$(\,\hat \imath_1,\ldots,\hat \imath_{K})$ and in this spirit it can be referred to as a `Post-Section' procedure, see {\it e.g.}~\cite{taylor2015statistical,taylor2014exact,hastie2015statistical}.

\subsection{Irrepresentable Check on the Active sets}
\label{sec:Assumption}

We define the set of {\it Active Sets}~$\cA_K$ as all the  sequences~$i_1,\ldots,i_K$ of signed variables such that~$j_1,\ldots,j_K$ are pairwise different, where the~$j$'s are defined by~\eqref{e:findices}, namely 
\[
\cA_K:=\big\{(i_1,\ldots,i_K)\in[2p]^K\ :\ j_1,\ldots,  j_K\text{ are pairwise different}\big\}.
\]
Sometimes it would be useful to consider~$\cA_{K+1}$, the set of active sets of size~$K+1$.
We introduce the notion of `Irrepresentable Check', which is the only assumption on the design and the selected active set in most of our results.
  
   \begin{definition}[Irrepresentable Check]\label{def:EIC}
An active set~$(i_1,\ldots,i_K)\in\cA_K$ is said to satisfy the {\it Irrepresentable Check} condition if 
\eq
\label{hyp:IrrAlongThePath}
\tag{${\cA}_{\mathrm{Irr.}}$}
\forall k\in[K],\ \forall j\notin { T^{k}}:=\{j_1,\ldots, j_k\},\quad X_j^\top X_{T^{k}}\big( X_{T^{k}}^\top  X_{T^{k}}\big)^{-1} s_{k}<1,
\qe
where~$j_k$ and~$s_k$ are defined from~$i_k$ using~\eqref{e:findices}.
By a slight abuse of notation we will  denote by~\eqref{hyp:IrrAlongThePath}  the set of sequences~$(i_1,\ldots,i_K)$ that satisfy this property.
\end{definition}


In our procedures and theoretical results, we will limit our attention to sequences~$\hat \imath_1,\ldots, \hat \imath_K\,$ chosen by LARS that  satisfy~\eqref{hyp:IrrAlongThePath}.
 A particular case  is when the  property is true for all possible active sets.
This is equivalent  to the Irrepresentable Condition that we will now recall.

\begin{definition}[Irrepresentable Condition of order~$K$]
The design matrix~$X$ satisfies the Irrepresentable Condition of order~$K$ if and only if 
\eq
\label{eq:Irr_matrix}
\tag{${\mathrm{Irrep.}}$}
\forall S\subset[p]\ \text{s.t.}\ \# S\leq K,\quad
\max_{j\in[p]\setminus S}\max_{||v||_\infty\leq1} X_j^\top X_S\big(X_S^\top X_S\big)^{-1}v<1,
\qe
where~$X_j$ denotes the~$j^{\text{th}}$ column of~$X$ and~$X_S$ the sub-matrix of~$X$ obtained by keeping the columns indexed by~$S$.
\end{definition}
 
\begin{remark}
Note that the Irrepresentable Condition is a standard condition, as presented for instance, in~\cite[Page 53]{van2016estimation} and~\cite[Sec.~7.5.1]{Buhlmann_vandeGeer11}, and is also referred to as the {\it Mutual Incoherence Condition}~\citep{wainwright2009sharp}.
\end{remark}

\begin{remark}
This condition has been intensively studied in the literature and it is now well established that some random matrix models satisfy it with high probability.
For instance, one may refer to~\cite{wainwright2009sharp}, where it is shown that a design matrix~$X\in\bbR^{n\times p}$ whose rows are drawn independently with respect to a centered Gaussian distribution with variance-covariance matrix satisfying~\eqref{eq:Irr_matrix} (for instance the identity matrix) satisfies~\eqref{eq:Irr_matrix} with high probability when~$n\gtrsim K\log(p-K)$, where~$ \gtrsim$ denotes an inequality up to some multiplicative constant.
\end{remark}

In practice, the Irrepresentable Condition~\eqref{eq:Irr_matrix} is a strong requirement on the design~$X$ and, in addition, this condition cannot be checked in polynomial time.
One important feature of our results is that we do not require the Irrepresentable Condition~\eqref{eq:Irr_matrix} but only the weaker requirement of Irrepresentable Check~\eqref{hyp:IrrAlongThePath} on the selected active set.
Namely, we will assume that
\[
\tag{Assumption}
\text{For } K \text{ defined by }\eqref{e:K},\text{ }
(\,\hat \imath_1,\ldots,  \hat \imath_K)\text{ satifies }\eqref{hyp:IrrAlongThePath}.
\]
Given~$(\,\hat \imath_1,\ldots,  \hat \imath_K)$, note that this condition can be checked in polynomial time.

\begin{example}
Taking the (signed) variables entering the model with LARS in an iid Gaussian design and as response variable a centered Gaussian vector with iid entries from~$10,000$ Monte Carlo repetitions, Figure~\ref{fig:F1} illustrates the law of the maximal order~$K_{\mathrm{max}}$ for which the Irrepresentable Check condition holds.
For example, we found that for~$p=1,000$ and~$n=100$~$($with ratio~$n/p=0.1)$ resp.~$n=500$~$($with ratio~$n/p=0.5)$, the Irrepresentable Check condition~\eqref{hyp:IrrAlongThePath} of order~$K_{\mathrm{max}}$ holds when~$K_{\mathrm{max}}$ is about~$K_{\mathrm{max}}\simeq 0.16\times n=16$ respectively~$K_{\mathrm{max}}\simeq 0.12\times n=60$, see Figure~\ref{fig:F1}.
\end{example}

\begin{figure}[!h]
\includegraphics[width=0.39\textwidth]{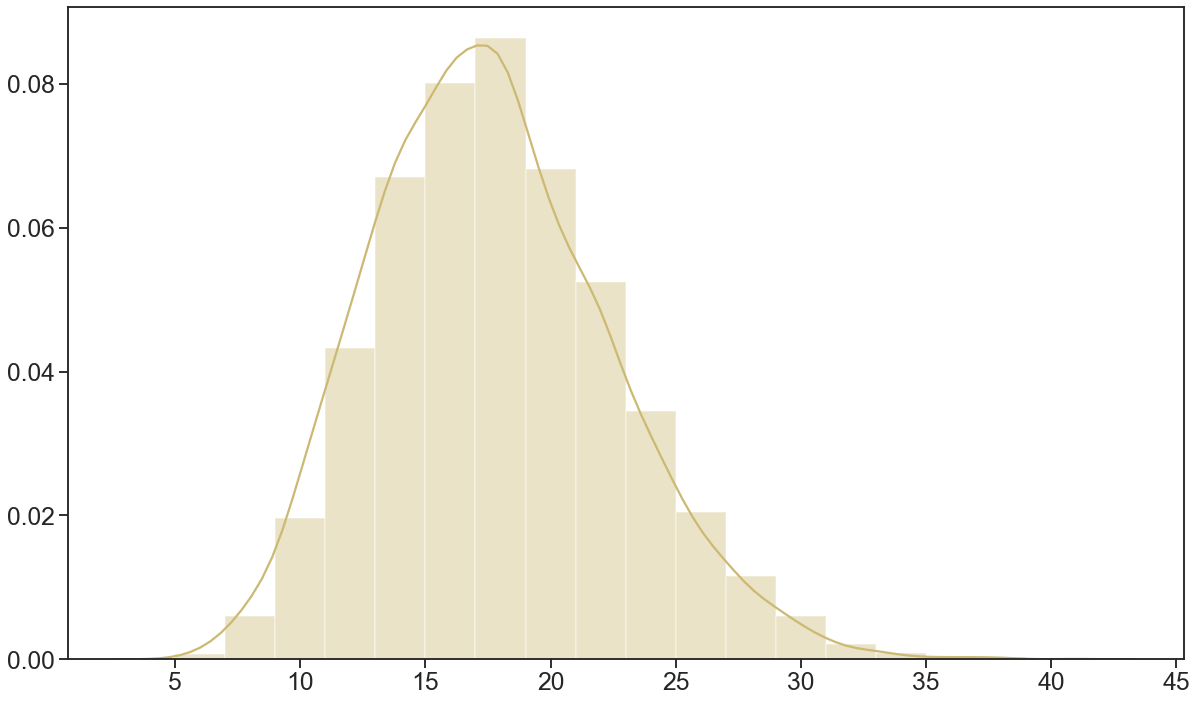}\quad 
\includegraphics[width=0.39\textwidth]{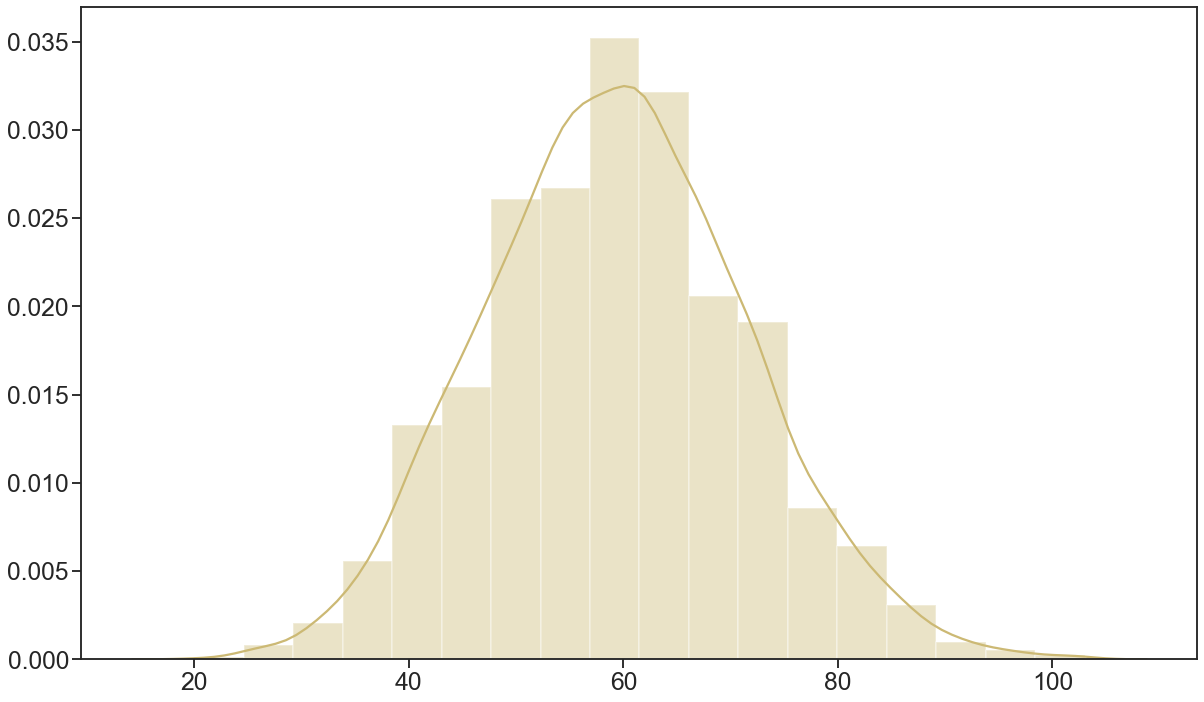}
\caption{The law of the maximal order~$K$ for which Irrepresentable Check holds, taking the (signed) variables entering the model with LARS from an iid Gaussian design and with response variable a centered Gaussian vector with iid entries using~$10,000$ Monte Carlo repetitions. 
There are~$p=1,000$ predictors and~$n=100$ (left)~$n=500$ (right) observations, and we observe that~$K_{\mathrm{max}}\in[10,27]$ (left) and~$K_{\mathrm{max}}\in[39, 81]$ (right) for 95\% of the values.}
\label{fig:F1}
\end{figure}

\subsubsection{Irrepresentable Check: An equivalent formulation}
\label{sec:AssumptionBis}
Now, we can define 
\eq
\label{eq:thetai}
\forall (i_1,\ldots, i_k)\in[2p]^k,\quad
\theta_j(i_1,\ldots, i_k) := \big(R_{j,i_1} \cdots R_{j,i_{k}}  \big)M^{-1}_{i_1,\ldots, i_{k}}  
(1, \ldots, 1),
\qe
where~$(1, \ldots, 1)$ is the column vector of size~$k$ whose entries are equal to one;~$\sigma^2M_{i_1,\ldots, i_{k}}$ is the variance-covariance matrix of the vector~$ (Z_{i_1}, \cdots , Z_{i_{k}})$
and~$(R_{j,i_1} \cdots R_{j,i_{k}})$ is a row vector of size~$k$.
Note that~$M_{i_1,\ldots, i_{k}}$ is the submatrix of~$R$ obtained by keeping the columns and the rows indexed by~$\{i_1,\ldots, i_{k}\}$, namely 
\eq
\label{e:M}
M_{i_1,\ldots, i_{k}}:=(R_{i,j})_{i,j=i_1,\ldots, i_{k}}.
\qe
 Remark that
\[
\theta_j(i_1,\ldots, i_k)=\bbE\big[Z_j\ |\ Z_{i_1}=1,\ldots, Z_{i_k}=1\big],
\]
when~$\bbE Z=0$.
Then Proposition~\ref{prop:Irrepresentable} shows that the Irrepresentable Condition~\eqref{eq:Irr_matrix} of order~$K$ is equivalently given by  
\eq
\label{eq:Irrepresentable}
\forall k\leq K,\ \forall (i_1,\ldots, i_k)\in[2p]^k,\ \forall j\notin\{i_1,\ldots, i_k\},\quad\theta_j(i_1,\ldots, i_k) <1,
\qe
where~$\theta_j(i_1,\ldots, i_k)$ is given by~\eqref{eq:thetai}.

\newpage

\begin{proposition}
\label{prop:Irrepresentable}
Let~$X$ and~$R$ be defined by~\eqref{eq:white}. Then, the following assumptions are equivalent:
\begin{itemize}
\item the design matrix~$X$ satisfies~\eqref{eq:Irr_matrix} of order~$K$,
\item the variance-covariance matrix~$R$ satisfies~\eqref{eq:Irrepresentable} of order~$K$.
\end{itemize}
Furthermore, they imply that for all~$(i_1,\ldots,i_{K})\in\mathcal A_K$ one has
\eq
\notag
\max\Big[
\max_{{j\neq i_1}}\theta_j(\imath_1),\,\ldots,\max_{j\neq i_1,\ldots, i_{K}}\theta_j(i_1,\ldots,i_{K})\Big]<1
\qe
which is an equivalent formulation of~$(i_1,\ldots,i_{K})$ satisfying~\eqref{hyp:IrrAlongThePath}.
\end{proposition}

\begin{proof}
Let~$S=\{j_1,\ldots,j_k\}\subset[p]$ and~$j\in[2p]\setminus S$.
Let~$\bar v=(\bar v_1,\ldots,\bar v_k)\in\{-1,1\}^k$ and define~$i_\ell = j_\ell+p(1-\bar v_\ell)/2$ for~$\ell\in[k]$.
Note that
\begin{align*}
\theta_j(i_1,\ldots, i_k) &= \big(R_{j,i_1} \cdots R_{j,i_{k}}  \big)M^{-1}_{i_1,\ldots, i_{k}}  
(1, \ldots, 1)\\
&= \Big[X_j^\top X_S\mathrm{Diag}(\bar v)\Big]M^{-1}_{i_1,\ldots, i_{k}}(1, \ldots, 1)\\
&=\Big[  X_j^\top X_S\mathrm{Diag}(\bar v)\Big]M^{-1}_{i_1,\ldots, i_{k}}\Big[\mathrm{Diag}(\bar v)\bar v\Big]\\
&= X_j^\top X_S\Big[\mathrm{Diag}(\bar v)M^{-1}_{i_1,\ldots, i_{k}}\mathrm{Diag}(\bar v)\Big]\bar v\\
&= X_j^\top X_S\big( X_S^\top  X_S\big)^{-1}\bar v.
\end{align*}
Now, observe that 
\[
\max_{||v||_\infty\leq1}v^\top\big( X_S^\top  X_S\big)^{-1} X_S^\top X_j
=\max_{\bar v\in\{-1,1\}^k} X_j^\top X_S\big( X_S^\top  X_S\big)^{-1}\bar v,
\]
showing the equivalence between the two assumptions.
\end{proof}

\begin{remark}
One may require that the design be `normalized' so that~$R_{i,i}=1$, namely its columns have unit Euclidean norm.
Under this normalization, one can check that~$R$ satisfies~\eqref{eq:Irr_matrix} of order~$K=1$.
Hence, up to some normalization, one can always assume~\eqref{eq:Irr_matrix} of order~$K=1$.
\end{remark}

\begin{remark}
When computing the LARS path, one has to compute the values 
\[
X_j^\top X_{\bar S^{k}}\big( X_{\bar S^{k}}^\top  X_{\bar S^{k}}\big)^{-1}\varepsilon^{k},
\] 
see for instance Algorithm~\ref{alg:LAR} or Algorithm~\ref{alg:LAR2}, where these values are given by~$\theta$, as shown by Proposition~\ref{prop:Irrepresentable}.
This implies that, in practice, along the LARS path, one witnesses the maximal order~$K$ for which Irrepresentable Check~\eqref{hyp:IrrAlongThePath} holds.
\end{remark}

\subsection{The estimator of the variance}
\label{sec:var_estimation}

 In our analysis, we introduce an estimate of the variance $  \hat \sigma^2$ to perform post-selection inference when the noise level $\sigma$ is unknown. The degree of freedom to the estimation of  the variance is~$n-K$. Let us fix, for the moment, $ j_1,\ldots, j_{K}$,   the indices that are  the putative indices for the  selected variables. Let~$P^\perp_{K}$ be the orthogonal projection on the orthogonal to~$H_K^{j_1,\ldots, j_{K}}:=\mathrm{Span}(X_{j_1}, \ldots, X_{ j_{K}})$. We define   
 \eq
      \label{e:def_sigmas}
     \hat \sigma^{j_1,\ldots, j_{K}}  := \frac{ ||P^\perp_{K} Y||_2 }{\sqrt{n-K}} .
    \qe
   By a slight abuse of notation, we can index the estimator above by the signed indexes~$i_1, \ldots, i_{K}$.
Eventually, we set 
\[
    \hat \sigma  := \hat \sigma^{\,\bar \imath_1,\ldots, \bar \imath_{K}},
\]
 the estimates of the standard deviation~$\sigma$.
 
 \subsection{Admissible Selection Procedures} \label{s:ap}\label{sec:exact_test_admissible}
 
 Note that choosing a model~$\hat S$ is equivalent to choosing a model size~$\hat m$ so that~
\eq
\label{e:selection_procedure}
\hat S =\big \{\bar  \imath_1,\bar  \imath_2,\ldots,\bar \imath_{\hat m}\big\}.
\qe
Our procedure is flexible on this point and allows any choice of~$\hat m$ as long as the following property~\eqref{e:stopping_rule} is satisfied:
\begin{quote}
{\normalsize
{\bf Stopping Rule:}
\it 
The estimated model size~$\hat m$ is a ‘{stopping time}':~$\hat m\in[K-1]$ and, for all~$a\in[K-1]$,
\eq
\tag{$\bf\mathcal A_{\mathrm{\,Stop}}$}
\label{e:stopping_rule}
\mathds 1_{\{\hat m\leq a\}}\text{~is a measurable function of~}(\lambda_1, \ldots, \lambda_a, \hat \imath_1,\ldots, \hat \imath_{K})\,.
\qe
}
\end{quote}

\noindent
In other words, the decision to select a model of size~$\{\hat m=a\}$ depends only on the first~$a$ variables entering the LARS.

\begin{remark}
We now give an example to show that {\eqref{e:stopping_rule}} implies some restriction.
Suppose, for example, that we want to decide wether the target~$\beta^0$ is two sparse or one sparse.
A natural decision rule is to look at large values of the second knot~$\lambda_2$, if~$\text{‘‘}\,\lambda_2 > \mathrm{(some\ threshold)}\text{''}$   choose~$ m=2$  otherwise choose~$m=1$.
  This rule does not satisfy~\eqref{e:stopping_rule}, since looking at~$\lambda_2$  we can choose only sizes~$m$ greater than or equal to 2.
\end{remark}

Denote by~$P_k(Y)$ (resp.~$P^\perp_k(Y)$) the orthogonal projection of the observation~$Y$ onto~$H_k$ (resp. the orthogonal of~$H_k$) for all~$k\geq 1$ where~$H_k$ are defined by~\eqref{eq:sequence_subspaces}. Given $h$ any measurable function, 
\[
\mathds 1_{\{\hat m\leq a\}}=h(P_a(Y)),
\]
determines a class of selection procedures satisfying~\eqref{e:stopping_rule}. These procedures decide whether to stop at~$\{\hat m=a\}$ based on the information given by~$P_a(Y)$. Once one has selected a model of size~$\hat m$, one may be willing to test if~$\hat S$ contains the true support~$S^0$ by considering the null hypothesis 
\[
\mathds H_0\,:\ \text{‘}S^0\subseteq \hat S\,\text{'},
\] 
namely there are no false negatives.
Equivalently, one aims at testing the null hypothesis
\eq
\label{e:nullH}
\mathds H_0\,:\ \text{‘}X\beta^0\in H_{\hat m}\,\text{'},
\qe
at an exact significance level~$\alpha\in(0,1)$, where~$(H_a)_{a=0}^{K-1}$ is defined by~\eqref{eq:sequence_subspaces}.

\section{Exact Controls using Least Angle Regression: Main Results} 
\label{sec3}

\subsection{Key notion: the ‘frozen' knots, their means and variances}
\label{sec:frozen_knots}

\subsubsection{Frozen knots}
Given~$K$ as defined in~\eqref{e:K} and fixed~$i_1,\ldots, i_{K+1}\in[2p]$, one may define
\begin{equation} \label{e:frozenZ}
\forall j \text{ s.t. } \theta_j(i_1, \ldots, i_{k})\neq1,\quad 
Z^{(i_1,\ldots, i_{k})}_j:=\frac{Z_j-\Pi_{ i_1,\ldots, i_{k}} (Z_j)}{1-\theta_j(i_1,\ldots, i_{k})},
\end{equation}
where
\eq
\label{e:Pi}
\Pi_{ i_1,\ldots, i_{k}} (Z_j) :=  \big(R_{j,i_1} \cdots R_{j,i_{k}}  \big)M^{-1}_{i_1,\ldots, i_{k}}  (Z_{i_1}, \ldots, Z_{i_{k}})
\qe
and~$\theta_j(i_1, \ldots, i_{k})$ is given by~\eqref{eq:thetai}. When~$\bbE Z=0$, one may remark that~$\Pi_{ i_1,\ldots, i_{k}} (Z_j)$ is the regression of~$Z_j$ on the vector~$ (Z_{i_1}, \cdots , Z_{i_{k}})$ whose variance-covariance matrix is~$\sigma^2M_{i_1,\ldots, i_{k}}$, namely
\[
\text{When~}\bbE Z=0,\quad
\Pi_{ i_1,\ldots, i_{k}} (Z_j)=\big(R_{j,i_1} \cdots R_{j,i_{k}}  \big)M^{-1}_{i_1,\ldots, i_{k}}   (Z_{i_1}, \ldots, Z_{i_{k}})=\mathds E\big[Z_j| Z_{i_1}, \cdots , Z_{i_{k}}\big].
\]

\noindent
From this point on, we introduce
\eq
\label{e:frozenLambda}
\forall k\geq 0,\quad
\lambda_{k+1}^{(i_1,\ldots, i_{k})}:=\max_{j: \theta_j(i_1, \ldots, i_{k})<1}Z^{(i_1,\ldots, i_{k})}_j,
\qe
and we emphasize that
\eq
\label{e:unfrozen}
\forall k\geq 0,\quad
\lambda_{k+1}\,\mathds 1_{\{\hat \imath_1=i_1,\ldots, \hat \imath_{k}=i_k\}}=\lambda_{k+1}^{(i_1,\ldots, i_{k})}\,\mathds 1_{\{\hat \imath_1=i_1,\ldots, \hat \imath_{k}=i_k\}},
\qe
as proven in Appendix~\ref{app:recursion} (Eq.~\eqref{e:lambda_max}) and Proposition~\ref{prop:reccur}.
We are now  able to define the ‘‘{\it frozen}''  values of the knots:
\eq \label{e:frozen}
 \lambda^f_1 : = Z_{i_1}, \ldots , \lambda^f_{K+1} :=   Z^{i_1,\ldots i_K}_{i_{K+1}}.
  \qe
They are the Gaussian random variables that coincide  with $\lambda_1,\lambda_2,\cdots ,  \lambda_{K+1}$  when the random variables defined by the signed indices~$\hat  \imath_1,\hat  \imath_2, \ldots,  \hat \imath_{K+1}~$ take the particular  values~$i_1, i_2,\ldots, i_{K+1}$.

\begin{remark}
An interesting feature of the LARS knots is that they have a simple expression in terms of the partition given by the identity
\[
\sum_{(i_1,\ldots,i_{K})\in\cA_{K}}\mathds 1_{\{\hat \imath_1=i_1,\ldots,\hat \imath_{K}=i_{K}\}}=1\quad\text{almost\ surely}.
\]
As we have seen in~\eqref{e:frozen},
\[
\forall k\in [K],\quad
\lambda_k=\sum_{(i_1,\ldots,i_{k}) \in\cA_{k}}\mathds 1_{\{\hat \imath_1=i_1,\ldots,\hat \imath_{k}=i_{k}\}}\underbrace{Z^{i_1,\ldots i_{k-1}}_{i_{k}}}_{=:\lambda^f_k},
\]
giving the definition of the {\it frozen} knots~$\lambda^f_k$ above.
\end{remark}

\subsubsection{Mean and centering of the frozen knots}
Now, write 
\eq
\label{e:P_seq}
\forall y\in\mathds R^n,\quad
P^{(i_1,\ldots, i_{k})}(y)=\big(X_{j_1} \cdots X_{j_{k}}  \big)M^{-1}_{j_1,\ldots, j_{k}}  (X_{j_1}^\top, \ldots, X_{j_{k}}^\top)\,y
\qe
for the orthogonal projection of $y$ onto~$\mathrm{Span}(X_{j_1} ,\ldots, X_{j_{k}} )$.
Recall that $P_k$ is the orthogonal projection onto~$H_k$, which is a random subspace.
Recall also that, conditional on the event $\{\hat \imath_1=i_1,\ldots, \hat \imath_{k}=i_k\}$ the subspace~$H_k$ is fixed. Note that 
\begin{align*}
P_k\,\mathds 1_{\{\hat \imath_1=i_1,\ldots, \hat \imath_{k}=i_k\}}&=P^{(i_1,\ldots, i_{k})}\,\mathds 1_{\{\hat \imath_1=i_1,\ldots, \hat \imath_{k}=i_k\}},\\
P^\perp_k\,\mathds 1_{\{\hat \imath_1=i_1,\ldots, \hat \imath_{k}=i_k\}}&=\big(\mathrm{Id}_n-P^{(i_1,\ldots, i_{k})}\big)\,\mathds 1_{\{\hat \imath_1=i_1,\ldots, \hat \imath_{k}=i_k\}},\\
\text{and for all }i\in [2p],\quad \Pi_{ i_1,\ldots, i_{k}} (Z_i)
&= s\langle X_j, P^{(i_1,\ldots, i_{k})}(Y)\rangle,
\end{align*}
where~$i=j+p({1-s})/2$.
The mean~$m^f_k$ and standard deviation~$\sigma\,\rho^f_k$ of~$\lambda^f_k$ are important values defined for all~$k\in[K]$:
\begin{align*}
m_{k}\,\mathds 1_{\{\hat \imath_1=i_1,\ldots, \hat \imath_{k}=i_k\}}&
=m_{k}^f\,\mathds 1_{\{\hat \imath_1=i_1,\ldots, \hat \imath_{k}=i_k\}},
\\
\rho_{k}\,\mathds 1_{\{\hat \imath_1=i_1,\ldots, \hat \imath_{k}=i_k\}}&=\rho_{k}^{f}\,\mathds 1_{\{\hat \imath_1=i_1,\ldots, \hat \imath_{k}=i_k\}},
\end{align*}
with 
\begin{align}
\label{e:frozen_m}
m_{k}^f
&=
\frac{
s_k\langle X_{j_k},\big(\mathrm{Id}_n-P^{(i_1,\ldots, i_{k-1})}\big)X\beta^0\rangle
}{1-\theta_{i_{k}}(i_1,\ldots,i_{k-1})},
\\
\label{e:frozen_rho}
\rho_{k}^{f}
&=
\frac{
\sqrt{
\langle X_{j_k},\big(\mathrm{Id}_n-P^{(i_1,\ldots, i_{k-1})}\big)X_{j_k}\rangle
}}{1-\theta_{i_{k}}(i_1,\ldots,i_{k-1})},
\end{align}
and this definition is equivalent to~\eqref{eq:def_mu_proj} and~\eqref{eq:def_rho}, see Section~\ref{subsec:main}.
Recall that~$\bar\mu^0=\bar R\beta^0$ as defined in~\eqref{e:mu0} and note that
\eq
\label{e:meanH0}
m^f_k=0
\quad\Leftrightarrow\quad 
\bar\mu^0_{i_k}-\Pi_{ i_1,\ldots, i_{k-1}} (\bar\mu^0_{i_k})=0
\quad\Leftrightarrow\quad 
\langle X_{j_k},\big(\mathrm{Id}_n-P^{(i_1,\ldots, i_{k-1})}\big)X\beta^0\rangle=0,
\qe
which is true when the true support~$S^0$ of~$\beta^0$ is included in~$\bar S^k$, defined by~\eqref{eq:sequence_supports}.
 This proves the next proposition.

\begin{proposition}
\label{prop:H0}
For fixed~$0\leq a\leq K-1$, conditional on the selection event~$\{\hat \imath_1 =i_1, \ldots, \hat \imath_{K}=i_{K}\}$, the hypothesis 
\eq
\notag
\mathds H_0\,:\ \text{‘‘}X\beta^0\in H_{a}\text{''}
\qe
implies that~$m^f_k=0$ for all~$a<k\leq K$, namely~$(Z_{i_{a+1}}^{(i_1,\ldots, i_{a})},\ldots,Z_{i_{K}}^{(i_1,\ldots, i_{K-1})})$ is centered.
\end{proposition}

\noindent
This proposition is important for defining the hypothesis under consideration, see also Remark~\ref{rem:selective_testing}.

\subsubsection{A key result: The characterization of the selection event}
\label{sec:selection_event}
Regarding the joint law of the frozen knots, one has the following important proposition whose proof can be found in Section~\ref{proof:reccur}.

\begin{proposition}
\label{prop:reccur}
Let~$(i_1,\ldots,i_K,i_{K+1})\in\mathcal A_{K+1}$, that is, a fixed active set of size~$K+1$.
\begin{itemize}
\item  If~$(i_1,\ldots,i_K)$ satisfies~\eqref{hyp:IrrAlongThePath} then
\begin{align*}
&\big\{\hat \imath_1=i_1,\ldots,\hat \imath_{k+1}=i_{k+1}\big\}\\
&=\big\{
\lambda_{k+1}^{(i_1,\ldots, i_{k})}=Z_{i_{k+1}}^{(i_1,\cdots, i_{k})}
\leq
Z_{i_k}^{(i_1,\cdots, i_{k-1})}\leq\cdots\leq
Z_{i_2}^{(i_1)}\leq
Z_{i_1}\big\}\\
&
=\big\{
\lambda_{k+1}^{(i_1,\ldots, i_{k})}=Z_{i_{k+1}}^{(i_1,\cdots, i_{k})}
\leq
Z_{i_k}^{(i_1,\cdots, i_{k-1})}\leq\cdots\leq
Z_{i_{a+1}}^{(i_1,\cdots, i_{a})}\leq
Z_{i_{a}}^{(i_1,\cdots, i_{a-1})}=\lambda_{a}^{(i_1,\ldots, i_{a-1})}\big\}\\
&
\quad
\bigcap
\big\{\lambda_{k+1}=Z_{i_{k+1}}^{(i_1,\cdots, i_{k})},\ldots,\lambda_{a}=Z_{i_{a}}^{(i_1,\cdots, i_{a-1})}\big\}
\bigcap
\big\{\hat \imath_1=i_1,\ldots,\hat \imath_{a}=i_{a}\big\},
\end{align*}
for any~$0\leq a< k\leq K$ with the convention~$\lambda_0=\infty$.
\item It holds that 
\[
(Z^{(i_1,\ldots, i_{k})}_j)_{ j\neq i_1, \ldots, i_{k}}\indep Z_{i_{k}}^{(i_1,\ldots, i_{k-1})}
\indep Z_{i_{k-1}}^{(i_1,\ldots, i_{k-2})}
\indep \cdots 
\indep Z_{i_{2}}^{(i_1)}
\indep Z_{i_1}
\]
are mutually independent, for any~$k\in[K]$. 
\newpage
Furthermore, if $X\beta^0\in H_K$ then
\eq
\label{e:indep_sigmas}
{\hat \sigma^{i_1,\ldots, i_{K}}}
\indep 
\Big(\frac{Z^{(i_1,\ldots, i_{K})}_j}{\hat \sigma^{i_1,\ldots, i_{K}}}\Big)_{ j\neq i_1, \ldots, i_{K}}
\indep {Z_{i_{K}}^{(i_1,\ldots, i_{K-1})}}
\indep \cdots 
\indep {Z_{i_{2}}^{(i_1)}}
\indep {Z_{i_1}}.
\qe
\item  If~$(i_1,\ldots,i_K)$ satisfies~\eqref{hyp:IrrAlongThePath} then 
\begin{align*}
&\big\{\hat \imath_1=i_1,\ldots,\hat \imath_{K}=i_{K}\big\}\\
&=\Big\{
\lambda_{K+1}^{(i_1,\ldots, i_{K})}
\leq
Z_{i_K}^{(i_1,\cdots, i_{K-1})}\leq\cdots\leq
Z_{i_2}^{(i_1)}\leq
Z_{i_1}\Big\}\\
&=\Bigg\{
{\Lambda_{K+1}^{(i_1,\ldots, i_{K})}:=\frac{\lambda_{K+1}^{(i_1,\ldots, i_{K})}}{\hat \sigma^{i_1,\ldots, i_{K}}}}
\leq
\frac{Z_{i_{K}}^{(i_1,\ldots, i_{K-1})}}{\hat \sigma^{i_1,\ldots, i_{K}}}
\leq\cdots\leq
\frac{Z_{i_{2}}^{(i_1)}}{\hat \sigma^{i_1,\ldots, i_{K}}}
\leq
\frac{Z_{i_1}}{\hat \sigma^{i_1,\ldots, i_{K}}}
\Bigg\}\,.
\end{align*}
\end{itemize}
\end{proposition}

\begin{remark} 
{\bf Is Proposition~\ref{prop:reccur} a new polyhedral lemma?}
The characterization of the selection event for the inference of a single testing statistic has been known as the ‘polyhedral lemma' in the literature, see for instance \cite[Figure~6.9]{hastie2015statistical} and references therein. This result is the cornerstone of selective inference with sparse models. It is based on two ingredients: First, the selection event can be expressed as a polyhedra; Second, conditional on the selection event, any linear statistics is distributed according to a truncated Gaussian with independent truncation bounds. 

A first remark is that the polyhedral lemma is shown for one linear statistic and, as far as we known, there is no polyhedral lemma for multiple linear statistics. The interesting point is that our result (Proposition~\ref{prop:reccur}) can be seen as a polyhedral lemma for multiple linear statistics. Under \eqref{hyp:IrrAlongThePath}, the selection event $\{\hat \imath_1,\ldots, \hat \imath_K\}$ corresponds to a polyhedra described by the $Z_{i_k}^{(i_1,\cdots, i_{k-1})}$ variables in the third point of Proposition \ref{prop:reccur}. Our main result shows that the joint law of these multiple linear statistics are the Gaussian distribution restricted to the polyhedra $\{\ell_1\geq\ldots\geq \ell_K\geq\lambda_{K+1}\}$, see Theorem~\ref{thm:Main1}.

Note that the selection event has to include $\lambda_{K+1}$. As discussed above, our polyhedral lemma (Theorem \ref{thm:Main1} and Proposition~\ref{prop:reccur}) shows that, conditional on the selection event, $\lambda_1,\ldots,\lambda_K$ are distributed with respect to a Gaussian distribution restricted to the polyhedra $\{\ell_1\geq\ldots\geq \ell_K\geq\lambda_{K+1}\}$. If one does not include $\lambda_{K+1}$ in the selection event, then one has to integrate this latter conditional law with respect to the distribution of $\lambda_{K+1}$ which is not known. 
\end{remark}

\medskip

\begin{proposition}
Assume that the design~$X$ is such that the Irrepresentable Condition~\eqref{eq:Irr_matrix} of order~$K$ holds.
Almost surely, one has
\begin{itemize}
\item 
Among all possible sets~$(i_1,\ldots,i_K)\in\cA_K$, there is one and only one such that 
\eq
\label{eq:combinatorial}
\max_{i_{K+1}\neq i_1,\ldots,i_K}Z_{i_{K+1}}^{(i_1,\cdots, i_{K})}\leq Z_{i_K}^{(i_1,\cdots, i_{K-1})}\leq\cdots\leq Z_{i_2}^{(i_1)}\leq Z_{i_1}.
\qe
\item 
This set is the set selected by LARS, namely~$\hat \imath_1=i_1,\ldots,\hat \imath_{K}=i_{K}$,
\item 
and, for all~$(i_1,\ldots,i_K)\in\cA_K$, 
\[
\mathds P\big(\hat \imath_1=i_1,\ldots,\hat \imath_{K}=i_{K}\big)
=
\mathds P\Big(\max_{i_{K+1}\neq i_1,\ldots,i_K}Z_{i_{K+1}}^{(i_1,\cdots, i_{K})}\leq Z_{i_K}^{(i_1,\cdots, i_{K-1})}\leq\cdots\leq Z_{i_2}^{(i_1)}\leq Z_{i_1}\Big).
\]
\end{itemize}
\end{proposition}
\begin{proof}
Note that~\eqref{eq:Irr_matrix}  implies~\eqref{hyp:IrrAlongThePath} by Proposition~\ref{prop:Irrepresentable}.
Then apply the first point of Proposition~\ref{prop:reccur} to conclude.
\end{proof}

\newpage

\noindent
Finding the set~$\{\hat \imath_1=i_1,\ldots,\hat \imath_{K}=i_{K}\}$ selected by LARS may be related to a combinatorial search testing~\eqref{eq:combinatorial} all possible candidates~$(i_1,\ldots,i_K)\in\cA_K$.
Under the Irrepresentable Condition, the support selected by LARS is given by~\eqref{eq:combinatorial}, which can be seen as the extension of~\eqref{eq_simplest} introducing~{\bf [Q1]} in Section~\ref{sec:Q1}.

\subsection{Main results: Joint law and construction of post-selection tests}
\label{subsec:main}

We assume that~$K$ is defined as in~\eqref{e:K}.
Except in Section~\ref{sec:FalseNeg},~$\sigma^2$ is  assumed to  be known.
Let~$(\,\hat \imath_1,\ldots, \hat \imath_K)$ be the first signed variables entering along the~LARS path.
In this section, we are interested in the joint law of the LARS knots~$(\lambda_1,\ldots, \lambda_K)$ conditional on~$\lambda_{K+1}$ and~$(\,\hat \imath_1,\ldots, \hat \imath_K)$.
To determine this joint law, we need to make precise the centering parameters~$m_k$, by (see also~\eqref{e:frozen_m})
\eq
\label{eq:def_mu_proj}
m_k:=
\frac{
\mu^0_{\hat \imath_k}- \big(R_{\hat \imath_k,\hat \imath_1} \cdots R_{\hat \imath_{k},\hat \imath_{k-1}}  \big)
M^{-1}_{\hat \imath_1,\ldots, \hat \imath_{k-1}}
\big(\mu^0_{\hat \imath_1} ,\cdots ,\mu^0_{\hat \imath_{k-1}} \big)
}{1-\theta_{\hat \imath_{k}}^{k-1}}
\qe
the first standard deviation~$\sigma\rho_1$ with~$\rho_1:=\sqrt{R_{\hat \imath_1,\hat \imath_1}}$, and the others~$\sigma\rho_k$ by
(see also~\eqref{e:frozen_rho})
\begin{align}
\rho_\ell
&
:=\frac{ \sqrt{ R_{\hat \imath_\ell,\hat \imath_\ell} - \big(R_{\hat \imath_\ell,\hat \imath_1} \cdots R_{\hat \imath_{\ell},\hat \imath_{\ell-1}}  \big)M^{-1}_{\hat \imath_1,\ldots, \hat \imath_{\ell-1}} \big(R_{\hat \imath_\ell,\hat \imath_1} ,\cdots ,R_{\hat \imath_{\ell},\hat \imath_{\ell-1}} \big)}} { 1-\theta^{\ell-1}_{\hat \imath_\ell}}
\quad \text{for }2\leq \ell\leq K+1,
\label{eq:def_rho}
\end{align}
where~
\[
\theta^{\ell-1}:=\theta(\,\hat \imath_1,\ldots, \hat \imath_{\ell-1}),\quad \text{for }2\leq \ell\leq K+1, 
\]
is defined by~\eqref{eq:thetai} and~$M_{\hat \imath_1,\ldots, \hat \imath_{\ell-1}}$ is defined by~\eqref{e:M}.


\subsubsection{Proof of Theorem~\ref{thm:Main1}}
\label{proof:Main1}
From the definition of the Gaussian random variable~$ Z_{i_k}^{(i_1,\ldots, i_{k-1})}~$ in~\eqref{e:frozenZ} one can deduce that its mean~$m_k$ is given by~\eqref{eq:def_mu_proj} and its standard deviation~$v_k$ by~\eqref{eq:def_rho}, considering putative indices for the selected variables. By the second point of Proposition~\ref{prop:reccur}, we know that these variables are independent.
We deduce that their joint density~$(Z_{i_1}, Z_{i_{2}}^{(i_1)}, \ldots, Z_{i_{K}}^{(i_1,\ldots, i_{K-1})})$ is 
\[
\prod_{k=1}^K \varphi_{m_k,v_k^2}(\ell_k),
\]
 with respect to Lebesgue measure.
For now on, we condition on  $\mathcal E:=\{\hat \imath_1 =i_1, \ldots,\hat \imath_{K}=i_{K},\lambda_{K+1}\}$ and we assume that~$(i_1,\ldots, i_K)$ satisfies~\eqref{hyp:IrrAlongThePath}. By the first equality of the third point of Proposition \ref{prop:reccur} we known that~$\mathcal E=\big\{\lambda_{K+1}\leq Z_{i_{K}}^{(i_1,\ldots, i_{K-1})} \leq  \cdots \leq Z_{i_1}\big\}$, and on the event~$\mathcal E$,
\eq
\label{e:equal}
(Z_{i_1}, Z_{i_{2}}^{(i_1)}, \ldots, Z_{i_{K}}^{(i_1,\ldots, i_{K-1})})
=
(\lambda_1,\lambda_2,\ldots,\lambda_K).
\qe
Conditional on~$\mathcal E$, the joint density of~$(Z_{i_1}, Z_{i_{2}}^{(i_1)}, \ldots, Z_{i_{K}}^{(i_1,\ldots, i_{K-1})})$ is proportional to
\eq
 \label{zaz0}
\Big(
\prod_{k=1}^K \varphi_{m_k,v_k^2}(\ell_k)
\Big)
\,
\mathds{1}_{\{\ell_1\geq \ell_{2}\geq\cdots\geq\ell_{K}\geq\lambda_{K+1}\}},
\qe
 with respect to Lebesgue measure, and by~\eqref{e:equal} it is the conditional density of the knots.

\subsubsection{Construction of the Generalized Spacing test}
\label{sec:fabc}
A useful consequence of Theorem~\ref{thm:Main1} is that one can explicitly describe the joint law of the LARS knots after having selected a support~$\hat S$ of size~$\hat m$ with any procedure satisfying~{\eqref{e:stopping_rule}}.
In the sequel, we write
\eq
\label{eq:pij}
F_i:=\Phi_{i}({\lambda_i}) :=\Phi\Big(\frac{\lambda_i}{\sigma{\rho_i}}\Big)
\quad\text{and}\quad \cP_{i,j}:=\Phi_{i}\circ\Phi_{j}^{-1},\quad \text{for }i,j\in[K+1],
\qe
 where
 $\lambda_0=\infty$ and~$F_0=1$ by convention.

\begin{proposition}
\label{cor:Main1}
Let~$a\in\bbN$ be such that~$0\leq a\leq K-1$.
Let~$\hat m$ be a selection procedure satisfying~\eqref{e:stopping_rule}.
Under the conditions of Theorem \ref{thm:Main1}, under the null hypothesis 
\eq
\label{e:H0}
\mathds H_0\,:\ \text{‘‘}X\beta^0\in H_{a}\text{''},
\qe
and conditional on the selection event 
$\big\{\hat m=a, F_{a},F_{K+1},\hat \imath_1,\ldots, \hat \imath_K\big\}$, we have that $(F_{a+1},\ldots,F_K)$ is uniformly distributed on 
\begin{align*}
\cD_{a+1,K}:=&
\big\{
(f_{a+1},\ldots,f_K)\in\bbR^{K-a}:\ \\
&\cP_{a+1,a}(F_{a})\geq f_{a+1}\geq \cP_{a+1,a+2}(f_{a+2})
\geq \cdots
\geq \cP_{a+1,K}(f_K)\geq \cP_{a+1,K+1}(F_{K+1})
\big\},
\end{align*}
where the~$\cP_{i,j}$ are described in~\eqref{eq:pij}.
\end{proposition}
\noindent
A proof of this proposition can be found in Appendix~\ref{proof:cor4}.


\begin{remark}
The previous statement is consistent with the case~$a=0$ corresponding to~the {\rm global null} hypothesis~$\mathds H_0\,:\, \text{`}X\beta^0=0\text{'}$~$($or equivalently~$\bbE Z=0)$.
Therefore, if~$Z$ is centered, then, conditional on~$F_{K+1}$, one has that~$(F_1,\ldots,F_K)$ is uniformly distributed on 
\[
\cD_{1,K}:=
\big\{
(f_1,\ldots,f_K)\in\bbR^K\ :\ 
1\geq f_1\geq \cP_{1,2}(f_2)\geq \cdots\geq \cP_{1,K}(f_K)\geq \cP_{1,K+1}(F_{K+1}))
\big\}.
\]
\end{remark}

\begin{remark}
\label{rem:stat_order}
In the orthogonal case, where~$\bar R=\mathrm{Id}$, note that~$\theta_j(i_1,\ldots,i_\ell)=0$ for all~$\ell\geq1$ and all~$i_1,\ldots,i_\ell\neq j$,~$\rho_j=1$ and~$\cP_{i,j}(f)=f$.
We recover that~$\cD_{1,K}$ is the set of order statistics~
\[
1\geq f_1\geq f_2\geq\ldots\geq f_K\geq  \Phi(\lambda_{K+1}/\sigma).
\] 
In this case, the knots~$\lambda_i$ are Gaussian order statistics~$\lambda_1=Z_{\hat \imath_1}\geq \lambda_2=Z_{\hat \imath_2}\geq \ldots\geq \lambda_K=Z_{\hat \imath_K}\geq \lambda_{K+1}$ for the vector~$Z$.
\end{remark}

\noindent
From Theorem~\ref{thm:Main1}, we deduce several test statistics.
To this end, we introduce some notation defining
\begin{align}
\label{fun:iabc}
\cI_{ab}(s,t)&:=
\displaystyle
\int\displaylimits_{\cP_{{(a+1)},b}(t)}^{\cP_{{(a+1)},a}(s)}\!\!\!\mathrm df_{a+1}\!\!\!
\int\displaylimits_{\cP_{{(a+2)},b}(t)}^{\cP_{(a+2),(a+1)}(f_{a+1})}\!\!\!\!\mathrm df_{a+2}\!\!\!
\int\displaylimits_{\cP_{{(a+3)},b}(t)}^{\cP_{(a+3),(a+2)}(f_{a+2})}\!\!\!\!\!\mathrm df_{a+3}\,\,\,
\cdots
\!\!\!\int\displaylimits_{\cP_{{(b-1)},b}(t)}^{\cP_{(b-1),(b-2)}(f_{b-2})}\!\!\!\!\!\mathrm df_{b-1}
\\
&\text{for }0\leq a<b\ \text{and }s,t\in\bbR,\text{ with the convention that }\cI_{ab}=1\text{ when }b=a+1,
\notag
\end{align}
\noindent
and also
\begin{align}
\label{fun:fabc}
\bbF_{abc}(t)&:=
\1_{\{\lambda_c\leq t\leq \lambda_a\}}
\displaystyle
\int\displaylimits_{\Phi_{b}(\lambda_{c})}^{\Phi_{b}(t)}\!\!\cI_{ab}(
{F_a},f_b)\,\cI_{bc}(f_b,
{F_c})\,\mathrm df_b\\
\notag
&\text{for }0\leq a<b<c\leq K+1,\ t\in\bbR\, \text{ where }F_a=\Phi_{a}(\lambda_a)\text{ and }F_c=\Phi_{c}(\lambda_c).
\end{align}

\begin{remark}
On the numerical side, note that this quantity can be computed using {\it Quasi Monte Carlo} (QMC) methods as in~\cite[Chapter 5.1]{genz2009computation} or Appendix~\ref{app:cbc}.
The function~$\bbF_{abc}$ gives the CDF of~$\lambda_b$ conditional on~$\lambda_a,\lambda_c$ and on some selection event, as shown in the next proposition.
\end{remark}


\begin{proposition}
\label{prop:fabc}
Let~$a, b,$ and $c$ be such that~$0\leq a<b<c\leq K+1$.
Let~$(\lambda_1,\ldots, \lambda_K,\lambda_{K+1})$ be the first knots and let~$(\,\hat \imath_1,\ldots, \hat \imath_K)$ be the first variables entering along the~LARS path.
If~$(\,\hat \imath_1,\ldots, \hat \imath_K)$ satisfies~\eqref{hyp:IrrAlongThePath} and~$\hat m$ is chosen according to a procedure satisfying~{\eqref{e:stopping_rule}}, then under the null hypothesis
\[
\mathds H_0\,:\ \text{‘‘}X\beta^0\in H_{a}\,\text{''},
\] 
it holds that
\eq
\label{eqFabc}
\bbP\big[\lambda_b\leq t\ |\ \hat m=a,\lambda_a,\lambda_c, \hat \imath_{1}, \ldots,\hat \imath_{c-1}\big]=\frac{\bbF_{abc}(t)}{\bbF_{abc}(\lambda_a)}.
\qe
\end{proposition}

\noindent
A proof of this proposition can be found in Appendix~\ref{proof:fabc}.

\begin{remark}
Note that the deterministic choice~$\hat m=a$, for a fixed~$a\in[K-1]$, is a procedure satisfying~{\eqref{e:stopping_rule}} and Proposition \ref{prop:fabc} holds.
This shows that if~$(\,\hat \imath_1,\ldots, \hat \imath_K)$ satisfies~\eqref{hyp:IrrAlongThePath}, then, under the null hypothesis~$\mathds H_0\,:\ \text{‘‘}X\beta^0\in H_{a}\,\text{''}$,
\eq
\label{e:fabc}
\bbP\big[\lambda_b\leq t\ |\ \lambda_a,\lambda_c, \hat \imath_{1}, \ldots,\hat \imath_{c-1}\big]=\frac{\bbF_{abc}(t)}{\bbF_{abc}(\lambda_a)},
\qe
for any~$0\leq a<b<c\leq K+1$.
\end{remark}

\subsubsection{Proof of Theorem~\ref{nt:pvalue}}
\label{proof:nt:pvalue}
From Proposition \ref{prop:fabc} we know that under~$\mathds H_0$ and conditional on the selection event~$\big\{\hat m=a\big\}$, Eq.~\eqref{eqFabc} gives the conditional CDF of~$\lambda_b$.
As a consequence, 
  \[
  \frac{
  \bbF_{abc}(\lambda_b)    }{    \bbF_{abc}(\lambda_a)  }  \sim  \mathcal  U(0,1).
  \]
  Finally, considerations of the distribution under the alternative  show that to obtain a~$p$-value, we must consider the complement to 1 of  the quantity above.

\subsubsection{Monte Carlo simulations, Spacing tests, and Generalized Spacing tests}
Theorem~\ref{nt:pvalue} is illustrated numerically in Figure \ref{fig:law}.
Note that we have a perfect fit with the uniform law: the conditional law of the LARS knots obtained theoretically is numerically validated\footnote{A reproducible experiment given in a Python notebook is available at \url{https://github.com/ydecastro/lar_testing/blob/master/Law_LAR.ipynb}}. This test statistic generalizes previous test statistics that appeared in `\textit{Spacing Tests}', as presented in~\cite[Chapter~5]{hastie2015statistical} for instance, and will be referred to as the {\it Generalized Spacing test}.

\begin{remark}
\label{rem:ST}
If one takes~$a=0$,~$b=1$, and~$c=2$ then
\[
\hat\alpha_{012}=1-\frac{\Phi_{1}(\lambda_1)-\Phi_{1}(\lambda_2)}{\Phi_{1}(\lambda_0)-\Phi_{1}(\lambda_2)}
=\frac{1-\Phi_{1}(\lambda_1)}{1-\Phi_{1}(\lambda_2)}.
\]
Similarly, taking~$b=a+1$ and~$c=a+2$,
\[
\hat\alpha_{a(a+1)(a+2)}
=\frac
{\Phi_{{a+1}}(\lambda_{a+1})-\Phi_{{a+1}}(\lambda_a)}
{\Phi_{{a+1}}(\lambda_{a+2})-\Phi_{{a+1}}(\lambda_a)}.
\]
which is the {\it conservative spacing test}, see \cite[Theorem~2]{taylor2014exact}.
\end{remark}

\newpage

\begin{figure}[t]
\includegraphics[width=0.9\textwidth]{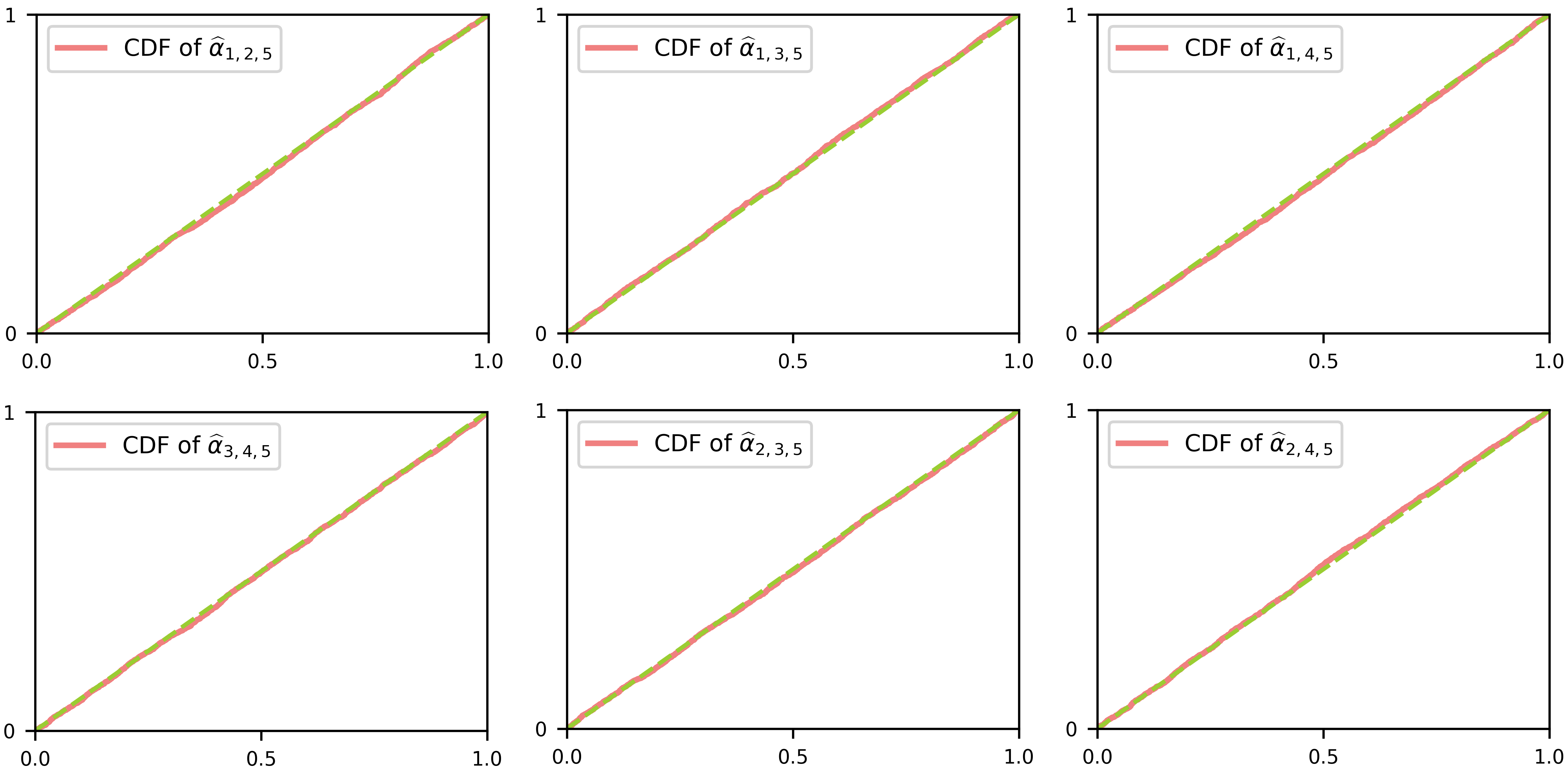}
\caption{Observed empirical law of~$\hat\alpha_{abc}$ over~$5,000$ Monte Carlo repetitions with~$n=200$ and~$p=300$.
We considered a design~$X\in\bbR^{n\times p}$ with independent column vectors uniformly distributed on the sphere and an independent~$y\in\bbR^{n}$ with i.i.d. standard Gaussian entries, and we computed the indices~$(\,\hat \imath_1,\ldots,\hat \imath_n)$ and the knots~$(\lambda_1,\ldots,\lambda_n)$ entering the model with LARS.
The empirical CDF of the~$\hat\alpha_{abc}$ are displayed.
We observe a perfect fit with the uniform distribution: the conditional law of the LARS knots obtained theoretically is numerically validated.}
\label{fig:law}
\end{figure}

\subsection{Exact false negative testing after model selection}
\label{sec:exact_model_selection}
We return to the case of a general design.
Given~$\alpha\in(0,1)$ and using Theorem~\ref{nt:pvalue}, one can consider the following exact testing procedure at level~$\alpha$ on false negatives, see the pseudo-code in Algorithm~\ref{Alg:TestingAStop}.
The theoretical guarantee of this algorithm is given by the next proposition.
It shows that conditional on the event that `there are no false negatives', namely~$\text{`}X\beta^0\in H_{\hat m}\,\text{'}$, the observed significance~$\hat \alpha$ obeys the uniform law and hence~$\1_{\{\hat \alpha\leq\alpha\}}$ is a testing procedure with level exactly~$\alpha$.

\pagebreak[3]

\begin{corollary}
\label{prop:unconditional_test}
Let~$(\lambda_1,\ldots, \lambda_K,\lambda_{K+1})$ be the first knots and let~$(\,\hat \imath_1,\ldots, \hat \imath_K)$ be the first variables entering along the~LARS path.
If~$(\,\hat \imath_1,\ldots, \hat \imath_K)$ satisfies~\eqref{hyp:IrrAlongThePath} and~$\hat m$ is chosen according to a procedure satisfying~{\eqref{e:stopping_rule}}, then, conditional on the null hypothesis
\[
\mathds H_0\,:\ \text{‘‘}X\beta^0\in H_{\hat m}\,\text{''},
\] 
it holds that
\[
\hat\alpha_{\hat m(\hat m+1)(K+1)}:=1
-\frac{\bbF_{\hat m(\hat m+1)(K+1)}(\lambda_{\hat m+1})}{\bbF_{\hat m(\hat m+1)(K+1)}(\lambda_{\hat m})}
\sim\cU(0,1),
\]
that is, it is uniformly distributed over~$(0,1)$.
\end{corollary}

\pagebreak[3]

\begin{proof}
By Theorem \ref{nt:pvalue}, the conditional law of~$\hat\alpha_{a(a+1)(K+1)}$ with respect to~$\{\hat m\leq a\}$ is the uniform distribution.
Note that the conditional law~\eqref{zaza5} does not depend on~$a,b=a+1,c=K+1$, hence this law is unconditional on~$\hat m$.
\end{proof}

\noindent
When the variance~$\sigma^2$ is unknown, one can `Studentize' this test, as presented in the next section.
The reader may consult Section~\ref{sec:Student} for a definition and check that the quantities~$\hat\beta_{abc},\tilde\bbF,\Lambda_k$ do not require~$\sigma$ to be computed.

\begin{corollary}
\label{prop:unconditional_t_test}
Let~$(\lambda_1,\ldots, \lambda_K,\lambda_{K+1})$ be the first knots and let~$(\,\hat \imath_1,\ldots, \hat \imath_K)$ be the first variables entering along the~LARS path.
If~$(\,\hat \imath_1,\ldots, \hat \imath_K)$ satisfies~\eqref{hyp:IrrAlongThePath} and~$\hat m$ is chosen according to a procedure satisfying~{\eqref{e:stopping_rule}}, then, conditional on the null hypothesis
\[
\mathds H_0\,:\ \text{‘‘}X\beta^0\in H_{\hat m}\,\text{''},
\] 
it holds that
\[
\hat\beta_{\hat m(\hat m+1)(K+1)}:=1
-\frac{\tilde\bbF_{\hat m(\hat m+1)(K+1)}(\Lambda_{\hat m+1})}{\tilde\bbF_{\hat m(\hat m+1)(K+1)}(\Lambda_{\hat m})}
\sim\cU(0,1),
\]
that is, it is uniformly distributed over~$(0,1)$.
\end{corollary}

\begin{proof}
By Theorem \ref{t:pvalue}, the conditional law of~$\hat\beta_{a(a+1)(K+1)}$ with respect to~$\{\hat m\leq a\}$ is the uniform distribution. Note that the conditional law~\eqref{e:talpha} does not depend on~$a,b=a+1,c=K+1$, hence this law is unconditional on~$\hat m$.
\end{proof}

\subsection{Exact Testing Procedure for False Negatives with Variance Estimation}
\label{sec:FalseNeg}\label{sec:Student}

From the results of Section~\ref{subsec:main}, one can present a method to select a model and propose an exact test of false negatives in the case of a general design, when the variance is unknown.
We introduce a new exact testing procedure that can be deployed when {\eqref{e:stopping_rule}} holds, namely an `admissible' selection procedure is used to build~$\hat S$. We start by a preliminary result whose proof is in Appendix~\ref{proof:frozen}.
 
   \begin {proposition}\label{p:frozen}
Let~$(\lambda_1,\ldots, \lambda_K,\lambda_{K+1})$ be the first knots of LARS and let~$(\,\hat \imath_1,\ldots, \hat \imath_K)$ be the first variables entering along the~LARS path.
If~$(\,\hat \imath_1,\ldots, \hat \imath_K)$ satisfies~\eqref{hyp:IrrAlongThePath}, then  
\begin{itemize}
\item conditional on $\{\hat\imath_1,\ldots, \hat \imath_K,\lambda_{K+1}\}$, the random variables~$( \lambda _1,\ldots,  \lambda_{K})$ and $\hat \sigma$ are independent;
\item conditional on $\{\hat\imath_1,\ldots, \hat \imath_K\}$ and under the null hypothesis~$\mathds H_0\,:\ \text{`}X\beta^0\in H_{K}\text{'}$, the random variables~$(\lambda_{K+1}/\hat\sigma)$ and $\hat \sigma$ are independent;
\item conditional on $\{\hat\imath_1,\ldots, \hat \imath_K,\lambda_{K+1}\}$ and under the null hypothesis~$\mathds H_0\,:\ \text{`}X\beta^0\in H_{K}\text{'}$, the distribution of~$ (\lambda _1,\ldots, \lambda_{K})$ is given by Theorem~\ref{thm:Main1}, while the distribution of~$\hat \sigma/\sigma$ is the same as the random variable
\[
(n-K)^{-\frac12}\,\Big({{\,\displaystyle\sum_{\ell=1}^{n-K}w_\ell^2}}\Big)^{\frac12}\,,
\]
where $W:=(w_1,\ldots,w_{n-K})$ is a ‘truncated' standard Gaussian vector with the truncation given by
\[
\|\mathrm{Diag}(\mathds{1}_p-\theta^K)^{-1}\times X^\top UW\|_\infty=\lambda_{K+1}/\sigma\,,
\]
where $U\in\mathds R^{n\times (n-K)}$ is any matrix such that $UU^\top=\mathrm{Id}_n-P^{(\hat\imath_1,\ldots, \hat \imath_K)}$, $\theta^K:=(\theta_{j}(\hat\imath_1,\ldots, \hat \imath_K))_{j\in[p]}$, and with the convention $0/0=0$.
\end{itemize}
\end{proposition}
  
  \begin{remark}
  \label{rem:centered}
Under the null hypothesis $\mathds H_0\,:\ \text{`}X\beta^0\in H_{K}\text{'}$, the Gaussian vectors~$P_{K}^\perp(Y)$~$($see~\eqref{e:def_sigmas}$)$, defining the variance estimate~$\hat \sigma^2$, is centered. This null hypothesis means that the true support is included in the set of the~$K$ first indices chosen by LARS.
One may choose~$K$ large enough to guarantee this null hypothesis.
  \end{remark}
 
\newpage

\noindent
Recall that, up to some positive numerical constant, the probability density function of the multivariate~$t$-distribution with~$(n-K)$ degrees of freedom, mean~$m=(m_1,\ldots,m_K)$ and variance-covariance matrix $\mathrm{Diag}(\rho_1^2,\ldots,\rho_K^2)$ is given by
\[
\tilde\varphi(t_1,\ldots,t_K):=
\Bigg[1+\frac1{n-K}{\sum_{k=1}^K \Big(\frac{t_k-m_k}{\rho_k}\Big)^2}\Bigg]^{-\frac{n}2}.
\]
We have an analogue to Theorem~\ref{thm:Main1} giving the joint law of 
\eq
\label{def_Lambda}
\Lambda_k:=\frac{\lambda_{k}}{\hat\sigma}\quad\text{for}\quad k=1,\ldots, K+1,
\qe
where~$\hat\sigma$ is given by~\eqref{e:def_sigmas} has~$n-K$ degrees of freedom, see Proposition~\ref{p:frozen}. 

\begin{theorem}[Conditional Joint Law of the~Studentized LARS knots]
\label{thm:Main2}
Let $(\lambda_1,\ldots, \lambda_K,\lambda_{K+1})$ be the first knots and let~$(\,\hat \imath_1,\ldots, \hat \imath_K)$ be the first variables entering along the~LARS path. If~$(\,\hat \imath_1,\ldots, \hat \imath_K)$ satisfies~\eqref{hyp:IrrAlongThePath} then, under the null hypothesis
\[
\mathds H_0\,:\ \text{`}X\beta^0\in H_{K}\text{'}
\] 
and conditional on the selection event~$\{\,\hat \imath_1,\ldots, \hat \imath_K,\Lambda_{K+1}\}$, the vector~$(\Lambda_{1},\ldots,\Lambda_K)$ obeys a law with the density~$($w.r.t. Lebesgue measure$)$ 
\eq
\notag
{\mathrm P}^{-1}_{(\,\hat \imath_1,\ldots,\hat \imath_K,\Lambda_{K+1})}
\tilde\varphi(t_1,\ldots,t_K)
\mathds{1}_{\{t_1\geq t_{2}\geq\cdots\geq t_{K}\geq\Lambda_{K+1}\}},
\qe
{at point }$(t_1, t_{2},\ldots,t_{K})$, where~${\mathrm P}_{(\,\hat \imath_1,\ldots,\hat \imath_K,\Lambda_{K+1})}$ is a normalizing constant,~$m_k$ and~$\rho_k$ are as in~\eqref{eq:def_mu_proj} and~\eqref{eq:def_rho}.
\end{theorem}

\begin{proof}[Proof of Theorem~\ref{thm:Main2}]
Let us fix some values~$i_1, \ldots, i_{K}$. From the definition of the Gaussian random variable~$ Z_{i_k}^{(i_1,\ldots, i_{k-1})}~$ in~\eqref{e:frozenZ}, one can deduce that its mean~$m_k$ is given by~\eqref{eq:def_mu_proj} and its standard deviation~$v_k:=\sigma\rho_k$ by~\eqref{eq:def_rho}, considering putative indices for the selected variables. By the proof of Proposition~\ref{p:frozen}, we know that these variables are independent of~$\hat \sigma^{i_1,\ldots, i_{K}}$. We deduce that the vector~$(Z_{i_1}/\hat \sigma^{i_1,\ldots, i_{K}},\ldots,Z_{i_{K}}^{(i_1,\ldots, i_{K-1})}/\hat \sigma^{i_1,\ldots, i_{K}})$ has density a multivariate~$t$-distribution with~${n-K}$ degrees of freedom, mean~$m=(m_1,\ldots,m_K)$ and variance-covariance matrix~$\mathrm{Diag}(\rho_1,\ldots,\rho_K)$. Furthermore, by~\eqref{e:indep_sigmas} of Proposition~\ref{p:frozen}, we know that this vector is independent of  $({Z^{(i_1,\ldots, i_{K})}_j}/{\hat \sigma^{i_1,\ldots, i_{K}}})_{ j\neq i_1, \ldots, i_{K}}$, and, in particular, independent of $\Lambda_{K+1}^{(i_1,\ldots, i_{K})}:=\max_j\{{Z^{(i_1,\ldots, i_{K})}_j}/{\hat \sigma^{i_1,\ldots, i_{K}}}\}$. Recall that, conditional on 
 \[
 \mathcal E:=\{\hat \imath_1 =i_1, \ldots, \hat \imath_{K}=i_{K},\lambda_{K+1}\},
 \] and assuming that~$(i_1,\ldots, i_K)$ satisfies~\eqref{hyp:IrrAlongThePath}, Proposition \ref{prop:reccur} implies that 
 \[
 \mathcal E=\Bigg\{\Lambda_{K+1}^{(i_1,\ldots, i_{K})}\leq \frac{Z_{i_{K}}^{(i_1,\ldots, i_{K-1})}}{\hat \sigma^{i_1,\ldots, i_{K}}} \leq  \cdots \leq \frac{Z_{i_1}}{\hat \sigma^{i_1,\ldots, i_{K}}}\Bigg\}\,.
 \]
Furthermore, on the event~$\mathcal E$ we have
\eq
\notag
(Z_{i_1}/\hat \sigma^{i_1,\ldots, i_{K}},  \ldots, Z_{i_{K}}^{(i_1,\ldots, i_{K-1})}/\hat \sigma^{i_1,\ldots, i_{K}},\Lambda_{K+1}^{(i_1,\ldots, i_{K})})
=
(\Lambda_1,\ldots,\Lambda_{K},\Lambda_{K+1}).
\qe
Because of the independence above, this implies that the conditional distribution is the one claimed.
\end{proof}

\newpage

\noindent
 For~$0\leq a<b\leq K+1$, we introduce 
 \begin{align}
 \notag
 \tilde\varphi_{ab}(t_{a+1},\ldots,t_{b-1})
 &:=
\Bigg[1+\frac1{n-K}{\sum_{k=a+1}^{b-1} \Big(\frac{t_k}{\rho_k}\Big)^2}\Bigg]^{-\frac{{n-K}+b-a}2}\\
\notag
\tilde\cI_{ab}(s,t)
&:=
\int_{\{s\geq t_{a+1}\geq \ldots \geq t_{b-1}\geq t\}}\tilde\varphi_{ab}(t_{a+1},\ldots,t_{b-1})\mathrm d t_{a+1} \cdots \mathrm d t_{b-1},
 \end{align}
 with the convention~$\tilde\cI_{ab}(s,t)=1$ when~$b=a+1$: and also
\begin{align}
\label{fun:tilde_fabc}
\tilde\bbF_{abc}(t)
&:=
\1_{\{\Lambda_c\leq t\leq \Lambda_a\}}
\displaystyle
\int\displaylimits_{\Lambda_{c}}^{t}
\tilde\cI_{ab}({\Lambda_a},\ell_b)\,\tilde\cI_{bc}(\ell_b,{\Lambda_c})\,
\Bigg[1+\frac1{n-K}{\Big(\frac{\ell_b}{\rho_b}\Big)^2}\Bigg]^{-\frac{{n-K}+1}2}
\mathrm d\ell_b\\
\notag
&\text{for }0\leq a<b<c\leq K+1,\ t\in\bbR.
\end{align}
When~$m_{a+1}=\cdots=m_{c-1}=0$, the function~$\tilde\bbF_{abc}$ gives the CDF of~$\Lambda_b$ conditional on~$\Lambda_a,\Lambda_c$ and on some selection event, as shown below in Theorem~\ref{t:pvalue} and~\eqref{t_CDF}. For~$0\leq a<b<c\leq K+1$, we introduce the~$p$-value 
 \eq
\label{e:talpha}
 \hat  \beta_{abc} = \hat  \beta_{abc}( {\Lambda_a},{ \Lambda_b},{ \Lambda_c}, \hat \imath_1, \ldots, \hat \imath_{K} )= 1 - 
  \frac{\tilde\bbF_{abc}(\Lambda_b)}
 {\tilde\bbF_{abc}(\Lambda_a)}
\, 
\qe
On the numerical side, note that this quantity can be computed using \textit{Quasi Monte Carlo} (QMC) methods as in~\cite[Chapter 5.1]{genz2009computation}.

\subsubsection{Proof of Theorem~\ref{t:pvalue}}
\label{proof:t:pvalue}
Fix~$a$ such that~$0\leq a\leq K-1$ and consider any selection procedure~$\hat m$ satisfying~\eqref{e:stopping_rule}.
 From Proposition~\ref{prop:H0}, conditional on 
\[
 \mathcal F:=\big\{\hat \imath_1 =i_1, \ldots, \hat \imath_{K}=i_{K}, \Lambda_a,\Lambda_{K+1}\big\}\,
 \]
and under the null hypothesis~$\mathds H_0\,:\ \text{`}X\beta^0\in H_{a}\,\text{'}$, we know that~$m_{a+1}=\ldots=m_K=0$. From Theorem~\ref{thm:Main2} we know that the density of~$(\Lambda_{a+1},\Lambda_{a+2},\ldots,\Lambda_K)$ conditional on~$\mathcal F$ is given by 
    \eq
    \notag
    (const)  
    \Bigg[1+\frac1{n-K}{\sum_{k=a+1}^K \Big(\frac{t_k}{\rho_k}\Big)^2}\Bigg]^{-\frac{{n}-a}2}
    \, \1_{\Lambda_a \geq \ell_{a+1}\geq \cdots \geq \ell_K\geq  \Lambda_{K+1}}.
   \qe

\noindent From the definition of assumption~\eqref{e:stopping_rule}, and on the event~$\mathcal F$, we know that the indicator~$\1_{\{\hat m=a\}}$ is a measurable function of~$\lambda_1,\ldots, \lambda_{a-1}$, which are respectively equal to $Z_{i_1},\ldots,  Z_{i_{a-1}}^{(i_1,\ldots, i_{a-2})}$ on~$\mathcal F$ by~\eqref{e:frozen}. By \eqref{e:indep_sigmas} of Proposition~\ref{prop:reccur}, we deduce that $\1_{\{\hat m=a\}}$ is independent of~$(\lambda_{a+1}/\hat \sigma^{i_1,\ldots, i_{K}},\ldots,\lambda_K/\hat \sigma^{i_1,\ldots, i_{K}})$ and of $\Lambda_{K+1}:=\lambda_{K+1}/\hat \sigma^{i_1,\ldots, i_{K}}$ conditional on~$\mathcal F$. We deduce that the conditional density above is also the conditional density on the event 
 \[
 \mathcal G:=\big\{\hat m =a,\hat \imath_1 =i_1, \ldots, \hat \imath_{K}=i_{K}, \Lambda_a,\Lambda_{K+1}\big\}.
 \]
Now, a simple integration shows that  
\eq
\label{t_CDF}
\bbP\big[\Lambda_b\leq t\ |\ \hat m =a,\Lambda_a,\Lambda_c, \hat \imath_{1},\ldots,\hat \imath_{K}\big]=\frac{\tilde\bbF_{abc}(t)}{\tilde\bbF_{abc}(\Lambda_a)}.
\qe
As a consequence and under the same conditioning, one has
  \[
  \frac{
  \tilde\bbF_{abc}(\Lambda_b)    }{    \tilde\bbF_{abc}(\Lambda_a)  }  \sim  \mathcal  U(0,1).
  \]
  Finally, considerations of the distribution under the alternative  show that  to obtain a~$p$-value we must consider the complement to~$1$ of  the quantity above.

 \subsubsection{$t$-Spacing tests and Generalized $t$-Spacing tests}

Consider the following testing procedures:
\eq
\label{zaza66}
\cT_{abc}:=\mathds1_{\{\hat\beta_{abc}\leq\alpha\}},
\qe
that rejects if the~$p$-value~$\hat\beta_{abc}$ is less than the level~$\alpha$ of the test.
This test statistic generalizes previous test statistics that appeared in~$t$-Spacing Tests, as presented in~\cite{azais2018power} for instance, and will be referred to as the {\it Generalized~$t$-Spacing test} (GtSt).

\begin{remark}
\label{rem:tST}
If one takes~$a=0$,~$b=1$ and~$c=2$, then one gets
\[
\hat\beta_{012}=1-\frac{\bT_{1}(\Lambda_1)-\bT_{1}(\Lambda_2)}{\bT_{1}(\Lambda_0)-\bT_{1}(\Lambda_2)}
=\frac{1-\bT_{1}(\Lambda_1)}{1-\bT_{1}(\Lambda_2)}.
\]
Similarly, taking~$b=a+1$ and~$c=a+2$, one gets
\[
\hat\beta_{a(a+1)(a+2)}
=\frac
{\bT_{{a+1}}(\Lambda_{a+1})-\bT_{{a+1}}(\Lambda_a)}
{\bT_{{a+1}}(\Lambda_{a+2})-\bT_{{a+1}}(\Lambda_a)}.
\]
which is the {\it $t$-spacing test} as presented in~\cite{azais2018power}, where 
\eq
\label{e:Student}
\bT_k(\ell):=\int\displaylimits_{-\infty}^\ell\Bigg[1+\frac1{n-K}{\Big(\frac{\ell}{\rho_k}\Big)^2}\Bigg]^{-\frac{{n-K}+1}2}\mathrm d\ell
\qe
is, up to some positive numerical constant, the CDF of a centered~$t$-Student distribution with variance~$\rho_k^2$ and~$n-K$ degrees of freedom.
\end{remark}

\subsection{Power studies} \label{subsec:ortho}

\subsubsection{Power when the design is orthogonal}
One may investigate the power of these tests at detecting false negatives, namely, the alternatives given by: there exists~$k\in S^0$ such that~$k\notin\{\bar\imath_1,\ldots,\bar\imath_a\}$.
In particular, what is the most powerful test among these latter~\eqref{zaza6} testing procedures? A comprehensive study for the case of orthogonal designs is given by Theorem \ref{thm:zazPower}.

\subsubsection{Numerical studies on the power for the general design case}
\label{sec:num_studies}
In the case of an orthogonal design, Theorem~\ref{thm:zazPower} shows that the test based on~$\hat\alpha_{a,a+1,K+1}$ is uniformly more powerful than tests based on~$\hat\alpha_{x,y,z}$ with~$a\leq x<y<z\leq K+1$.
Numerical experiments on the power of these tests are presented in Figure~\ref{fig:laws} and they witness the same phenomenon for Gaussian designs.
It presents the CDF of the~$p$-value~$\hat \alpha_{abc}$ under the null and under two~$2$-sparse alternatives, one with low signal and one with~$5$ times more signal.
The numerical results show that all the tests are exact (leftmost panel) and the test~$\mathcal S_{125}$ is the most powerful.
A detailed presentation of this is given in Section~\ref{sec_power_num}.

\begin{figure}[!t]
\includegraphics[width=0.3\textwidth]{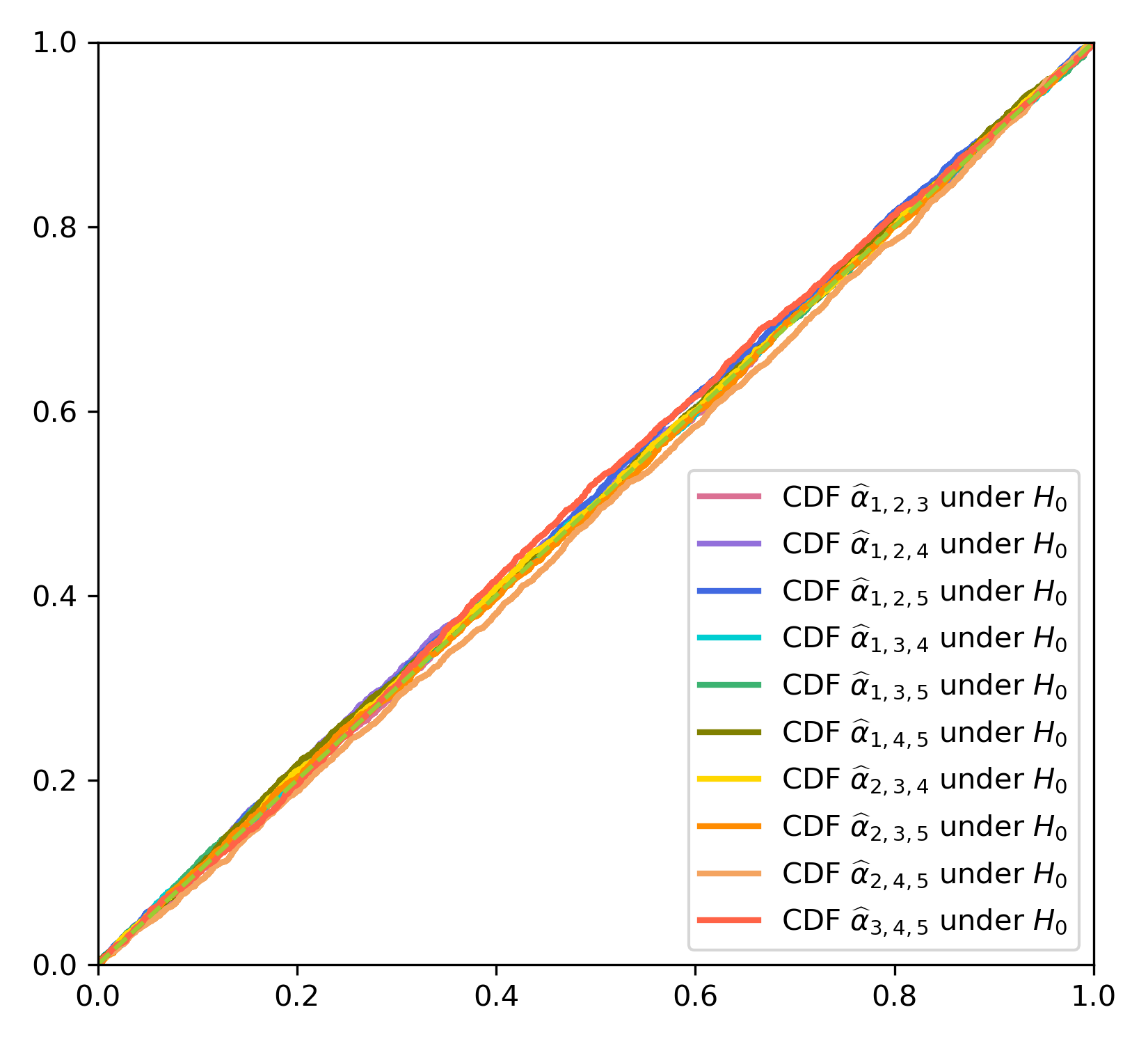}
\includegraphics[width=0.3\textwidth]{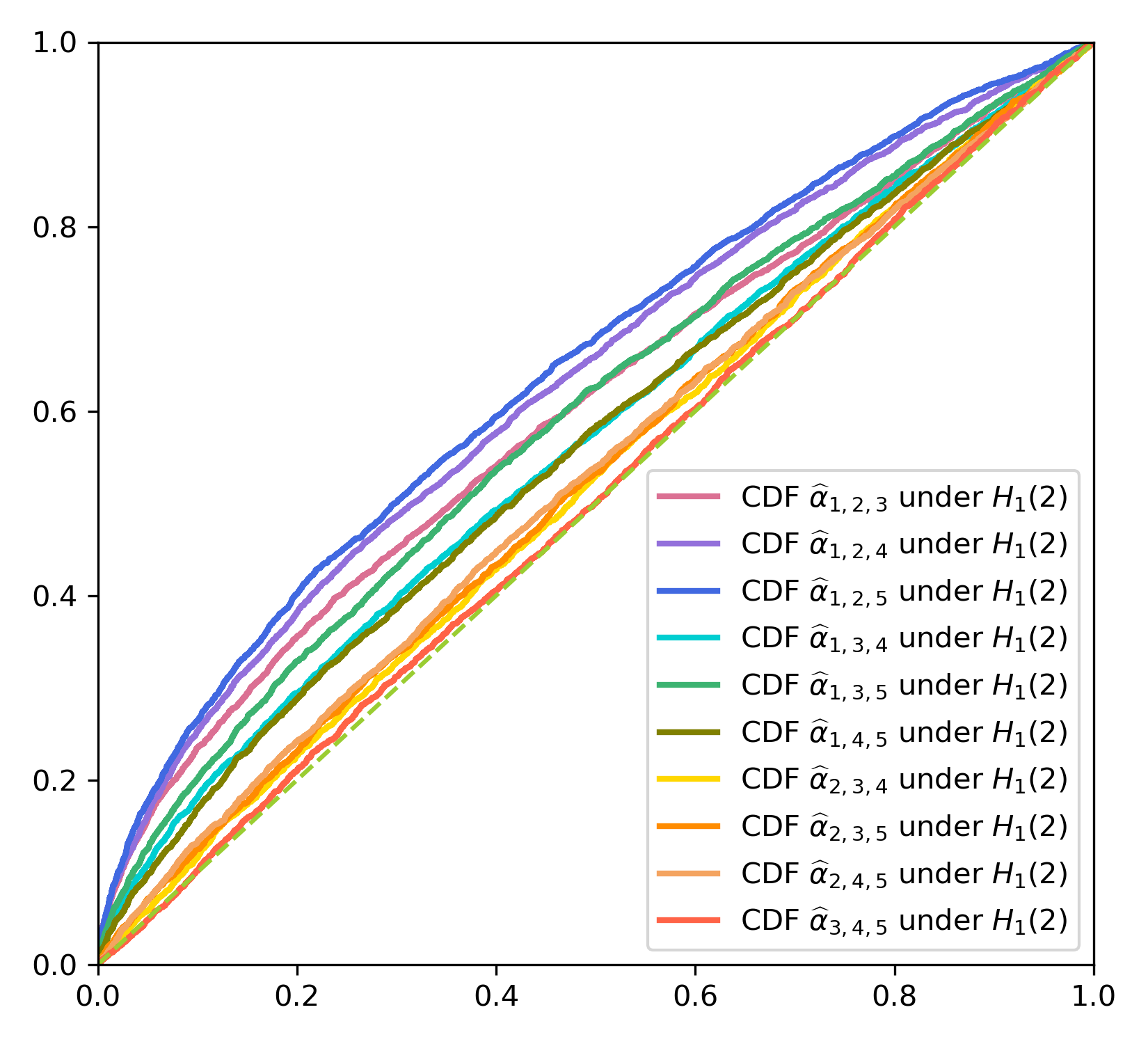}
\includegraphics[width=0.3\textwidth]{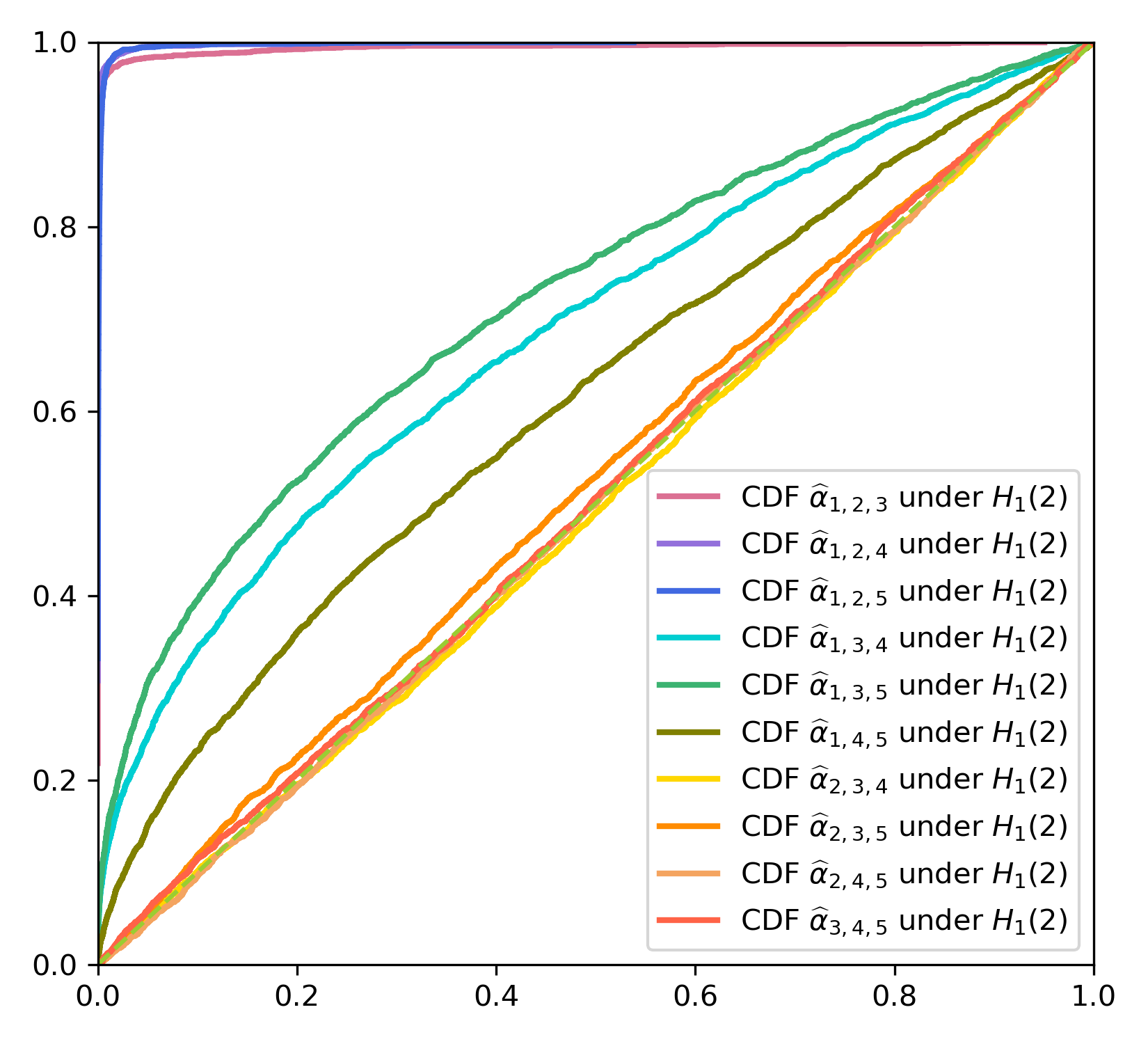}
\caption{CDF of~$p$-values~$\hat\alpha_{abc}$ over~$3,000$ Monte Carlo iterations and a random design~$X\in\mathds R^{n\times p}$ given by~$p=300$ independent column vectors uniformly distributed on the Euclidean sphere~$\mathds S^{199}$ ($n=200$).
Central panel represents alternative composed by~$2$-sparse vector, right panel alternative composed by~$2$-sparse vector~$5$ times larger while left panel corresponds to the null.}
\label{fig:laws}
\end{figure}

\newpage

\noindent
More precisely, one has, as proved in the orthogonal case by Theorem~\ref{thm:zazPower} (and its proof),  that 
\begin{itemize}
\item~$\hat \alpha_{125} \preccurlyeq \hat \alpha_{124}\preccurlyeq \hat\alpha_{123}$;
\item~$\hat\alpha_{125} \preccurlyeq  \hat\alpha_{135} \preccurlyeq \hat\alpha_{235} \preccurlyeq \hat\alpha_{234}$;
\item~$\hat\alpha_{125} \preccurlyeq \hat\alpha_{135} \preccurlyeq \hat\alpha_{145}\preccurlyeq \hat\alpha_{245}$.
\end{itemize}
where~$\preccurlyeq$ denotes stochastic ordering.
In the proof of Theorem~\ref{thm:zazPower}, it was shown that 
\[
\hat \alpha_{ab(c+1)}\preccurlyeq \hat \alpha_{abc}\text{ and }
\hat \alpha_{a(b-1)c}\preccurlyeq \hat \alpha_{abc}\text{ and }
\hat \alpha_{(a-1)bc}\preccurlyeq \hat \alpha_{abc},
\]
for orthogonal designs.

\subsection{Control of False Discovery Rate in the Orthogonal Design case}
\label{sec:FDR}

\subsubsection{Presentation in the general case}
\label{hyp:FDR} 

 For the sake of readability, we will assume, for the moment,  that~$\sigma$ is known. We understand that the law of test statistics are parametrized by the hypotheses~$(m_k)_{k\in[K]}$, where~$m_k$ is given by~\eqref{eq:def_mu_proj}.

We recall that we write~$\bar\mu^0=X^\top X \beta^0$ and~$\bar\mu^0_i$ for its~$i$th coordinate.
Assuming that the predictors are normalised, in the general case, this quantity is the sum of~$\beta_i^0$ and a linear combination of the~$\beta_j^0$'s whose predictors~$X_j$ are highly correlated with the predictor~$X_i$.
Now, given the variables~$\bar\imath_1,\ldots,\bar\imath_{k}\in[p]$ and signs~$\varepsilon_{1},\ldots,\varepsilon_{{k}}\in\{\pm1\}^{k}$, we denote by~$(\Pi_{ \bar \imath_1,\ldots, \bar \imath_{k-1}}^\perp (\bar\mu^0))_{\bar\imath_{k}}$ the orthogonal projection given by
\eq
\label{eq:ortho_mu}
(\Pi_{ \bar \imath_1,\ldots, \bar \imath_{k-1}}^\perp (\bar\mu^0))_{\bar\imath_{k}}
:=\varepsilon_{{k}}X_{\bar \imath_k}^\top \Big[\mathrm{Id}_{n}-X_{\bar S^{k-1}}\big(X_{\bar S^{k-1}}^\top X_{\bar S^{k-1}}\big)^{-1}\!\!X_{\bar S^{k-1}}^\top\Big] X\beta^0.
\qe

The tested null hypotheses are conditional on  some sub-sequence of variables $(\bar \imath_1,\ldots,\bar \imath_{K+1})\in[p]^{K+1}$ and signs~$\varepsilon_{1},\ldots,\varepsilon_{{K+1}}\in\{\pm1\}^{K+1}$ entering the model.
The~$p$-values under consideration are 
\begin{align}
&\circ\ \hat p_1:=\hat\alpha_{0,1,2} \mbox{ is the }p\mbox{-value testing }\mathds H_{0,1}: ‘‘m_1=0\," \mbox{ namely } \bar\mu^0_{\bar \imath_1}=0 \notag ; \\
&\circ\ \hat p_2:=\hat\alpha_{1,2,3} \mbox{ is the }p\mbox{-value testing }\mathds H_{0,2}: ‘‘m_2=0\," \mbox{ namely } (\Pi_{{1}}^\perp (\bar\mu^0))_{\bar\imath_2}=0 \notag ; \\
\label{e:pval}
&\circ\ \hat p_3:=\hat\alpha_{2,3,4} \mbox{ is the }p\mbox{-value testing }\mathds H_{0,3}: ‘‘m_3=0\," \mbox{ namely }  (\Pi_{{2}}^\perp (\bar\mu^0))_{\bar\imath_3}=0 ;\\
&\circ\ \mbox{and so on...} \notag
\end{align}

We write~$I_0$ of the set $I_0=\big\{k\in[K]\ :\ \mathds H_{0,k}\text{ is true}\big\}$. Given a subset~$\hat R\subseteq [K]$ of hypotheses that we consider as rejected, we call {\it false positive} ($\mathrm{FP}$) and {\it true positive} ($\mathrm{TP}$) the quantities~$\mathrm{FP}=\mathrm{card}(\,\hat R\cap I_0)$ and~$\mathrm{TP}=\mathrm{card}(\,\hat R\setminus I_0)$. Denote by~$\hat p_{(1)}\leq\ldots\leq \hat p_{(K)}$ the~$p$-values ranked in a nondecreasing order.
Let~$\alpha\in(0,1)$ and consider the Benjamini--Hochberg procedure, see for instance~\cite{benjamini1995controlling}, defined by a rejection set~$\hat R\subseteq[K]$ such that~$\hat R=\emptyset$ when~$\{k\in[K]\ :\ \hat p_{(k)}\leq\alpha k/K\}=\emptyset$ and 
\begin{equation}\label{e:rchap}
\hat R=\{k\in[K]\ :\ \hat p_k\leq \alpha\hat k/K\}\quad\mathrm{where}\quad \hat k=\max\big\{k\in[K]\ :\ \hat p_{(k)}\leq\alpha k/K\big\}.
\end{equation}
Recall the definition of the FDR as the mean of the False Discovery Proportion ($\mathrm{FDP}$), namely
\[
\mathrm{FDR}:=\mathds E\Big[\underbrace{\frac{\mathrm{FP}}{\mathrm{FP}+\mathrm{TP}}\mathds 1_{\mathrm{FP}+\mathrm{TP}\geq1}}_{\mathrm{FDP}}\Big],
\]
where the expectation is unconditional on the sequence of variables entering the model, while the hypotheses that are being tested are conditional on the  sequence of variables entering the model.
This FDR can be understood by invoking the following decomposition
\[
\mathrm{FDR}=\sum_{(\imath_1,\ldots,\imath_K)\in[p]^K}
\!\!\!\!\!\!
\bar\pi_{ (\imath_1,\ldots, \imath_K)}\,\mathds E\big[\mathrm{FDP}|\bar\imath_1=i_1,\ldots,\bar\imath_K=i_K\big],
\]
where~$\bar\pi_{ (\imath_1,\ldots, \imath_K)}=\mathds P\big\{\bar \imath_1=\imath_1,\ldots, \bar \imath_K=\imath_K\big\}$.

\subsubsection{Control of the FDR by the Benjamini–Hochberg procedure in the orthogonal design case}

We now consider the case of an orthogonal design where~$X^\top X=\mathrm{Id}_p$ and the set of~$p$-values
is given by~\eqref{e:pval}.
Note that~$I_0$  is simply the set  of null coordinates of~$\beta^0$.
Remark  also that the Irrepresentable Condition~\eqref{eq:Irr_matrix} of order~$p$ holds and so does Empirical Irrepresentable Check~\eqref{hyp:IrrAlongThePath}, see Proposition~\ref{prop:Irrepresentable}.

\begin{theorem}
 \label{thm:independentTest}
 Assume that the design is orthogonal, i.e.~$X^\top X=\mathrm{Id}_p$, and let~$K\in[p]$.
Let~$(\bar \imath_1,\ldots, \bar \imath_K)$ be the first variables entering along the~LARS path.
Consider the~$p$-values given by~\eqref{e:pval} and the set~$\hat R$ given by~\eqref{e:rchap}.
Then
\[
\mathds E\big[\mathrm{FDP}|\bar\imath_1=i_1,\ldots,\bar\imath_K=i_K\big]\leq\alpha,
\]
and so~$\mathrm{FDR}$ is bounded above by~$\alpha$.
 \end{theorem}
 \noindent
The proof of this result is given in Appendix \ref{proof:yoFDR}. One interpretation of post-selection type may be given as follows: if one looks at all the experiments giving the same sequence of variables entering the model~$\{\bar\imath_1=i_1,\ldots,\bar\imath_K=i_K\}$ and if one considers the Benjamini–Hochberg procedure for the hypotheses described in Section~\ref{hyp:FDR}, then the FDR is exactly controlled by~$\alpha$.

\section{Testing procedures: Numerical studies}
\label{sec:num}
\subsection{Power in the non-orthogonal case}
\label{sec_power_num}
To study the power in the case of a {\it non-orthogonal} design, we built a Monte-Carlo experiment with:
\begin{itemize}
\item a model with~$n=200$ observations and~$p=300$ predictors,
\item a random design matrix~$X$ given by~$300$ independent column vectors uniformly distributed on the Euclidean sphere~$\mathds S^{199}$,
 \item and we ran~$3,000$ Monte Carlo experiments.
 \item The results are presented in Figure~\ref{fig:laws}.
\end{itemize}
The computation of the function~$\bbF_{abc}$ given by~\eqref{fun:fabc} requires multivariate integration tools. All our test statistics can be efficiently computed using Quasi Monte Carlo methods (QMC) for Multi-Variate Normal (MVN) and~$t$ (MVT) distributions, see the book~\cite{genz2009computation} for a comprehensive treatment of this topic or Appendix~\ref{app:cbc} for a short overview of the method we used. We compute spacings of length at most~$4$, which implies that~$c\leq 5$ when~$a=1$ in our experimental framework. 

A Python notebook and codes are given at \url{https://github.com/ydecastro/lar_testing}. The base function is 
\begin{center}
{\tt observed\_significance\_CBC(lars, sigma, start, end, middle)}
\end{center}
in the file {\tt multiple\_spacing\_tests.py}. It gives the~$p$-value~$\hat\alpha_{(\mathtt{start})(\mathtt{middle})(\mathtt{end})}$ of the knots and indices given by {\tt lars} and an estimate of (or the true) standard deviation given by {\tt sigma}. We ran~$3,000$ repetitions of this function to get the laws displayed in Figure~\ref{fig:laws}. It presents the CDF of the~$p$-value~$\hat \alpha_{abc}$ under the null and under two~$2$-sparse alternatives, one with low signal and one with~$5$ times more signal. The results show, in our particular case, that all the tests are exact and the test~$\mathcal S_{125}$ is the most powerful, see Section~\ref{sec:num_studies} for further details.

\subsection{{A comparison of FDR control and power on simulated data}}
\label{sec:simulated}
We take the experiments introduced in \cite[Section 5]{javanmard2019false}. As in this reference, we consider a linear model with design $X$ with independent rows drawn with respect to $\mathcal N_p(0,\Sigma)$. The covariance $\Sigma\in\mathds R^{p\times p}$ is such that $\Sigma_{ij} = r^{|i-j|}$, for some parameter $r\in(0,1)$. We then normalize the columns of~$X$ to have unit Euclidean norm. We draw a $k$-sparse vector $\beta^0\in \mathds R^p$ by choosing a support of size $k$ at random with values $\{\pm A\}$ uniformly at random, where $A>0$ denotes the absolute value of the amplitudes. The Gaussian noise term $\eta$ is drawn from $\mathcal N_n(0,\mathrm{Id}_{n})$.

\noindent
We compare the performances of three procedures:
\begin{itemize}
\item {\bf [Knockoff]} Knockoff filters for FDR control \citep{barber2015controlling} and we use {\tt knockoff+} as implemented on \url{https://web.stanford.edu/group/candes/knockoffs/}; 
\item {\bf [FCD]} False Discovery Control via Debiasing \cite[Section 5]{javanmard2019false} and we use the implementation of debiased lasso presented on the webpage \url{https://web.stanford.edu/~montanar/sslasso/} with the theoretical value $\bar\lambda=2\sqrt{(2\log p)/n}$ for the regularizing parameter. When the sample size is larger than the number of predictors ($n\geq p$), the debiasing step in FCD is superfluous as the decorrelating matrix ($M$) can be the inverse of the sample covariance. So, in this case, we start with an unbiased estimator upfront (which is Ordinary Least Squares OLS). The FCD then becomes thresholding the test statistics $|T_i|$ obtained from OLS;
\item {\bf [GtSt-BH]} Generalized $t$-Spacing tests on successive entries of the LARS path combined with a Benjamini–Hochberg procedure \cite{benjamini1995controlling} based on the sequence of spacings $\hat  \beta_{012},\hat  \beta_{123},\ldots,\hat  \beta_{a(a+1)(a+2)},\ldots$ 
as described in Section  \ref{sec:Student} with nominal value~$\alpha= 0.1$;
\end{itemize}
we numerically investigate the effects of the level of sparsity, the magnitude of the signal, the correlation between the features, and the empirical power. In all simulations, we set the target level FDR to $\alpha= 0.1$.

\subsubsection{The effect of the amplitude of the signal}
We chose $n = 200$ (sample size), $p = 100$ (predictors), $k = 20$ (sparsity), $\eta = 0.1$ (features correlation) and varied the amplitude within the set $A\in [25]$. We computed the FDR and power by averaging over $3,000$ realizations of the noise and generations of the coefficients of the vector $\beta^0$. The results are plotted in Figure~\ref{f:effect_ampl}. Recall that, in the case $n\geq p$, FCD is a thresholded OLS and it might be considered as the best test here. One may note that it presents the best features (low FDR and high power). GtSt controls the FDR below the nominal value~$\alpha= 0.1$ with a slightly lower power than FCD. Knockoff+ has controlled FDR and matches the power of FCD.
\begin{figure}[!h]
\vspace*{-0.43cm}
\begin{center}
\includegraphics[width=0.42\textwidth]{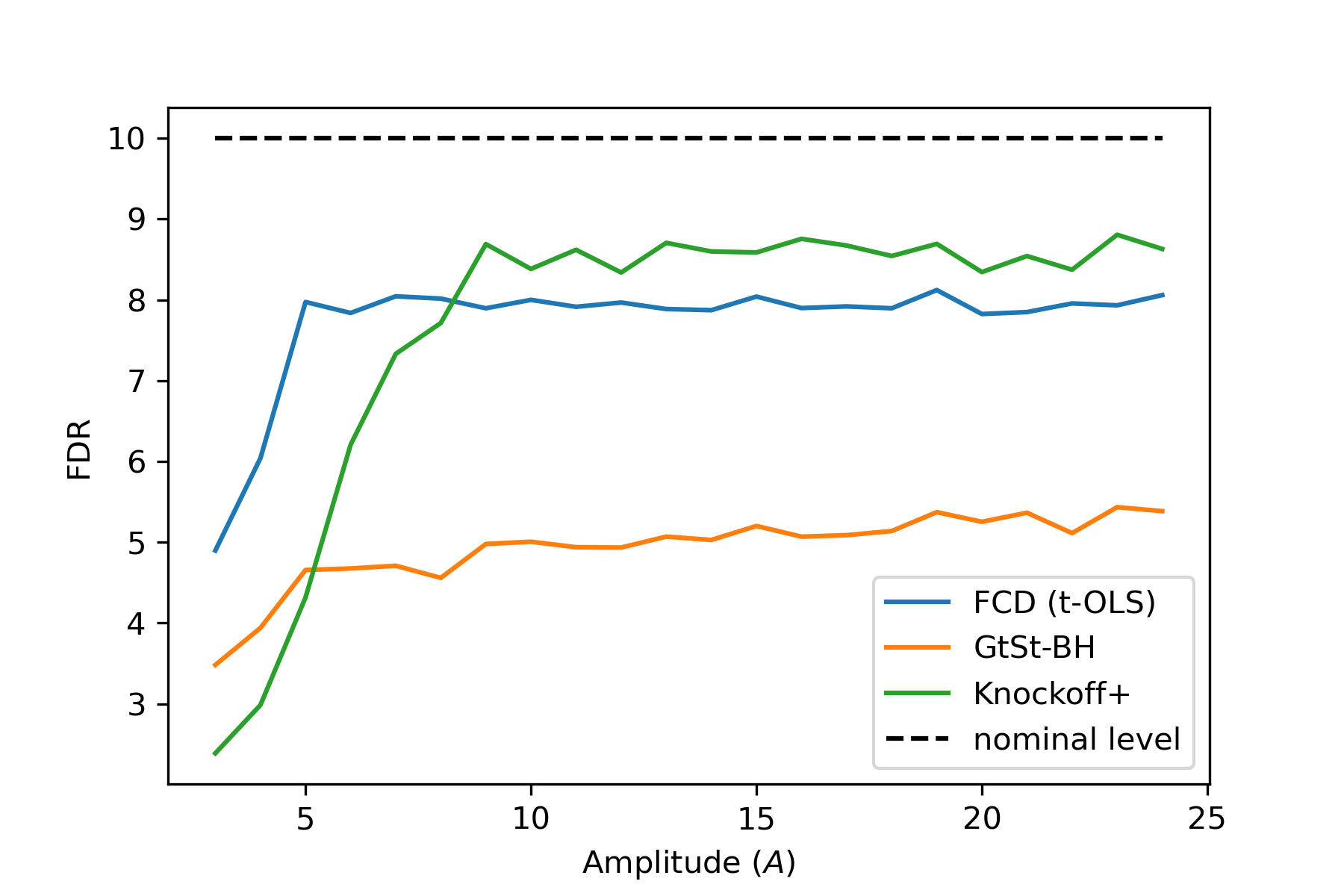}
\includegraphics[width=0.42\textwidth]{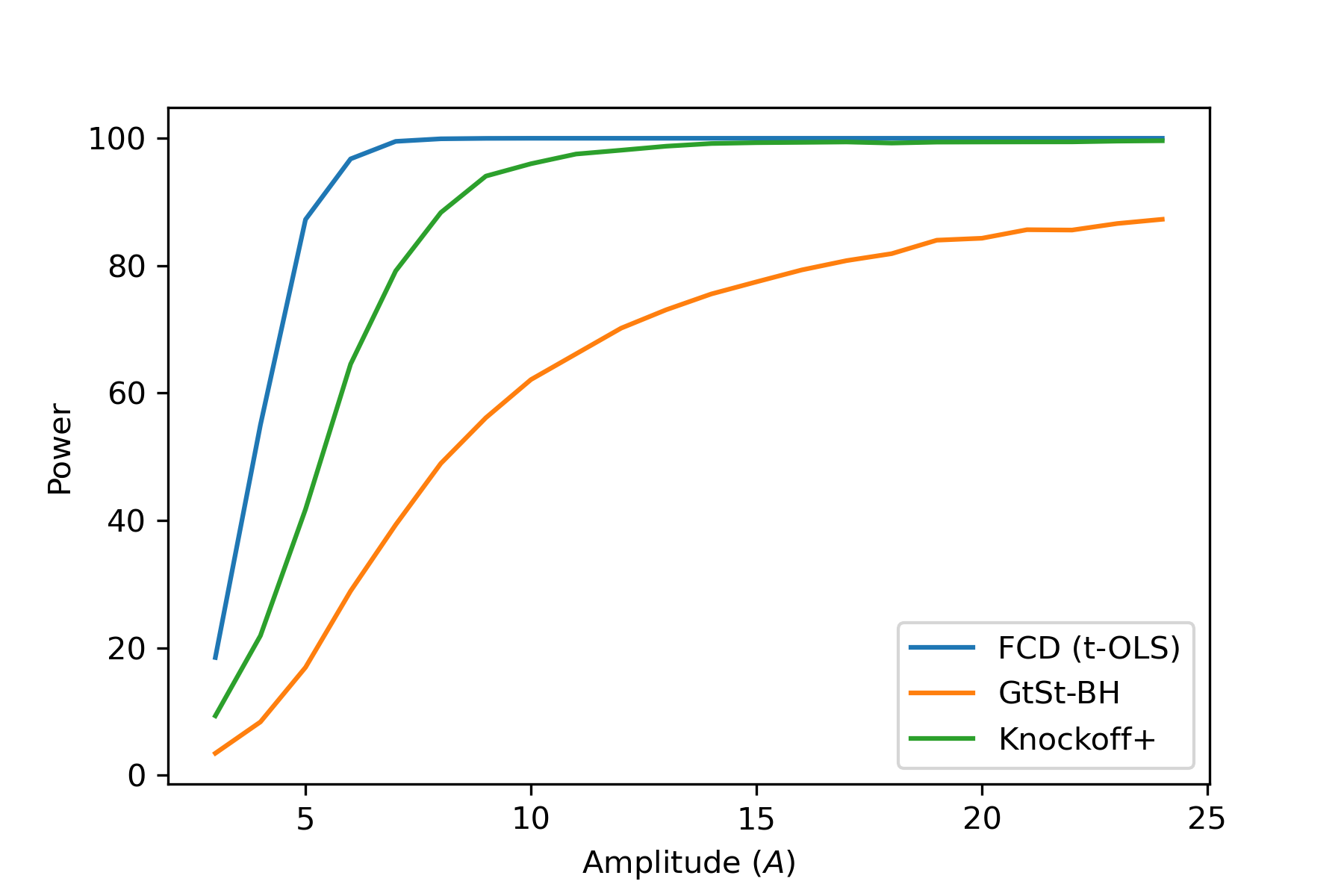}
\end{center}
\caption{Comparison of the FDR  control and power for GtSt, FCD and Knockoff+ when the amplitude of the signal varies from $1$ (low signal) to $25$ (strong signal) over $3,000$ trials.} \label{f:effect_ampl}
\end{figure}

\subsubsection{Effect of feature correlation}
We test the effect of correlations between the features with $n =200$, $p = 100$, $k = 20$, and $A = 10$. Recall that the rows of the design matrix $X$ are generated from an~$\mathcal N_p(0,\Sigma)$ distribution, with $\Sigma_{ij} = \eta^{|i-j|}$, and then the columns of $X$ are normalized to have unit norm. We vary the  parameter~$\eta$ within the set $\{0.1, 0.15, 0.2, \dotsc, 0.75, 0.8, 0.85\}$. For each value of $\eta$, we compute the FDR and power by averaging over $3,000$ realizations of the noise and design matrix $X$. The results are displayed in Figure~\ref{f:effect_corr}. One may note that, in the case $n\geq p$, FCD is a thresholded OLS and might be considered as the best estimation here (with low FDR and high power). Knockoff+ has controlled FDR and matches the power of FCD for small feature correlations. GtSt controls the FDR below the nominal value~$\alpha= 0.1$ with a lower power than Knockoff+. 

\begin{figure}[!h]
\vspace*{-0.43cm}
\begin{center}
\includegraphics[width=0.42\textwidth]{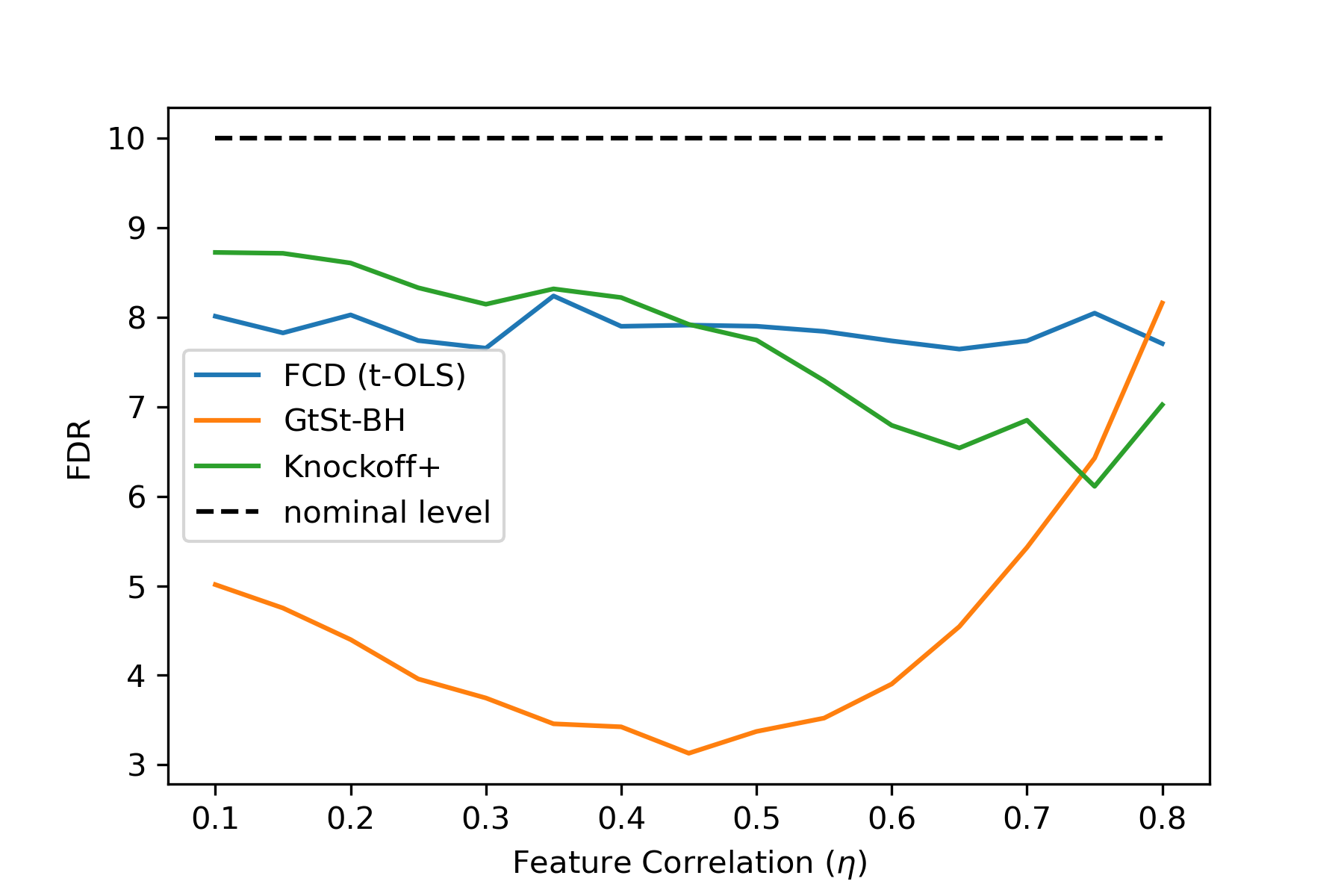}
\includegraphics[width=0.42\textwidth]{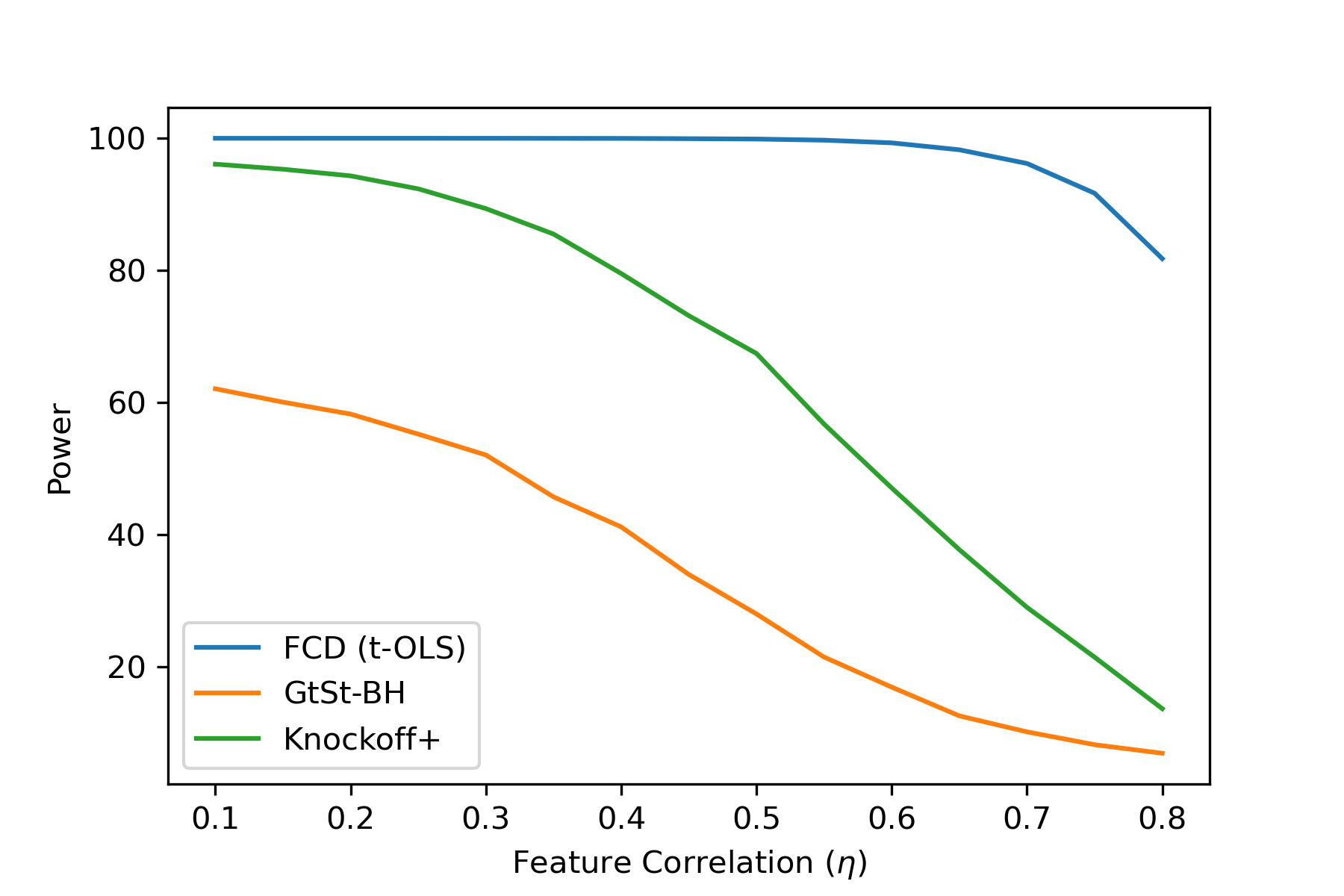}
\end{center}
\caption{Comparison of the FDR  control and power for GtSt, FCD and Knockoff when the correlation between features varies from $0.1$ (low correlation) to $0.85$ (strong correlation) over $3,000$ trials.} \label{f:effect_corr}
\end{figure}

\subsubsection{Effect of sparsity}
We set $n = 200$, $p = 100$, $A = 10$, and $\eta = 0.1$, and varied the level of sparsity of the coefficients within the set $k\in(10,40)$. The power and the FDR are computed by averaging over $3,000$ trials of the noise and generations of the coefficients of the vector $\beta^0$. The results are displayed in Figure~\ref{f:effect_sparsity}. Knockoff+ has controlled FDR and matches the power of FCD. GtSt controls the FDR below the nominal value~$\alpha= 0.1$ with a lower power than Knockoff+. 

\begin{figure}[!h]
\vspace*{-0.43cm}
\begin{center}
\includegraphics[width=0.42\textwidth]{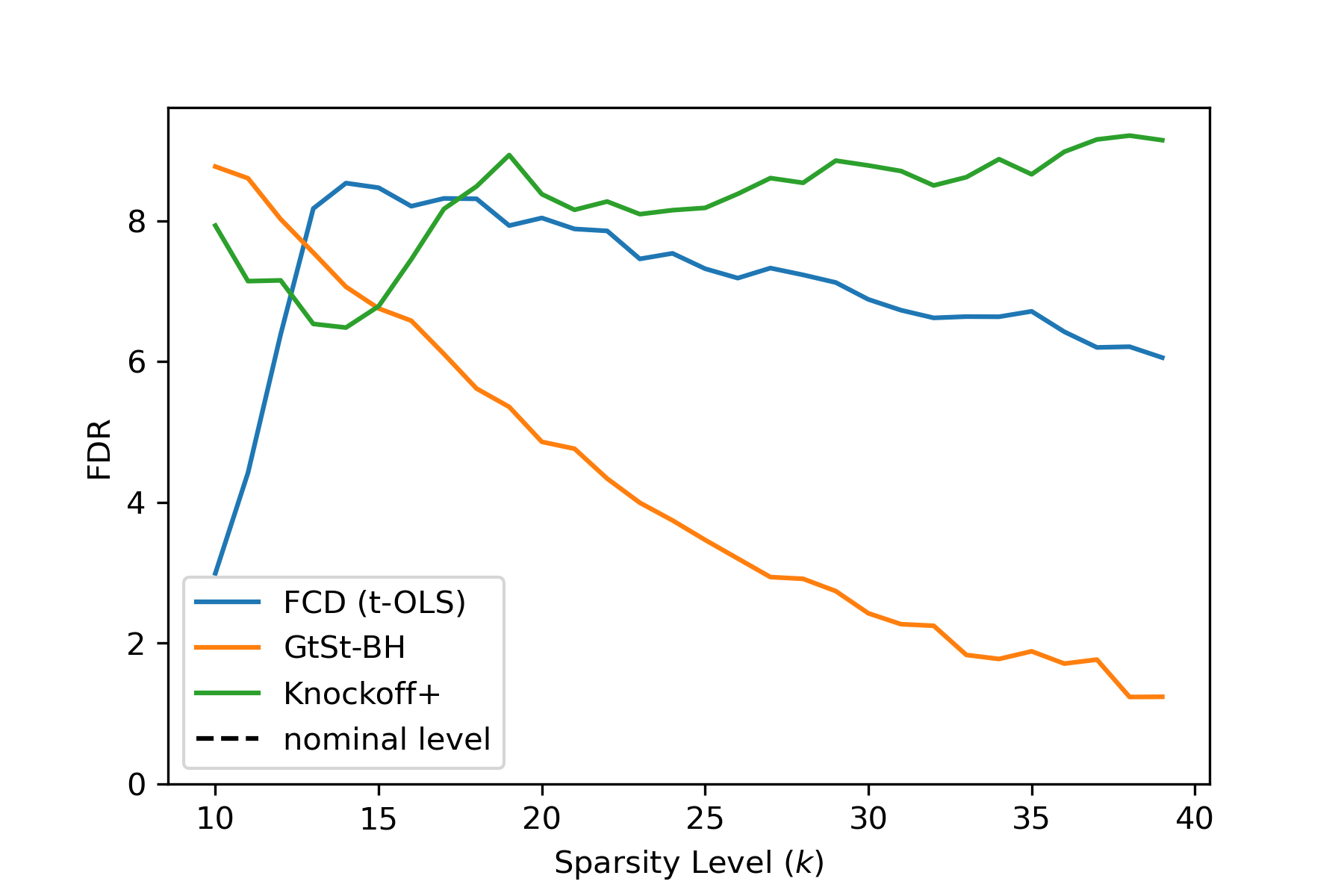}
\includegraphics[width=0.42\textwidth]{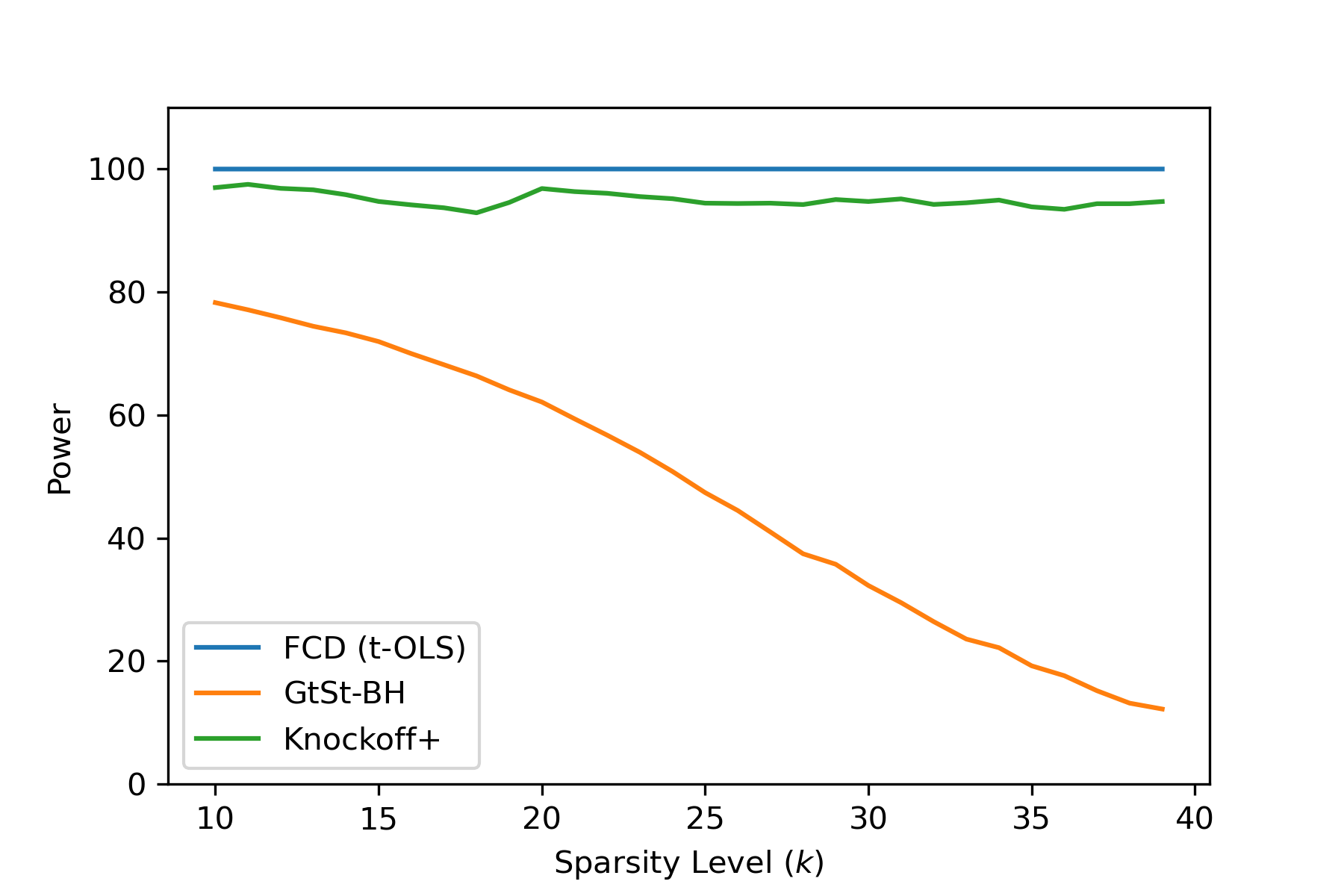}
\end{center}
\caption{Comparison of the FDR  control and power for GtSt, FCD and Knockoff when the feature correlation varies from $1$ (very sparse) to $35$ (slightly sparse) over $3,000$ trials.} \label{f:effect_sparsity}
\end{figure}

\subsection{FDR on real data}
\label{realdata}

A detailed presentation in a Python notebook is available at \url{https://github.com/ydecastro/lar_testing/blob/master/multiple_spacing_tests.ipynb}. We consider a data set  about HIV drug resistance extracted from~\cite{barber2015controlling}  and~\cite{rhee2006genotypic}. The experiment consists in identifying mutations of the genes of the HIV that are involved with drug resistance. The data set  contains  about~$p=200$ and~$n=700$ observations.  Since some protocol was used to remove some genes and some individuals, the exact numbers depend  on the considered drug.  

\begin{figure}[h]
\vspace*{-0.1cm}
\begin{center}
\includegraphics[width=0.99\textwidth]{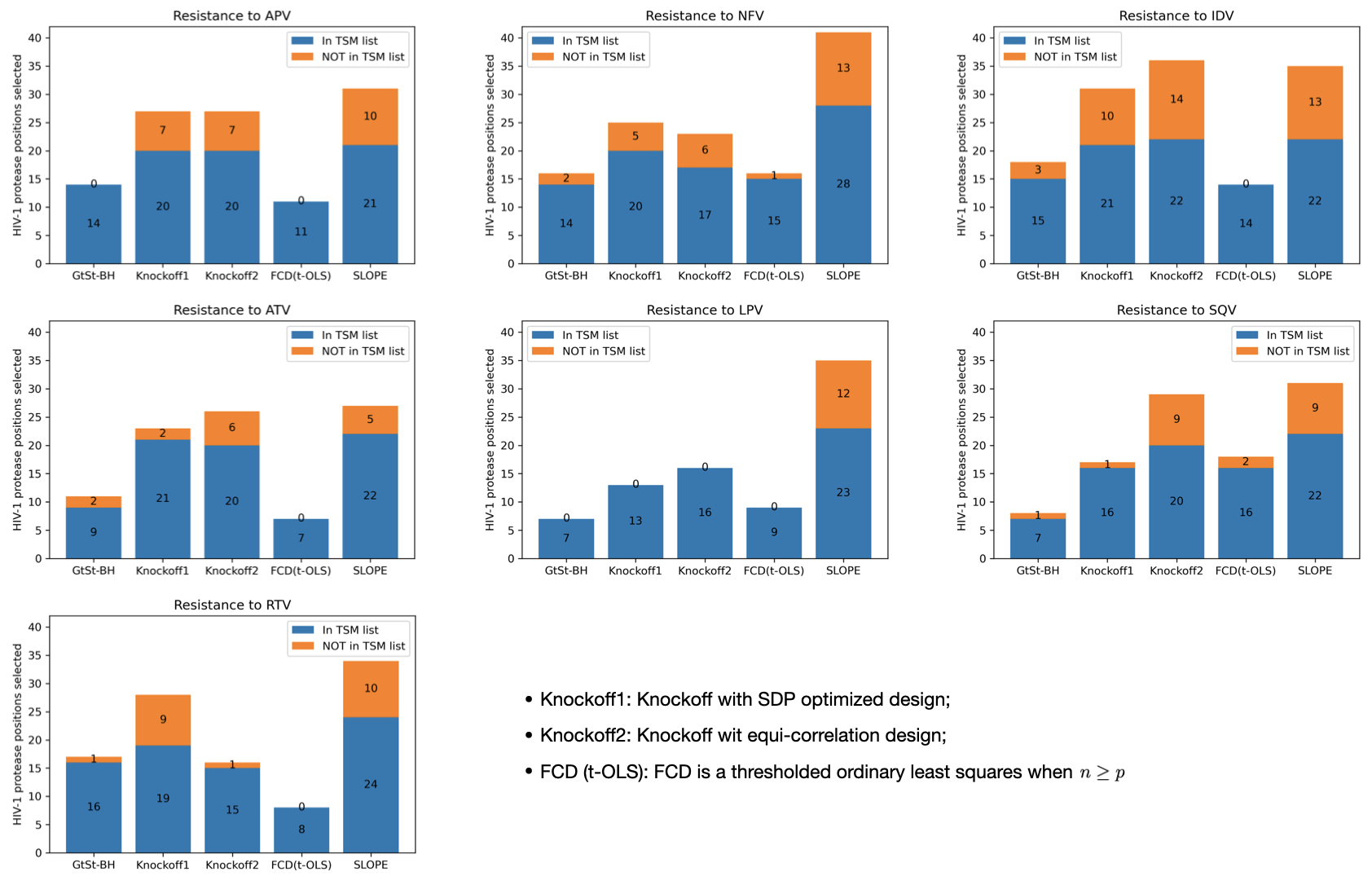}
\end{center}
\caption{Comparison of  the number of true and false positives for procedures:  GtSt-BH, Knockoff (with SDP optimisation designs in Knockoff1 and equi-correlation designs in Knockoff2),
FDC and Slope. In the procedures, the FDR aimed at is~$\alpha=20\%$. Blue indicates protease positions that appear in the TSM panel for the PI class of treatments, given in \cite[Tabel 1]{rhee2006genotypic}, while orange indicates positions selected by the method that do not appear in the TSM list. The total number of HIV-1 protease positions appearing in the TSM list is $34$.} \label{f:candes}
\end{figure}

\noindent
The methods considered are {\bf [Knockoff]}, {\bf [FCD]}, {\bf [GtSt-BH]}, and:
\begin{itemize}
\item {\bf [Slope]}  Slope for FDR control, as presented in \citep{bogdan2015slope}.
\end{itemize}
The comparison is displayed in Figure \ref{f:candes}. It appears that GtSt-BH and FCD procedures are more conservative but they give a better control of the False Discovery Proportion (FDP).  SLOPE and Knockoff are more powerful but their FDP is greater than the expected FDR $ \alpha= 0.2$ (in $6$ experiments out of $7$ for SLOPE, in $3$ out of $7$ for Knockoff1, and in $5$ out of $7$ for Knockoff2). 

\section*{Acknowledgements}
The authors would like to thank anonymous referees for their time, comments and useful remarks on preliminary versions of this paper. The authors are in debt to {Quentin Duchemin} for his valuable comments and remarks on this work.

 
 \cleardoublepage

 \begin{center}
 \LARGE Supplement to\\ \textbf{Multiple Testing and Variable Selection along\\ the path of the Least Angle Regression}
 \end{center}

 \bigskip
 
 \section{Representing the~LARS knots}
\label{sec:LARknots}

\subsection{The equivalent formulations of the~LARS algorithm}
\label{sec:LAR3}
We present here three equivalent formulations of the~LARS that are a consequence of the analysis provided in Appendices~\ref{sec:LARknots} and \ref{app:FirstSteps}. One new formulation is given by Algorithm~\ref{alg:LAR3}.
{\tiny
\begin{algorithm}[!h]
\setcounter{AlgoLine}{0}
  \SetAlgoLined 
  
  \KwData{Correlations vector~$\bar Z$ and variance-covariance matrix~$\bar R$.} 
    
  \KwResult{Sequence~$((\lambda_k,\bar \imath_k, \varepsilon_k))_{k\geq1}$ where~$\lambda_1\geq\lambda_2\geq\ldots>0$ are the knots, and~$\bar \imath_1,\bar \imath_2,\ldots$ are the variables that enter the model with signs~$\varepsilon_1,\varepsilon_2,\ldots$~$(\varepsilon_k=\pm1)$.}
    \BlankLine
    
    {\nonl \tcc{Initialize computing~$(\lambda_1,\bar \imath_1,\varepsilon_1)$ and defining a ‘residual'~$\bar N^{(1)}$.}}
    Set~$k=1$,~$\displaystyle \lambda_1:=\max|\bar Z|$,~$\displaystyle \bar \imath_1:=\arg\max|\bar Z|$ and~$\varepsilon_1=\bar Z_{\bar \imath_1}/\lambda_1\in\pm1$, and~$\bar N^{(1)}:=\bar Z$.\\
          \BlankLine

      {\nonl \tcc{Note that~$((\lambda_\ell,\bar \imath_\ell, \varepsilon_\ell))_{1\leq \ell\leq k-1}$ and~$\bar N^{(k-1)}$ have been defined at the previous step.}}
      Set~$k\leftarrow k+1$ and compute the least-squares fit 
      \[
      \bar \theta_j:=\big(\bar R_{j,\bar \imath_1} \cdots \bar R_{j,\bar \imath_{k-1}}  \big) M^{-1}_{\bar \imath_1,\ldots, \bar \imath_{k-1}} (\varepsilon_1,\ldots,\varepsilon_{k-1})\,,\quad j=1,\ldots,p\,,
      \]
      where~$ M_{\bar \imath_1,\ldots, \bar \imath_{k-1}}$ is the sub-matrix of~$\bar R$ keeping the columns and the rows indexed by~$\{\bar \imath_1,\ldots, \bar \imath_{k-1}\}$.

      For~$0<\lambda\leq\lambda_{k-1}$ compute the ‘‘residuals''~$\bar N^{(k)}(\lambda)=(\bar N^{(k)}_1(\lambda),\ldots,\bar N^{(k)}_p(\lambda))$ given by
      \[
      \bar N_j^{(k)}(\lambda):=\bar N_j^{(k-1)}-(\lambda_{k-1}-\lambda)\bar\theta_j\,,\quad j=1,\ldots,p\,,
      \]
and pick
\begin{align*}
\lambda_{k}	&:=\max\big\{\beta>0\, ;\  \exists\, j\notin \{\bar \imath_1,\ldots, \bar \imath_{k-1}\},\ 
			\mathrm{s.t.}\ |\bar N_j^{(k)}(\beta)|=\beta\big\}\ \mathrm{and}\ 
\bar \imath_k		:=\argmax_{j\notin \{\bar \imath_1,\ldots, \bar \imath_{k-1}\}} |\bar N_j^{(k)}(\lambda_k)|\,,\\
\varepsilon_k	&:=\bar N_{\bar \imath_k}^{(k)}(\lambda_k)/\lambda_k\in\pm1\ \mathrm{and}\ 
\bar N^{(k)}	:=\bar N^{(k)}(\lambda_k)\,.
\end{align*}
 Then, iterate from {\bf 2}.
          \BlankLine
\caption{LARS algorithm (standard formulation)}
\label{alg:LAR}
\end{algorithm}
}

{\tiny
\begin{algorithm}[!h]
\setcounter{AlgoLine}{0}
  \SetAlgoLined 
  
  \KwData{Correlations vector~$\bar Z$ and variance-covariance matrix~$\bar R$.} 
  
  \KwResult{Sequence~$((\lambda_k,\bar \imath_k, \varepsilon_k))_{k\geq1}$ where~$\lambda_1\geq\lambda_2\geq\ldots>0$ are the knots, and~$\bar \imath_1,\bar \imath_2,\ldots$ are the variables that enter the model with signs~$\varepsilon_1,\varepsilon_2,\ldots$~$(\varepsilon_k=\pm1)$.}
    \BlankLine
    
      {\nonl \tcc{Initialize computing~$(\lambda_1,\bar \imath_1,\varepsilon_1)$.}}
    Define~$Z=(\bar Z,-\bar Z)$ and~$R$ as in~\eqref{eq:white}, and set~$k=1$,~$\displaystyle \lambda_1:=\max Z~$,~$\displaystyle \hat \imath_1:=\arg\max Z$,~$\bar \imath_1=\hat \imath_1 \mod p$ and~$\varepsilon_1=1-2(\,\hat \imath_1-\bar \imath_1)/p\in\pm1$.\\
          \BlankLine

      {\nonl \tcc{Note that~$((\lambda_\ell,\hat \imath_\ell))_{1\leq \ell\leq k-1}$ have been defined at the previous step/loop.}}
      Set~$k\leftarrow k+1$ and compute 
       \[
 \lambda_k=  \max_{\{j :\, \theta_j( \hat \imath_1, \ldots,\hat \imath_{k-1}) <1\}} 
 \bigg\{
 \frac{Z_j - \Pi_{\hat \imath_1,\ldots,\hat \imath_{k-1}}  (Z_j)}{ 1-\theta_j(\,\hat \imath_1,\ldots,\hat \imath_{k-1})}\bigg\}
 \quad \mathrm{and}\quad \hat \imath_k 		= \!\!\!\!\!\!\!\!\!\!\!\!\argmax_{\{j :\, \theta_j( \hat \imath_1, \ldots,\hat \imath_{k-1}) <1\}} 
 \bigg\{
 \frac{Z_j - \Pi_{\hat \imath_1,\ldots,\hat \imath_{k-1}}  (Z_j)}{ 1-\theta_j(\,\hat \imath_1,\ldots,\hat \imath_{k-1})}\bigg\}\,,
 \]
    where 
    \begin{align*}
    \Pi_{\hat\imath_1,\ldots,\hat \imath_{k-1}} (Z_j) & :=  \big(R_{j,\hat\imath_1} \cdots R_{j,\hat\imath_{k-1}}  \big)M^{-1}_{\hat\imath_1,\ldots,\hat\imath_{k-1}}  
(Z_{\hat\imath_1}, \ldots, Z_{\hat\imath_{k-1}})\\
\theta_j(\,\hat \imath_1,\ldots, \hat \imath_{k-1}) &:= \big(R_{j,\hat \imath_1} \cdots R_{j,\hat \imath_{k-1}}  \big)M^{-1}_{\hat \imath_1,\ldots, \hat \imath_{k-1}}  
(1, \ldots, 1)
    \end{align*}
and set~$\bar \imath_k=\hat \imath_k \mod p$ and~$\varepsilon_k	=1-2(\,\hat \imath_k-\bar \imath_k)/p\in\pm1$.
 Then, iterate from {\bf 2}.
          \BlankLine
\caption{LARS algorithm (‘‘projected'' formulation)}
\label{alg:LAR2}
\end{algorithm}
}

\newpage

\subsection{A new formulation of Least Angle Regression algorithm}
\label{sec:larrec}
The Least Angle Regression (LARS) algorithm has been introduced in the seminal article~\cite{efron2004least}. 
{In the context of linear regression in high dimensions, the LARS algorithm can be used to identify a subset of potential covariates. The LARS outputs a piecewise affine solutions path, and the {\it knots}~$\lambda_1 \geq \lambda_2\geq \dots >0$ are the change points of the LARS path that are built by tracking the~$\ell_{\infty}$ of the residual. At each knot, the LAR algorithm adds to the active set of variables the covariate the most correlated with the actual residual. In that way, the descent direction is always equiangular to all variables present in the current active set.} . This sequence of knots is closely related to the sequence of knots of LASSO~\citep{tibshirani1996regression}, as they differ by only one rule: ‘‘{\it Only in the LASSO case, if a nonzero coefficient crosses zero before the next variable enters, drop it from the active set and recompute the current joint least-squares direction}'', as mentioned in~\cite[Page 120]{hastie2015statistical} or~\cite[Theorem 1]{efron2004least} for instance. 

{\normalsize
\begin{algorithm}[!t]
\captionsetup{labelfont={sc,bf}, labelsep=newline}
\caption{LARS algorithm (‘‘recursive'' formulation)}
\label{alg:LAR3}
  \SetAlgoLined 
  {\nonl {\tcc{Given a response~$Y$ and a design~$X$, set~$\bar Z=X^\top Y$ and~$\bar R=X^\top X$}}}

  \KwData{Correlations vector~$\bar Z$ and variance-covariance matrix~$\bar R$.}
                 \BlankLine  
          
  \KwResult{Sequence~$((\lambda_k,\bar \imath_k, \varepsilon_k))_{k\geq1}$ where~$\lambda_1\geq\lambda_2\geq\ldots>0$ are the knots, and~$\bar \imath_1,\bar \imath_2,\ldots$ are the variables that enter the model with signs~$\varepsilon_1,\varepsilon_2,\ldots$~$(\varepsilon_k=\pm1)$.}
    \BlankLine
  
    \SetKwFunction{FMain}{Rec}
  \SetKwProg{Fn}{Function}{:}{}
    \SetKwProg{Re}{Return}{:}{}
  
     {\nonl \tcc{Define the recursive function \FMain{} that would be applied repeatedly. The inputs of  \FMain{} are~$Z$ a vector,~$R$ a SDP matrix and~$T$ a vector.}}
          \BlankLine  
    
    \DontPrintSemicolon
      {\nonl  \Fn{\FMain{$R$,~$Z$,~$T$}}{
             	{\nonl  Compute}
{\nonl  \begin{align*}
	\displaystyle\lambda&=\max_{\{j:\, T_j<1\}} \Bigg\{\frac{Z_j}{1-T_j}\Bigg\}\,, \displaystyle &{\bf i}&=\arg\max _{\{j:\, T_j<1\}} \Bigg\{\frac{Z_j}{1-T_j}\Bigg\}\,, &{\bf x}&= \frac{R_{{\bf i}}}{R_{{\bf i} {\bf i}}}\,.
	\end{align*}
		}
	       	{\nonl  Update}
		{\nonl  \begin{align*}
		R&\leftarrow R-{\bf x}R_{{\bf i}}^\top\,, &Z&\leftarrow Z-{\bf x} Z_{{\bf i}}\,,&T&\leftarrow T+{\bf x}(1-T_ {{\bf i}})\,.
		\end{align*}
		}
        {\nonl  \KwRet ($R$,~$Z$,~$T$,~$\lambda$,~${\bf i}$)\;}}}

           \BlankLine  
           
        {Set~$k=0$,~$T=0$,~$Z=(\bar Z,-\bar Z)$ and~$R=\left[\begin{array}{cc} \bar R & -\bar R \\-\bar R & \bar R\end{array}\right]$.}
          \BlankLine
          
      {Update~$k\leftarrow k+1$ and compute 
      \[
      (R, Z, T, \lambda_k, \hat\imath_k)=\text{\FMain{$R$,~$Z$,~$T$} }
      \]
      Set~$\bar \imath_k=\hat \imath_k \mod p$ and~$\varepsilon_k	=1-2(\,\hat\imath_k-\bar \imath_k)/p\in\pm1$.}

\end{algorithm}
\
      \vspace{-0.5cm}
}


{
\begin{theorem}
\label{thm:equivalent}
Let~$n,p$ be integers. Given a vector~$Y\in\mathbb R^n$ and matrix~$X\in\mathds R^{n\times p}$ of rank~$r$ then Algorithm~\ref{alg:LAR} (LARS standard formulation), Algorithm~\ref{alg:LAR2} (LARS projected formulation) and  Algorithm~\ref{alg:LAR3} (LARS recursive formulation) output the same sequence~$((\lambda_k,\bar \imath_k, \varepsilon_k))_{k=1}^r$ from the input given by~$\bar Z=X^\top Y$ and~$\bar R=X^\top X$, where~$\lambda_1\geq\lambda_2\geq\lambda_2\ldots\geq0$ are the knots, and~$\bar \imath_1,\bar \imath_2,\ldots,\bar \imath_r$ are the variables that enter the model with signs~$\varepsilon_1,\varepsilon_2,\ldots,\varepsilon_r$~$(\varepsilon_k=\pm1)$.
\end{theorem}
\noindent The formulation of Algorithm~\ref{alg:LAR} (LARS standard formulation), Algorithm~\ref{alg:LAR2} (LARS projected formulation) and the proof of Theorem~\ref{thm:equivalent} are given in Appendices~\ref{sec:LARknots} and \ref{app:FirstSteps}.
}

\subsection{Initialization: First Knot}
The first step of the~LARS algorithm (Step {\bf 1} in Algorithm~\ref{alg:LAR}) seeks the most correlated predictor with the observation. In our formulation, introduce the first residual~$N^{(1)}:=Z$ and observe that~$N^{(1)}:=(\bar N^{(1)},-\bar N^{(1)})$. We define the first knot~$\lambda_1>0$ as 
\[
\lambda_1=\displaystyle\max Z\quad \mathrm{and}\quad\hat \imath _1 = \argmax Z\,.
\]  
\noindent
One may see that this definition is consistent with~$\lambda_1$ in Algorithm~\ref{alg:LAR} and note that~$\hat \imath_1$ and~$(\bar \imath_1,\varepsilon_1)$ are related as in~\eqref{eq:signed_indices}.

The~LARS  algorithm is a forward algorithm that selects a new variable and maintains a residual at each step.  We also define
\eq
\label{eq:recu_res}
  N^{(2)}(\lambda)=N^{(1)} - (\lambda_1-\lambda) \theta(\,\hat \imath_1)\,,\quad 0<\lambda\leq\lambda_1\,,
\qe
and one can check that~$N^{(2)}(\lambda)=(\bar N^{(2)}(\lambda),-\bar N^{(2)}(\lambda))$ where~$\bar N(\lambda)$ is defined in Algorithm~\ref{alg:LAR}. It is clear that the coordinate~$\hat \imath_1$ of~$  N^{(2)}(\lambda)$  is equal to~$\lambda$. On the other hand~$  N^{(1)}=Z$ attains  its maximum  at the single point~$\hat \imath_1$. By continuity  this last property is kept   for~$\lambda$ in a left neighborhood of~$\lambda_1$. We search for the first value of~$\lambda$ such that this property is not met, {\it i.e.} the~largest value of~$\lambda$ such that 
$$
 \exists j \neq \hat \imath_1 \mbox{ such that }     N^{(2)}(\lambda)  = \lambda\,,
~$$
 as in Step~${\bf 3}$ of Algorithm~\ref{alg:LAR}. We call this value~$\lambda_2$ and one may check that this definition is consistent with~$\lambda_2$ in Algorithm~\ref{alg:LAR}. 
 
 Now, we can be more explicit about the expression of~$\lambda_2$. Indeed, we make the following discussion on the values of~$\theta_j(\,\hat \imath_1)~$.
 \begin{itemize}
 \item If~$\theta_j(\,\hat \imath_1)  \geq 1$ ,  since~$N^{(1)}_j  <N^{(1)}_{\hat \imath_1}$ for~$j\neq \hat \imath_1$ there is no hope to achieve the equality  between~$N^{(2)}_j(\lambda)~$ and~$N^{(2)}_{\hat \imath_1}(\lambda)=\lambda$ for~$0<\lambda\leq\lambda_1$ in view of~\eqref{eq:recu_res}.
 \item   Thus we limit our attention  to the~$j$'s such that~$\theta_j(\,\hat \imath_1) < 1$. We have equality~$N_j^{(2)}(\lambda)  = \lambda$  when
\[
 \lambda =   \frac{N^{(1)}_j -\lambda_1 \theta_j(\,\hat \imath_1)}{ 1-\theta_j(\,\hat \imath_1)}.
\]
 \end{itemize}
So we can also define the second  knot~$\lambda_2$ of the~LARS  as 
\[
 \lambda_2 
  =\max_{j: \theta_j(\,\hat \imath_1) 
  < 1} 
  \bigg\{\frac{Z_j -\Pi_{\hat \imath_1}(Z_j)}{ 1-\theta_j(\,\hat \imath_1)}\bigg\}\,.
\]
where~$\Pi_{i_1}(Z_j):=Z_{i_1}\theta_j(i_1)$. Remark that~$\Pi_{i_1}(Z_j)=\bbE(Z_j\ |\ Z_{i_1})$ is the regression of~$Z_j$ on~$Z_{i_1}$ when $\bbE Z=0$.

\subsection{Recursion: Next Knots}
\label{app:recursion}
The loop~$({\bf 2}\rightleftarrows{\bf 3})$ in Algorithm~\ref{alg:LAR} builds iteratively the knots~$\lambda_1,\lambda_2\ldots$ of the~LARS algorithm and some ‘‘residuals''~$\bar N^{(1)},\bar N^{(2)},\ldots$ defined in Step {\bf 3}. We will present here an equivalent formulation of these knots. 

Assume that~$k\geq2$ and we have build~$\lambda_1,\ldots,\lambda_{k-1}$ and selected the ‘‘signed'' variables~$\hat \imath_1, \ldots,\hat \imath_{k-1}$. Introduce~$ N^{(k-1)}:=(\bar N^{(k-1)},-\bar N^{(k-1)})$ and define
 \[
 N^{(k)}(\lambda)=N^{(k-1)}-(\lambda_{k-1}-\lambda) \theta(\,\hat \imath_1,\ldots, \hat \imath_{k-1})\,,\quad 0<\lambda\leq\lambda_{k-1}\,.
  \]
 Check that~$ \theta_j(\,\hat \imath_1,\ldots, \hat \imath_{k-1})=(\bar \theta_j,-\bar \theta_j)$ where we recall that we define 
 \[
      \bar \theta_j:=\big(\bar R_{j,\bar \imath_1} \cdots \bar R_{j,\bar \imath_{k-1}}  \big) M^{-1}_{\bar \imath_1,\ldots, \bar \imath_{k-1}} (\varepsilon_1,\ldots,\varepsilon_{k-1})\,,\quad j=1,\ldots,p\,,
      \]
      at Step~${\bf 2}$ and it holds that~$\hat \imath_\ell$ and~$(\bar \imath_\ell,\varepsilon_\ell)$ are related as in~\eqref{eq:signed_indices}. From this equality, we deduce that it holds~$N^{(k)}(\lambda)=(\bar N^{(k)}(\lambda),-\bar N^{(k)}(\lambda))$. One may also check that the coordinates~$\hat \imath_1,\ldots,\hat \imath_{k-1}$ of~$  N^{(k)}(\lambda)$  are equal to~$\lambda$.

Again  if we want to solve~$  N_j^{(k)}(\lambda) = \lambda 
$ for some~$j$, we have to limit our attention  to~$j$'s such that~$  \theta_j(\,\hat \imath_1,\ldots, \hat \imath_{k-1}) <1$. Solving this latter equality yields to
 \[
 \lambda_k=  \max_{j : \theta_j( \hat \imath_1, \ldots,\hat \imath_{k-1}) <1} 
 \bigg\{\frac{N^{(k-1)}_j -\lambda_{k-1} \theta_j(\,\hat \imath_1,\ldots,\hat \imath_{k-1})}{ 1-\theta_j(\,\hat \imath_1,\ldots,\hat \imath_{k-1})}
 \bigg\}
 \,.
 \]
 This expression is consistent with~$\lambda_k$ in Algorithm~\ref{alg:LAR}. 
 
 \newpage
 
 Now, we can give an other expression of~$\lambda_k$ that will be useful in the proofs of our main theorems. Note that the residuals satisfy the relation
 \eq
 \label{eq:value_nk}
N^{(k)}=N^{(k-1)}-(\lambda_{k-1}-\lambda_k) \theta(\,\hat \imath_1,\ldots, \hat \imath_{k-1})\,,
 \qe
 and that~$N^{(k-1)}_j = \lambda_{k-1}$ for~$j=\hat \imath_1,\ldots,\hat \imath_{k-1}$. The following lemma permits  a drastic simplification of the expression of the knots. Its proof is given in Appendix~\ref{proof:ydc1}.
 
    \begin{lemma}  \label{l:jm}
It holds 
  \[
 {N^{(k-1)} -\lambda_{k-1} \theta(\,\hat \imath_1,\ldots,\hat \imath_{k-1})}=Z - \Pi_{\hat \imath_1,\ldots,\hat \imath_{k-1}}  (Z)
  \]
where we denote~$\Pi_{ i_1,\ldots, i_{k-1}} (Z) = (\Pi_{ i_1,\ldots, i_{k-1}} (Z_1),\ldots,\Pi_{ i_1,\ldots, i_{k-1}} (Z_{2p}))$ and note that, for all~$j\in[2p]$, one has $
\Pi_{ i_1,\ldots, i_{k-1}} (Z_j) =  \big(R_{j,i_1} \cdots R_{j,i_{k-1}}  \big)M^{-1}_{i_1,\ldots, i_{k-1}}  
(Z_{i_1}, \ldots, Z_{i_{k-1}})$. 
\end{lemma}

\noindent 
Using Lemma~\ref{l:jm} we deduce that~$\lambda_k$ in Algorithm~\ref{alg:LAR} is consistent with
 \eq
 \label{e:lambda_max}
 \lambda_k=  \max_{j : \theta_j( \hat \imath_1, \ldots,\hat \imath_{k-1}) <1} 
\bigg\{ \frac{Z_j - \Pi_{\hat \imath_1,\ldots,\hat \imath_{k-1}}  (Z_j)}{ 1-\theta_j(\,\hat \imath_1,\ldots,\hat \imath_{k-1})}\bigg\}\,.
 \qe
 where~$
\Pi_{\hat\imath_1,\ldots,\hat \imath_{k-1}} (Z_j) =  \big(R_{j,\hat\imath_1} \cdots R_{j,\hat\imath_{k-1}}  \big)M^{-1}_{\hat\imath_1,\ldots,\hat\imath_{k-1}}  
(Z_{\hat\imath_1}, \ldots, Z_{\hat\imath_{k-1}})$. When~$\bbE Z=0$, one may remark that $\Pi_{ i_1,\ldots, i_{k-1}} (Z_j)$ is the regression of~$Z_j~$ on the vector~$ (Z_{i_1}, \cdots , Z_{i_{k-1}})$ whose variance-covariance matrix is~$M_{i_1,\ldots, i_{k-1}}$. This analysis leads to an equivalent formulation of the~LARS algorithm (Algorithm~\ref{alg:LAR}). We present this formulation in Algorithm~\ref{alg:LAR2}.

 \begin{remark} \label{ZAZA2}
  Note  that Algorithm~\ref{alg:LAR} implies that~$\hat \imath_1,\ldots,\hat \imath_k$ are pairwise different, but also  that they differ modulo~$p$. 
  \end{remark}


\section{First Steps to Derive the Joint Law of the~LARS knots}
\label{app:FirstSteps}

\subsection{Law of the First Knot}
One has the following lemma governing the law of~$\lambda_1$.
\begin{lemma}
\label{lem:FirstKnotLaw}
It holds that 
\begin{itemize}
\item~$Z_{i_1}$ is independent of~$(Z_j^{(i_1)})_{j\neq i_1}$,
\item If 
$\theta_j(i_1)<1$ for all~$j\neq i_1$ then $\big\{\hat \imath_1=i_1\big\}=\big\{\lambda_2^{(i_1)}\leq Z_{i_1}\big\}$,
\item If 
$\theta_j(i_1)<1$ for all~$j\neq i_1$ then, conditional on~$\{\hat \imath_1=i_1\}$ and~$\lambda_2$,~$\lambda_1$ is a truncated Gaussian random variable with mean~$\E(Z_{i_1})$ and variance~$\rho_1^2:=R_{\hat \imath_1,\hat \imath_1}$ subject to be greater than~$\lambda_2$.
\end{itemize}
\end{lemma}
\begin{proof}
The first point is  a consequence or the properties of Gaussian regression. Now, observe that 
\begin{align*}
\big\{\lambda_2^{(i_1)}\leq Z_{i_1}\big\}&\Leftrightarrow\big\{\forall j\neq i_1\,,\ \frac{Z_j-Z_{i_1}\theta_j(i_1)}{1-\theta_j(i_1)}\leq Z_{i_1}\big\}\\
&\Leftrightarrow\big\{\forall j\neq i_1\,,\ {Z_j-Z_{i_1}\theta_j(i_1)}\leq Z_{i_1}-Z_{i_1}\theta_j(i_1)\big\}\\
&\Leftrightarrow\big\{\forall j\neq i_1\,,\ {Z_j}\leq Z_{i_1}\big\}\\
&\Leftrightarrow\big\{\hat \imath_1=i_1\big\}\,,
\end{align*}
as claimed. The last statement is a consequence of the two previous points.
\end{proof}

\subsection{Recursive Formulation of the~LARS}
One has the following proposition whose proof can be found in Section~\ref{proof:moulinette}. As we will see in this section, this intermediate result as a deep consequence, the~LARS algorithm can be stated in a recursive way applying the same function repeatedly, as presented in Algorithm~\ref{alg:LAR3}. 

\begin{proposition}
\label{prop:moulinette}
Set
\[
\tau_{j,i_k}
:=\frac
{R_{j,i_k} - \big(R_{j, i_1} \cdots R_{j,i_{k-1}}  \big)M^{-1}_{i_1,\ldots, i_{k-1}} \big(R_{i_k, i_1} ,\cdots ,R_{ i_k, i_{k-1}} \big)}
{(1-\theta_j(i_1,\ldots, i_{k-1}))(1-\theta_{i_k}(i_1,\ldots, i_{k-1}))}\,,
\]
and observe that~$\tau_{j,i_k}$ is the covariance between~$Z^{(i_1,\ldots, i_{k-1})}_j$ and~$Z^{(i_1,\ldots, i_{k-1})}_{i_k}$. Furthermore, it holds 
\eq
\label{eq:sign_tau}
\frac{\tau_{j,i_k}}{\tau_{i_k,i_k}}=1-\frac{1-\theta_j(i_1,\ldots, i_{k})}{1-\theta_j(i_1,\ldots, i_{k-1})}
\qe 
and 
\eq
\label{eq:rec_Zj}
\forall j\neq i_1,\ldots,i_{k}\,,\quad
Z^{(i_1,\ldots, i_{k})}_j=
\frac{Z^{(i_1,\ldots, i_{k-1})}_j-Z^{(i_1,\ldots, i_{k-1})}_{i_k}\tau_{j,i_k}/\tau_{i_k,i_k}}
{1-\tau_{j,i_k}/\tau_{i_k,i_k}}\,.
\qe
\end{proposition}

\noindent
Now, we present Algorithm~\ref{alg:LAR3}. Define~$R(0):=R$,~$Z(0)=Z$ and~$T(0)=0$. For~$k\geq1$ and~$i_1, \ldots, i_{k}\in[2p]$, introduce
\begin{align*}
R(k)&:=\Big({R_{j,\ell} - \big(R_{j,i_1} \cdots R_{j,i_{k}}  \big)M^{-1}_{i_1,\ldots, i_{k}} \big(R_{\ell, i_1} ,\cdots ,R_{ \ell, i_{k}} \big)}\Big)_{j,\ell}\\
Z(k)&:=Z-\Pi_{ i_1,\ldots, i_{k}} (Z)\\
T(k)&:=(\theta_j(i_1,\ldots, i_{k}))_j\,,
\end{align*}
and note that~$R(k)$ is the variance-covariance matrix of the Gaussian vector~$Z(k)$. The key property is following. Let~$v_1,  \ldots, v_k,$ be~$k$ linearly independent vectors of an Euclidean space and let~$u$ be any vector of the space. Set
\[
v:=  P_{(v_1, \ldots, v_{k-1})} ^\perp v_k,
\] 
the  projection of~$v_k$ orthogonally to~$v_1, \ldots, v_{k-1}$.Then 
\[ 
 P_{(v_1, \ldots, v_k)}^\perp u = P_{v}^\perp P_{(v_1, \ldots, v_{k-1})} u\,.
\]
Using this result we deduce that
\begin{align}
Z(k)
&=\Pi^\perp_{ i_1,\ldots, i_{k}} (Z)\notag\\
&=\Pi^\perp_{ i_{k}}(\Pi^\perp_{ i_1,\ldots, i_{k-1}} (Z))\notag\\
&=\Pi^\perp_{ i_{k}}(Z(k-1))\notag\\
&=Z(k-1)-\Pi_{ i_{k}}(Z(k-1))\notag\\
&=Z(k-1)-{\bf x}(k) Z(k-1)\label{eq:LAR3_1}\,,
\end{align}
where~${\bf x}(k)=R_{i_k}(k-1)/R_{i_k,i_k}(k-1)$. It yields that 
\eq
\label{eq:LAR3_2}
R(k)=R(k-1)-{\bf x}(k) R_{i_k}(k-1)^\top\,.
\qe
Using~\eqref{eq:sign_tau} (or~\eqref{eq:rec101}), remark that 
\begin{align}
T(k)
&=T(k-1)-{\bf x}(k) (1-T_{i_k}(k-1))
\label{eq:LAR3_3}
\,.
\end{align}
These relations give a recursive formulation of the~LARS as presented in Algorithm~\ref{alg:LAR3}.

\section{Proofs}

\subsection{Proof of Proposition~\ref{prop:reccur}}
\label{proof:reccur}
\noindent
{\bf First and third points: } The first point works by induction. The initialization of the proof is given by the second point of Lemma~\ref{lem:FirstKnotLaw}. We will use Proposition~\ref{prop:moulinette} to prove the first point. We have
\begin{align}
&\lambda_{k+1}^{(i_1,\ldots, i_{k})}\leq Z_{i_k}^{(i_1,\ldots, i_{k-1})}\notag\\
&\Leftrightarrow
\forall j\neq i_1,\ldots, i_{k}\,,\ Z^{(i_1,\ldots, i_{k})}_j\leq Z_{i_k}^{(i_1,\ldots, i_{k-1})}\notag\\
&\Leftrightarrow
\forall j\neq i_1,\ldots, i_{k}\,,\ Z^{(i_1,\ldots, i_{k-1})}_j-Z^{(i_1,\ldots, i_{k-1})}_{i_k}\frac{\tau_{j,i_k}}{\tau_{i_k,i_k}}\leq Z_{i_k}^{(i_1,\ldots, i_{k-1})}-Z_{i_k}^{(i_1,\ldots, i_{k-1})}\frac{\tau_{j,i_k}}{\tau_{i_k,i_k}}\notag\\
&\Leftrightarrow
\forall j\neq i_1,\ldots, i_{k}\,,\ Z^{(i_1,\ldots, i_{k-1})}_j\leq Z_{i_k}^{(i_1,\ldots, i_{k-1})}\label{eq:rec_sign1}\\
&\Leftrightarrow
\lambda^{(i_1,\ldots, i_{k-1})}_{k}= Z_{i_k}^{(i_1,\ldots, i_{k-1})}\notag\,.
\end{align}
using~\eqref{eq:rec_Zj} and that~$1-{\tau_{j,i_k}}/{\tau_{i_k,i_k}}>0$ (which is a consequence of~\eqref{eq:sign_tau} and~\eqref{hyp:IrrAlongThePath})  in~\eqref{eq:rec_sign1}. By induction and using~\eqref{eq:rec_sign1}, it holds that
\begin{align*}
\big\{
&\lambda_{k+1}^{(i_1,\ldots, i_{k})}\leq Z_{i_k}^{(i_1,\ldots, i_{k-1})}\leq Z_{i_{k-1}}^{(i_1,\ldots, i_{k-2})}\leq\ldots\leq Z_{i_2}^{(i_1)}\leq Z_{i_1}\big\}\notag\\
&\Leftrightarrow\big\{\forall j\neq i_1,\ldots, i_{k}\,,\ 
Z^{(i_1,\ldots, i_{k-1})}_j\leq Z_{i_k}^{(i_1,\ldots, i_{k-1})}\leq Z_{i_{k-1}}^{(i_1,\ldots, i_{k-2})}\leq\ldots\leq Z_{i_2}^{(i_1)}\leq Z_{i_1}\big\}\\
&\Leftrightarrow\big\{\forall j\neq i_1,\ldots, i_{k-1}\,,\ 
Z^{(i_1,\ldots, i_{k-1})}_j\leq Z_{i_{k-1}}^{(i_1,\ldots, i_{k-2})}\leq\ldots\leq Z_{i_2}^{(i_1)}\leq Z_{i_1}\\
&\quad\quad \text{and }\forall j\neq i_1,\ldots, i_{k}\,,\  Z^{(i_1,\ldots, i_{k-1})}_j\leq Z_{i_k}^{(i_1,\ldots, i_{k-1})}
\big\}\\
&\Leftrightarrow\big\{\lambda^{(i_1,\ldots, i_{k-1})}_{k}\leq Z_{i_{k-1}}^{(i_1,\ldots, i_{k-2})}\leq\ldots\leq Z_{i_2}^{(i_1)}\leq Z_{i_1} \text{ and }\lambda^{(i_1,\ldots, i_{k-1})}_{k}=Z_{i_k}^{(i_1,\ldots, i_{k-1})}
\big\}\\
&\ \  \vdots\\
&\Leftrightarrow\big\{\lambda^{(i_1,\ldots, i_{a-1})}_{a}\leq Z_{i_{a-1}}^{(i_1,\ldots, i_{a-2})}\leq\ldots\leq Z_{i_2}^{(i_1)}\leq Z_{i_1}\\
&\quad\quad  \text{ and }\lambda^{(i_1,\ldots, i_{k-1})}_{k}= Z_{i_k}^{(i_1,\ldots, i_{k-1})}
\leq\ldots\leq \lambda^{(i_1,\ldots, i_{a-1})}_{a}=Z_{i_a}^{(i_1,\ldots, i_{a-1})}
\big\}\tag{$s_a$}\label{eq:stopa}\\
&\ \  \vdots\\
&\Leftrightarrow\big\{\hat \imath_1=i_1,\ldots,\hat \imath_{k}=i_{k}\big\}\,.
\end{align*}
Now, observe that~$\hat i_{k+1}$ is the (unique) arg max of~$\lambda_{k+1}^{(i_1,\ldots, i_{k})}$ on the event~$\big\{\hat \imath_1=i_1,\ldots,\hat \imath_{k}=i_{k}\big\}$.

\noindent It yields that 
\begin{align*}
\big\{\hat \imath_1=i_1,\ldots,\hat \imath_{k+1}=i_{k+1}\big\}=\big\{
\lambda_{k+1}^{(i_1,\ldots, i_{k})}=Z_{i_{k+1}}^{(i_1,\cdots, i_{k})}
\leq 
Z_{i_k}^{(i_1,\cdots, i_{k-1})}\leq\cdots\leq 
Z_{i_2}^{(i_1)}\leq
Z_{i_1}\big\}\,,
\end{align*}
as claimed. Stopping at~$a$ as in~\eqref{eq:stopa} gives the second part of the statement. The third point of the proposition is a direct consequence of the first point.

\medskip

\noindent
{\bf Second point: }The proof of the second point can be lead by induction. The initialization of the proof is given by the first point of Lemma~\ref{lem:FirstKnotLaw}. Now, observe that~$Z_{i_{k}}^{(i_1,\ldots, i_{k-1})}, \cdots, Z_{i_{2}}^{(i_1)}, Z_{i_1}$ are measurable functions of~$(Z_{i_1},\ldots,Z_{i_{k}})$ and one may check that~$(Z_{i_1},\ldots,Z_{i_{k}})$ is independent of $(Z^{(i_1,\ldots, i_{k})}_j)_{ j\neq i_1, \ldots, i_{k}}$. We deduce that~$Z_{i_{k}}^{(i_1,\ldots, i_{k-1})}, \cdots, Z_{i_{2}}^{(i_1)}, Z_{i_1}$ are independent of~$(Z^{(i_1,\ldots, i_{k})}_j)_{ j\neq i_1, \ldots, i_{k}}$. One can also check that~$Z_{i_{k+1}}^{(i_1,\ldots, i_{k})}$ is independent of~$(Z^{(i_1,\ldots, i_{k})}_j)_{ j\neq i_1, \ldots, i_{k+1}}$.  Deduce that the variables $Z_{i_{K}}^{(i_1,\ldots, i_{K-1})}, \cdots, Z_{i_{2}}^{(i_1)}, Z_{i_1}$ are independent and independent of $(Z^{(i_1,\ldots, i_{K})}_j)_{ j\neq i_1, \ldots, i_{K}}$. This gives the first part of the second point of the proposition. 

Now observe that $(Z^{(i_1,\ldots, i_{K})}_j)_{ j\neq i_1, \ldots, i_{K}}$ is a linear function of $P^\perp_{K}Y$. More precisely, consider the vector $\tilde W:=(Z^{(i_1,\ldots i_K)}_{j})_{j\in[p]}$ and note that 
\[
(\hat \sigma^{i_1,\ldots, i_{K}})^2=\frac{\|P_K^\perp(Y)\|_2^2}{n-K}
\quad\mathrm{and}\quad
\tilde W=\mathrm{Diag}(\mathds{1}_p-\theta^K)^{-1}\times X^\top\big(P_K^\perp(Y)\big)\,,
\]
where $\theta^K:=(\theta_{j}(i_1,\ldots, i_K))_{j\in[p]}$, $P_K^\perp=\mathrm{Id}_n-P^{(i_1,\ldots, i_K)}$, and with the convention $0/0=0$. The key remark is that the direction and the norm of a centered Gaussian vector are independent, namely $\|P_K^\perp(Y)\|$ is independent of $P_K^\perp(Y)/\|P_K^\perp(Y)\|$ when $\mathds EP_K^\perp(Y)=0$. Observe that $\mathds EP_K^\perp(Y)=P_K^\perp(X\beta^0)$ and $\mathds EP_K^\perp(Y)=0$ is equivalent to $X\beta^0\in H_K$. We deduce that if $X\beta^0\in H_K$ then $\hat \sigma^{i_1,\ldots, i_{K}}$ is independent of $\tilde W/\sigma^{i_1,\ldots, i_{K}}$, as claimed.

\subsection{Proof of Proposition \ref{cor:Main1} } 
\label{proof:cor4}
Fix~$a$ such that~$0\leq a\leq K-1$ and consider any selection procedure~$\hat m$ satisfying~\eqref{e:stopping_rule}. From Theorem~\ref{thm:Main1} (more precisely~\eqref{zaz0}) we know that the density of~$(\lambda_{a+1},\lambda_{a+2},\ldots,\lambda_K)$ conditional on 
\[
 \mathcal F:=\big\{\hat \imath_1 =i_1, \ldots, \hat \imath_{K}=i_{K}, \lambda_a,\lambda_{K+1}\big\}\,
 \]
 is given by 
    \eq
    \notag
    (const)   \Big(\prod_{k=a+1}^K \varphi_{m_k,v_k^2}(\ell_k)\Big)\, \1_{\lambda_a \leq \ell_{a+1}\leq \cdots \leq \ell_K\leq  \lambda_{K+1}}\,.
   \qe
   From Proposition~\ref{prop:H0}, conditional on~$\mathcal F$ and under the null hypothesis described in Proposition \ref{cor:Main1}, we know that~$m_{a+1}=\ldots=m_K=0$. It implies that~$\Phi_k$ is the CDF of~$\lambda_k$ for~$a<k\leq K$. 
 
From the definition of a stopping time given by~\eqref{e:stopping_rule} and on the event~$\mathcal F$, we know that~$\1_{\{\hat m=a\}}$ is a measurable function of~$\lambda_1,\ldots, \lambda_{a-1}$ which are respectively equal to~$Z_{i_1},\ldots,  Z_{i_{a-1}}^{(i_1,\ldots, i_{a-2})}$ on~$\mathcal F$ by~\eqref{e:frozen} (as proven in Appendix~\ref{app:recursion}  and Eq.~\eqref{e:lambda_max}). By Proposition~\ref{prop:reccur} (more precisely~\eqref{e:indep_sigmas}), we also know that this function is independent of~$(\lambda_{a+1},\lambda_{a+2},\ldots,\lambda_K)$ conditional on~$\mathcal F$. Remark that its is also independent of~$\hat \sigma^{i_1,\ldots, i_{K}}$ conditional of~$\mathcal F$ for the same reason (it would be useful later, when we will build testing procedures when the variance is unknown). We deduce that the conditional density above is also the conditional density on the event 
 \[
 \mathcal G:=\big\{\hat m =a,\hat \imath_1 =i_1, \ldots, \hat \imath_{K}=i_{K},\lambda_a,\lambda_{K+1}\big\}\,.
 \]
From~$F_k=\Phi_k(\lambda_k)$ ({\it i.e.,} applying the CDF) we deduce by a change of variables that conditional on the selection event~$ \mathcal G$, the vector~$(F_{a+1},\ldots,F_K)$ is uniformly distributed on 
\begin{align*}
\cD_{a+1,K}:=&
\big\{
(f_{a+1},\ldots,f_K)\in\bbR^{K-a}:\ \\
&\cP_{a+1,a}(F_{a})\geq f_{a+1}\geq \cP_{a+1,a+2}(f_{a+2})
\geq \cdots
\geq \cP_{a+1,K}(f_K)\geq \cP_{a+1,K+1}(F_{K+1})
\big\}\,,
\end{align*}
where~$\cP_{i,j}$ are described in~\eqref{eq:pij}.

\subsection{Proof of Proposition~\ref{prop:fabc}}
\label{proof:fabc}

By Proposition~\ref{cor:Main1}, a simple integration shows that  
\eq
\notag
\bbP\big[\lambda_b\leq t\ |\ \hat m =a,\lambda_a,\lambda_c, \hat \imath_{1},\ldots,\hat \imath_{a},\hat \imath_{a+1}, \ldots,\hat \imath_{c-1},\hat \imath_{c},\ldots,\hat \imath_{K}\big]=\frac{\bbF_{abc}(t)}{\bbF_{abc}(\lambda_a)} \ ,
\qe
under the null hypothesis of Proposition \ref{prop:fabc} (which implies that~$m_{a+1}=\ldots=m_K=0$). Then note that the function~$\bbF_{abc}$ is defined by~$\sigma,\lambda_a,\lambda_c,\hat \imath_{a+1}, \ldots,\hat \imath_{c-1}$ only. We deduce that we can de-condition on~$\hat m =a, \hat \imath_{1},\ldots,\hat \imath_{a}, \hat \imath_{c},\ldots,\hat \imath_{K}$, which gives the result.

\subsection{Proof of Proposition~\ref{p:frozen}}
  \label{proof:frozen}
 Let us fix some values~$i_1, \ldots, i_{K+1}$. Recall  that the frozen values of the knots
\eq
\notag
 \lambda^f_1 : = Z_{i_1}, \ldots , \lambda^f_{K} :=Z^{(i_1,\ldots i_{K-1})}_{i_{K}},\lambda^f_{K+1} :=   Z^{(i_1,\ldots i_K)}_{i_{K+1}},
  \qe
  are Gaussian, independent, and 
 $
\hat \sigma^{i_1,\ldots, i_{K}}
\indep \big(Z_{i_{K+1}}^{(i_1,\ldots, i_{K})}/\hat \sigma^{i_1,\ldots, i_{K}}\big)
\indep Z_{i_{K}}^{(i_1,\ldots, i_{K-1})}
\indep \cdots 
\indep Z_{i_{2}}^{(i_1)}
\indep Z_{i_1}
$, see Proposition~\ref{prop:reccur} and~\eqref{e:indep_sigmas}. Let us condition by~$\{\hat\imath_1=i_1,\ldots, \hat \imath_K=i_K, \lambda_{K+1} = \ell_{K+1}\}$. Note that, on this event,~$ \lambda_h^f = \lambda_h$ ,~$h \in [K]$ and by Proposition~\ref{prop:reccur} this event is equivalent  to 
 \eq
 \label{e:se_prop}
\{ Z_{i_1}>\cdots >Z^{(i_1,\ldots i_{K-1})}_{i_{K}}> \max_j Z^{(i_1,\ldots i_K)}_{j}=\lambda_{K+1} = \ell_{K+1} \}.
\qe 
Because of the independence above, we get the conditional independence stated in the first and second point of the proposition. 

\noindent
For the last point, consider $\tilde W:=(Z^{(i_1,\ldots i_K)}_{j})_{j\in[p]}$ and note that 
\[
(\hat \sigma^{i_1,\ldots, i_{K}})^2=\frac{\|P_K^\perp(Y)\|_2^2}{n-K}
\quad\mathrm{and}\quad
\tilde W=\mathrm{Diag}(\mathds{1}_p-\theta^K)^{-1}\times X^\top\big(P_K^\perp(Y)\big)\,,
\]
where $\theta^K:=(\theta_{j}(i_1,\ldots, i_K))_{j\in[p]}$, $P_K^\perp=\mathrm{Id}_n-P^{(i_1,\ldots, i_K)}$, and with the convention $0/0=0$. Now, let $U\in\mathds R^{n\times (n-K)}$ be any matrix such that $UU^\top=P_K^\perp$ and define $W:=U^\top Y/\sigma$, then 
\[
\Big(\frac{\hat \sigma^{i_1,\ldots, i_{K}}}{\sigma}\Big)^2=\frac{\|W\|_2^2}{n-K}
\quad\mathrm{and}\quad
\tilde W=\sigma\times\mathrm{Diag}(\mathds{1}_p-\theta^K)^{-1}\times X^\top UW\,.
\]
Because of the independence above and \eqref{e:se_prop}, the distribution of~$\hat \sigma$ is independent of the other variables and such that 
\[
\Big(\frac{\hat \sigma}{\sigma}\Big)^2=\frac{\|W\|_2^2}{n-K}\quad\mathrm{with\ }W\mathrm{\ s.t.\ }\|\mathrm{Diag}(\mathds{1}_p-\theta^K)^{-1}\times X^\top UW\|_\infty=\ell_{K+1}/\sigma\,.
\]
This implies that the conditional distribution is the one claimed.


\newpage

 \subsection{Orthogonal Case: Proof of Theorem~\ref{thm:zazPower}}
 \label{proof:zazPower}
Let~$\mathcal I$ the set of admissible indexes 
~$$
 \mathcal I := \{a,b,c: a_0\leq a<b<c \leq K+1 \}\,.
~$$

\begin{figure}
\includegraphics[width=0.3\textwidth]{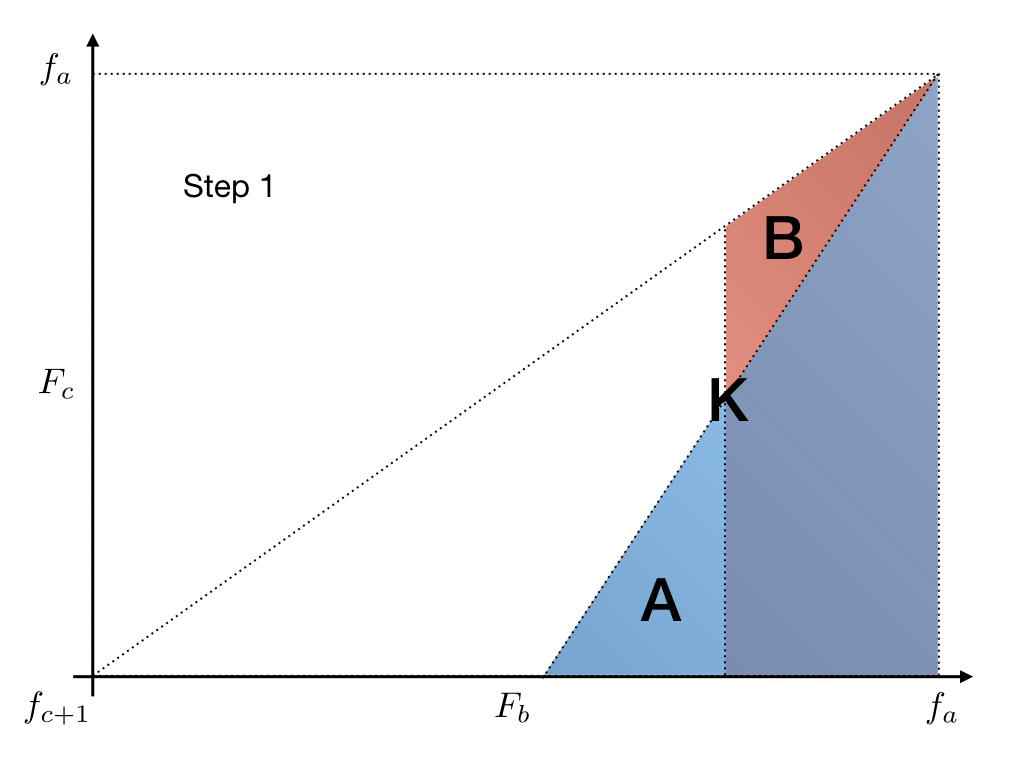}
\includegraphics[width=0.3\textwidth]{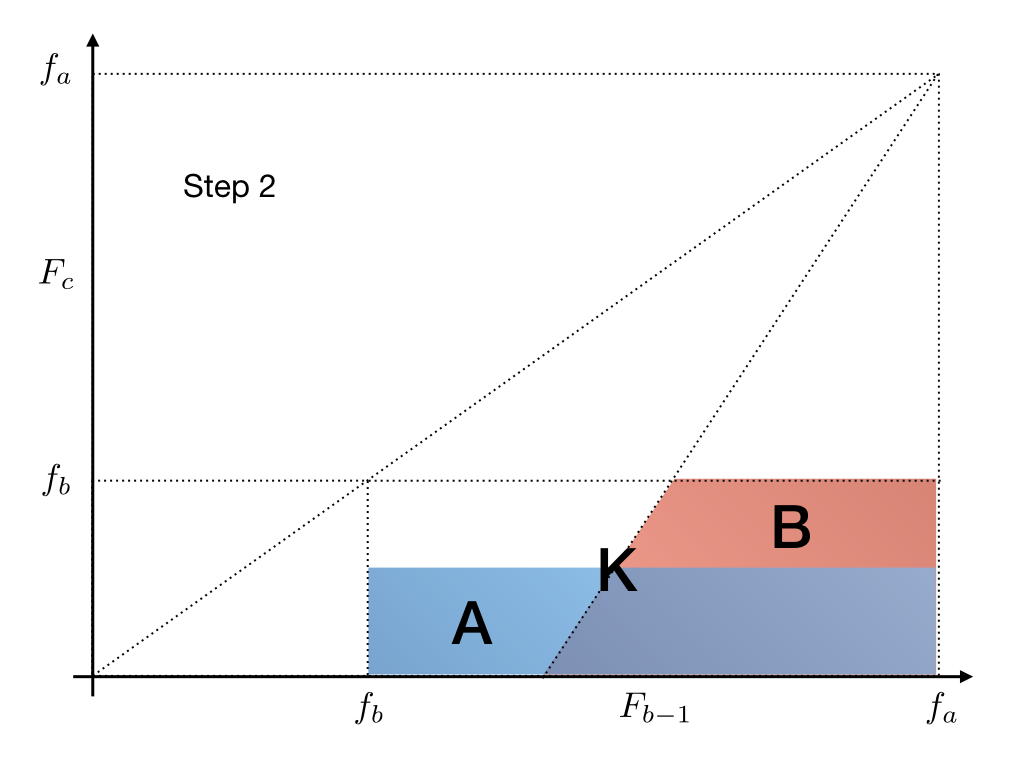}
\includegraphics[width=0.3\textwidth]{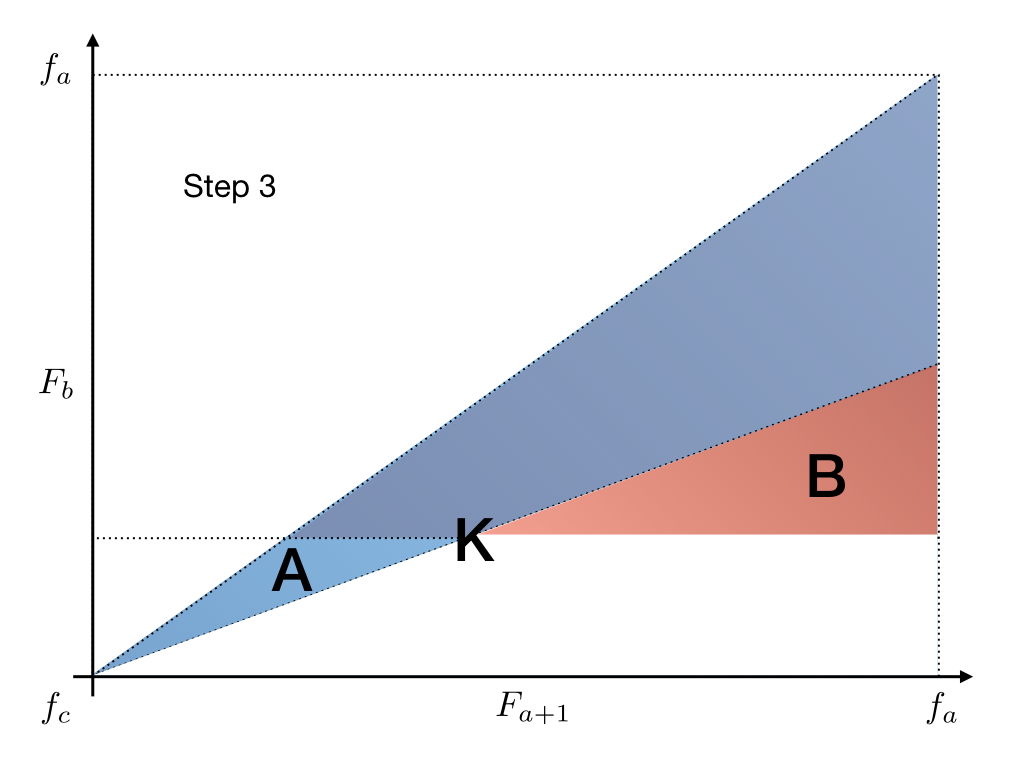}
\caption{Rejection domains associated to the different comparison sets appearing in steps of the proof of Theorem~\ref{thm:zazPower}.}
\label{fig:power}
\end{figure}

 \noindent
 {$\circ$ \bf Step 1: } We prove that, when the considered  indexes such that~$c+1\leq K+1$ belong to~$\mathcal I$,~$\cS_{a,b,c+1}$ is more powerful than~$\cS_{a,b,c}$. Our proof is conditional to~$F_a=f_a,F_{c+1} =f_{c+1}~$.  Note that~$(F_{a+1}, \ldots, F_c)~$ has for distribution the uniform distribution  on the simplex
~$$
  \mathcal S := \{ f_a>F_{a+1} > \cdots > F_c >f_{c+1} \}.
~$$ 
  This implies  by direct calculations that 
~$$
  \cI_{ab}(s,t) =\frac{ (s-t) ^{b-a-1}}{(b-a-1)!}
~$$
  and that 
\eq
\label{eq:BetaF}
   \frac{   \bbF_{abc}(\lambda_b)   }{\bbF_{abc}(\lambda_a)} = \bF_{\beta((b-a),(c-b))} \bigg( \frac{ F_b-F_c}{ f_a-F_c} \bigg).
\qe
  where~$\bF_\beta$ is the cumulative distribution  of the Beta distribution in reference.  Using monotony of this  function the~$\cS_{abc}$ test has for rejection region 
  \begin{equation} \label{ZAZA:e:r1}
  (F_b-F_c)  \geq z_1 (f_a-F_c) \Leftrightarrow  F_b \geq z_1 f_a +(1-z_1) F_c,
 \end{equation}
  where~$z_1$ is some threshold, depending on~$\alpha$, that belongs to~$(0,1)$.

  Similarly~$\cS_{ab(c+1)}$ has for rejection region 
  \begin{equation} \label{ZAZA:e:r2}
  F_b \geq z_2 (f_a-f_{c+1})  + f_{c+1},
  \end{equation}
  where~$z_2$ is some other threshold belonging to~$(0,1)$. We use the following lemmas.  
 
 \begin{lemma} \label{z:l:1}
  Let~$c \leq K$.  The density~$h_\mu$ of~$f_1, \ldots, f_c$, conditional on~$F_{c+1}$ with respect of the Lebesgue measure  under the alternative is  coordinate-wise  non-decreasing and given by~\eqref{zazcosh}.
  \end{lemma}
  \begin{proof}
Observe that it suffises to prove the result when~$\sigma=1$. Note that 
~$$
 \lambda_{c+1}^{i_1, \dots, i_c } =   \max_{j\in[p]\,,\   j\neq \bar \imath_1,\ldots,   j \neq  \bar \imath_c}  |Z_j|.
~$$
  Thus its density~$ p_{\mu^0,i_1,\ldots,i_c}$ does not depend  on~$\mu^0_{ \bar i_1},\ldots,\mu^0_{ \bar i_c}$. As a consequence the following variables have the same distribution ;
$ \lambda_{c+1}^{i_1 +\epsilon_1 p, \dots, i_c +\epsilon_c p}~$ , where~$\epsilon_1,\ldots,\epsilon_c$ take the value 0 or 1 and  indices are taken modulo~$p$.

Because of the independence of the different variables, the joint density, under the alternative hypothesis, of $\lambda_1, \ldots,  \lambda_{c+1}$ taken at~$ \ell_1, \ldots, \ell_{c+1},$  on the domain~$ \{\lambda_1>\cdots > \lambda_{c+1}\}$ takes the value
~$$
 (Const)   \sum'   \big(\varphi(\ell_1 -\mu^0_{j_1} ) +\varphi(\ell_1 +\mu^0_{j_1} ) \big),\ldots ,
   \big(\varphi(\ell_c -\mu^0_{j_c} ) +\varphi(\ell_c +\mu^0_{j_c} ) \big)
  p_{\mu^0,j_1,\ldots,j_K}(\ell_{k+1}).
~$$
Here the sum~$ \displaystyle \sum'$ is taken over all different~$j_1, \ldots, j_c$ belonging  to~$\llbracket 1, p\rrbracket$.
 
 Then the density, conditional on~$F_{c+1} = f_{c+1}$,   of
~$ F_1 ,\ldots,F_{c}~$  at~$f_1,\ldots,f_c$  takes the value
~
\eq
\label{zazcosh}
 \mathrm{(const)}  \sum'  \cosh( \mu_{j_1} f_1)\ldots \cosh( \mu_{j_c} f_c) \1_{f_1>\cdots> f_{c}> F_{c+1}},
\qe
 implying  that this density is coordinate-wise  non-decreasing. 
 \end{proof}
 
 \begin{lemma} \label{jma:lemma}  Let~$\nu_0$  the image  on the plane~$(F_b,F_c)$  on the uniform probability  on~$\mathcal S$: it is the distribution  under the null of~$(F_b,F_c)$. The two rejection regions :~$\mathcal R_1$ associated to 
\eqref{ZAZA:e:r1} and~$\mathcal R_2$ associated to 
\eqref{ZAZA:e:r2} 
have  of course the same probability~$\alpha$ under~$\nu_0$. Let~$\eta_{\mu^0}$ the density w.r.t.~$\nu_0$ of the distribution of~$(F_b,F_c)$ under the alternative. Then~$\eta_{\mu^0}$     is non decreasing coordinate-wise. 
  \end{lemma}
  \begin{proof}  
Integration yields that density of~$\nu_0$ w.r.t. the Lebesgue    measure taken at point~$(f_b,f_c)$  is 
\[
 \frac{  (f_a-f_b)^{b-a-1}   (f_b-f_c)^{c-b-1}  }{(b-a-1)!(c-b-1)!}.
\]
  The density of~$\nu_{\mu^0}$  w.r.t. Lebesgue measure  is 
  \begin{equation}\label{ZAZA:e:vm}
  \int _{f_b} ^{f_a} df_{a+1} \ldots   \int _{f_b} ^{f_b-2} df_{a+1}  \int _{f_b} ^{f_b-2} df_{b-1} 
   \int _{f_c} ^{f_b} df_{b+1} \ldots   \int _{f_c} ^{f_{c-2}} df_{c-1}   h_{\mu^0}(f_{a},\ldots ,f_c).
 \end{equation}
  Thus~$\eta_{\mu^0}$ which is the quotient  of these two quantities is just  a mean value  of~$h_{\mu^0}$ on the domain  of integration~$\mathcal D_{f_b,f_c}$  in~\eqref{ZAZA:e:vm}.
  
  Suppose that~$f_b$ and~$ f_c$  increase, then all the borns of the domain~$\mathcal D_{f_b,f_c}~$ increase also. By Lemma~\ref{z:l:1} the mean value  increases.
\end{proof}

\medskip
   
    \noindent
     {\bf We finish now the proof of  Step 1:}  For a given level~$\alpha$ let us consider the two rejection regions~$ R_{a,b,c}$ and~$R_{a,b,(c+1)}~$  of the two considered tests in the plane~$F_b,F_c$ and set
     \[
     A:=R_{a,b,c} \setminus R_{a,b,(c+1)} \text{ and } B:= R_{a,b,(c+1)}\setminus R_{a,b,c}\,,
     \] 
     see Figure \ref{fig:power}.  These two regions have  the same~$\nu_0$ measure. By elementary geometry there exist a point~$K=(K_b,K_c)$  in the plane such that 
      \begin{itemize}
      \item For every point of~$A$,~$F_b \leq K_b$,~$F_c \leq  K_c$,
      \item For every point of~$B$,~$F_b \geq K_b$,~$F_c \geq  K_c$,
\end{itemize}
      By transport of measure there exists a transport function~$\cT$  that  preserve the measure~$\nu_0$  and  that is one-to one~$A \to B$.   As a consequence  the transport  by~$\cT$ improve the probability under the alternative: the power of~$\cS_{a,b,c+1}$ is~larger than that of 
     ~$\cS_{a,b,c}$.  \bigskip

 \noindent
 {$\circ$ \bf Step 2: } We prove that, when the considered  indexes belong to~$\mathcal I$ such that~$a<b-1$,~$\cS_{a,(b-1),c}$ is more powerful than~$\cS_{a,b,c}$. 
 Our proof is conditional on~$F_a=f_a,F_b =f_b~$ and is located  
 in the plane~$(F_{b-1}, F_c)$. 
 
 
 \noindent
  The rejection region~$R_{a,b,c}$ takes the form 
~$F_c \leq  \frac{1}{1-z_1}f_b  -\frac{z_1}{1-z_1}f_a~$
   for some threshold~$z_1$ belonging to~$(0,1)$.\newline
The rejection region~$R_{a,(b-1),c}$ takes the form~$F_c \leq \frac{1}{1-z_2} F_{b-1} - \frac{z_2}{1-z_2} f_a~$ for some  other threshold~$z_2$ belonging to~$(0,1)$. \newline
These regions as well as the regions~$A$ and~$B$ and the point~$K$  are indicated in Figure \ref{fig:power}.

Transport of measure  and the convenient modification of Lemma \ref{jma:lemma} imply that the power of the test $\cS_{a,(b-1),c}$ is greater of equal than that of~$\cS_{a,b,c}$.  \bigskip

 \noindent
 {$\circ$ \bf Step 3: } We prove that, when the considered  indexes belong to~$\mathcal I$ such that~$a+1<b$,~$\cS_{a,b,c}$ is more powerful than~$\cS_{(a+1),b,c}$. 
  Our proof is conditional on~$F_a=f_a,F_{c} =f_{c}~$ and is located in the plane~$F_{a+1}, F_b$.\newline
The rejection region~$R_{a,b,c}$ takes the form~$F_b \geq z_1 f_a + (1-z_1)f_c~$ for some threshold~$z_1$belonging to~$(0,1)$. \newline
 The rejection region~$R_{a+1,b,c}$ takes the form~$F_b \geq z_2 F_{a+1} + (1-z_2)f_c ~$ for some  other threshold~$z_2$ belonging to~$(0,1)$.\newline
These regions as well as the regions~$A$ and~$B$ and the point~$K$  are indicated in Figure \ref{fig:power}.

Transport of measure  and  the convenient modification of Lemma \ref{jma:lemma} imply that the power of~$\cS_{a,b,c}$ is greater of equal  that that of~$\cS_{(a+1),b,c}$.  \bigskip

\noindent
Considering the three cases above,  we get the desired result.

\subsection{Proof of Lemma~\ref{l:jm}}
\label{proof:ydc1}
The proof works by induction. Let us check the relation for~$k=2$, namely
\[
 N^{(1)} -\lambda_{1} \theta(\,\hat \imath_1)=Z-Z_{\hat \imath_1} \theta(\,\hat \imath_1)=Z - \Pi_{\hat \imath_1}  (Z)\,.
\]
Now, let~$k\geq 3$. First, the three perpendicular theorem  implies  that  for every~$j,  i_1,\ldots, i_{k-1}~$,
\begin{align*}
	\theta_j(i_1, \ldots ,  i_{k-2})
   	&= \big(R_{j,i_1} \cdots R_{j,i_{k-1}}  \big)M^{-1}_{i_1,\ldots, i_{k-1}}   
		( 
 			\theta_{i_1}(i_1, \ldots ,  i_{k-2}),\ldots,\theta_{i_{k-1}} (i_1,\ldots, i_{k-2})  
  		)\,,\\
	\mathrm{and}\ \Pi_{ i_1,\ldots, i_{k-2}}  (Z_j)
	&= \big(R_{j,i_1} \cdots R_{j,i_{k-1}}  \big)M^{-1}_{i_1,\ldots, i_{k-1}}   
		( 
  			\Pi_{ i_1,\ldots, i_{k-2}}  (Z_{i_1}),\ldots,\Pi_{ i_1,\ldots, i_{k-2}}  (Z_{i_{k-1}})  
  		)\,.
\end{align*}
By induction, using~\eqref{eq:value_nk}, we get that
\begin{align}
N^{(k-1)}&=N^{(k-2)}-(\lambda_{k-2}-\lambda_{k-1}) \theta(\,\hat \imath_1,\ldots, \hat \imath_{k-2})\,,\notag
\\
&=
(N^{(k-2)}-\lambda_{k-2} \theta(\,\hat \imath_1,\ldots, \hat \imath_{k-2}))+\lambda_{k-1}\theta(\,\hat \imath_1,\ldots, \hat \imath_{k-2})\,,\notag\\
&=
Z  - \Pi_{\hat \imath_1,\ldots,\hat \imath_{k-2}}  (Z)+\lambda_{k-1} \theta(\,\hat \imath_1,\ldots, \hat \imath_{k-2})\,.\label{eq:induc}
\end{align}
Then, recall that~$N^{(k-1)}_j = \lambda_{k-1}$ for~$j=\hat \imath_1,\ldots,\hat \imath_{k-1}$ and remark that
\[
\lambda_{k-1} \theta_j(\,\hat \imath_1,\ldots,\hat \imath_{k-1}) 
= \big(R_{j,\hat \imath_1} \cdots R_{j,\hat \imath_{k-1}}  \big)M^{-1}_{\hat \imath_1,\ldots, \hat \imath_{k-1}} 
 (N^{(k-1)}_{\hat \imath_1},\ldots,N^{(k-1)}_{\hat \imath_{k-1}})\,.
\]
Using~\eqref{eq:induc} at indices~$j=\hat \imath_1,\ldots,\hat \imath_{k-1}$, we deduce that
\begin{align*}
\lambda_{k-1} \theta_j(\,\hat \imath_1,&\ldots,\hat \imath_{k-1}) \\
&
= \big(R_{j,\hat \imath_1} \cdots R_{j,\hat \imath_{k-1}}  \big)M^{-1}_{\hat \imath_1,\ldots, \hat \imath_{k-1}} 
 (N^{(k-1)}_{\hat \imath_1},\ldots,N^{(k-1)}_{\hat \imath_{k-1}})
\\
& = \big(R_{j,\hat \imath_1} \cdots R_{j,\hat \imath_{k-1}}  \big)M^{-1}_{\hat \imath_1,\ldots, \hat \imath_{k-1}} 
   (Z_{\hat \imath_1},\ldots,Z_{\hat \imath_{k-1}})
\\
 &\quad - \big(R_{j,\hat \imath_1} \cdots R_{j,\hat \imath_{k-1}}  \big)M^{-1}_{\hat \imath_1,\ldots, \hat \imath_{k-1}}   
	( 
  		\Pi_{\hat\imath_1,\ldots,\hat\imath_{k-2}}  (Z_{\hat \imath_1}),\ldots,\Pi_{\hat\imath_1,\ldots, \hat \imath_{k-2}}  (Z_{\hat \imath_{k-1}})  
  	)
\\
&\quad + \lambda_{k-1}\big(R_{j,\hat \imath_1} \cdots R_{j,\hat \imath_{k-1}}  \big)M^{-1}_{\hat \imath_1,\ldots, \hat \imath_{k-1}}   
( 
 \theta_{\hat \imath_1}(\,\hat \imath_1, \ldots , \hat\imath_{k-2}),\ldots,\theta_{\hat \imath_{k-1}} (\,\hat \imath_1,\ldots,\hat\imath_{k-2})  
  ) 
\\
& = \Pi_{\hat\imath_1,\ldots,\hat\imath_{k-1}}  (Z_{j})- \Pi_{\hat\imath_1,\ldots,\hat\imath_{k-2}}  (Z_j)
+\lambda_{k-1} \theta_j(\,\hat \imath_1,\ldots,\hat \imath_{k-2}) \,,
  \end{align*}
  Namely
  \[
  \Pi_{\hat \imath_1,\ldots,\hat \imath_{k-1}}  (Z)-\lambda_{k-1} \theta(\,\hat \imath_1,\ldots, \hat \imath_{k-1})
  =
  \Pi_{\hat \imath_1,\ldots,\hat \imath_{k-2}}  (Z)-\lambda_{k-1} \theta(\,\hat \imath_1,\ldots, \hat \imath_{k-2})\,.
  \]
Using again~\eqref{eq:induc} we get that
\begin{align*}
 N^{(k-1)}&= Z  - \Pi_{\hat \imath_1,\ldots,\hat \imath_{k-2}}  (Z)+\lambda_{k-1} \theta(\,\hat \imath_1,\ldots, \hat \imath_{k-2})\,,\\
 &=Z  - \Pi_{\hat \imath_1,\ldots,\hat \imath_{k-1}}  (Z)+\lambda_{k-1} \theta(\,\hat \imath_1,\ldots, \hat \imath_{k-1})\,,
\end{align*}
as claimed.

\subsection{Proof of Proposition~\ref{prop:moulinette}}
\label{proof:moulinette}
We denote
\begin{align*}
R_{j}&:=\big(R_{j, i_1}, \ldots, R_{j,i_{k-1}}  \big)\,,\\
R_{i_k}&:=\big(R_{i_k, i_1} ,\ldots ,R_{ i_k, i_{k-1}} \big)\,,\\
M&:=M_{i_1,\ldots, i_{k-1}}\,,\\
\bar M&:=M_{i_1,\ldots, i_{k}}=\left[\begin{array}{cc}M & R_{i_k} \\R_{i_k}^\top & R_{i_k,i_k}\end{array}\right]\,,\\
\bar R&:=\big(R_{j, i_1}, \ldots, R_{j,i_{k}}  \big)\,,\\
x&:=\frac{1-\theta_j(i_1,\ldots, i_{k-1})}{1-\theta_{i_k}(i_1,\ldots, i_{k-1})}\frac{\tau_{j,i_k}}{\tau_{i_k,i_k}}\,,
\end{align*}
and observe that
\begin{align}
x&=\frac
{R_{j,i_k} -R_j^\top M^{-1}R_{i_k}}
{R_{i_k,i_k} -R_{i_k}^\top M^{-1}R_{i_k}}\,,\notag\\
{\bar M}^{-1}&=
\left[\begin{array}{cc}\mathrm{Id}_{k-1} & -M^{-1}R_{i_k} \\0 & 1\end{array}\right]
\left[\begin{array}{cc}M^{-1} & 0 \\0 & \big(R_{i_k,i_k} -R_{i_k}^\top M^{-1}R_{i_k}\big)^{-1}\end{array}\right]
\left[\begin{array}{cc}\mathrm{Id}_{k-1} & 0 \\-R_{i_k}^{\top}M^{-1} & 1\end{array}\right]\,,\notag\\
{\bar M}^{-1}\bar R&=
\left[\begin{array}{c}M^{-1}\big(R_j-xR_{i_k}\big) \\ x\end{array}\right]\,,
\label{eq:Schur1}
\end{align}
using Schur complement of block~$M$ of the matrix~$\bar M$ and a LU decomposition. Note also that
\[
\frac{Z^{(i_1,\ldots, i_{k-1})}_j-Z^{(i_1,\ldots, i_{k-1})}_{i_k}\tau_{j,i_k}/\tau_{i_k,i_k}}
{1-\tau_{j,i_k}/\tau_{i_k,i_k}}
=\frac
{Z_j-\Pi_{ i_1,\ldots, i_{k-1}} (Z_j)-x\,(Z_{i_k}-\Pi_{ i_1,\ldots, i_{k-1}} (Z_{i_k}))}
{1-\theta_j(i_1,\ldots, i_{k-1})-x(1-\theta_{i_k}(i_1,\ldots, i_{k-1}))}\,.
\]
To prove~\eqref{eq:rec_Zj}, it suffices to show that the R.H.S term  above is equal to the following R.H.S term
\[
Z^{(i_1,\ldots, i_{k})}_j=
\frac{Z_j-\Pi_{ i_1,\ldots, i_{k}} (Z_j)}
{1-\theta_j(i_1,\ldots, i_{k})}\,.
\]
We will prove that numerators are equal and that denominators are equal. For denominators, 
\begin{align}
&1-\theta_j(i_1,\ldots, i_{k-1})-x(1-\theta_{i_k}(i_1,\ldots, i_{k-1}))\notag\\
&=1-\theta_j(i_1,\ldots, i_{k-1})-x+x\,\theta_{i_k}(i_1,\ldots, i_{k-1})\notag\\
&=1-(\underbrace{1\cdots 1}_{k\ \mathrm{times}})
\left[\begin{array}{c}M^{-1}\big(R_j-xR_{i_k}\big) \\ x\end{array}\right]\notag\\
&=1-\theta_j(i_1,\ldots, i_{k})\label{eq:rec101}\,,
\end{align}
using~\eqref{eq:Schur1}. Furthermore, it proves~\eqref{eq:sign_tau}. For the numerators, we use that
\begin{align*}
&Z_j-\Pi_{ i_1,\ldots, i_{k-1}} (Z_j)-x\,(Z_{i_k}-\Pi_{ i_1,\ldots, i_{k-1}} (Z_{i_k}))\\
&=Z_j-\Pi_{ i_1,\ldots, i_{k-1}} (Z_j)-x Z_{i_k}+ x \Pi_{ i_1,\ldots, i_{k-1}} (Z_{i_k})\\
&=Z_j-(Z_{i_1}\cdots Z_{i_k})
\left[\begin{array}{c}M^{-1}\big(R_j-xR_{i_k}\big) \\ x\end{array}\right]\\
&=Z_j-\Pi_{ i_1,\ldots, i_{k}} (Z_j)\,.
\end{align*}
using~\eqref{eq:Schur1}.

\subsection{Proof of Theorem \ref{thm:independentTest}}
\label{proof:yoFDR}
We rely on the {\bf Weak Positive Regression Dependency (WPRDS) property} to prove the result, one may consult~\cite[Page 173]{giraud2014introduction} for instance. We say that a function~$g:[0,1]^K\to \bbR^+$ is {\it nondecreasing} if for any~$p,q\in[0,1]^K$ such that~$p_k\geq q_k$ for every~$k=1,\ldots, K$, we have~$g(p)\geq g(q)$. We say that a Borel set~$\Gamma\in[0,1]^K$ is {\it nondecreasing} if~$g=\mathds 1_\Gamma$ is nondecreasing. In other words  if~$y \in \gamma$ and if~$z\geq0$, then~$y+z \in \gamma$. We say that the~$p$-values~$(\,\hat p_1=\hat\alpha_{0,1,K+1},\ldots,\hat p_K=\hat\alpha_{K-1,K,K+1})$ satisfy the WPRDS property if for any nondecreasing set~$\Gamma$ 
 and for all~$k^0\in I_0$, the function 
\[
u\mapsto \bbP_{\mu^0}\big[(\,\hat p_1,\ldots,\hat p_K)\in\Gamma\big|\hat p_{k^0}\leq u\big]\text{ is nondecreasing}
\]
where~$\mu^0=\beta^0$ in our  orthogonal design case, and we recall that
\[
I_0=\big\{k\in[K]\ :\ \mathds H_{0,k}\text{ is true}\big\}\,.
\] 
To prove Theorem \ref{thm:independentTest}, note that it is sufficient~\cite[Chapter~8]{giraud2014introduction} to prove that 
\eq
\label{yoyonondecresaing}
u\mapsto \overline\bbP\big[(\,\hat p_1,\ldots,\hat p_K)\in\Gamma\big|\hat p_{k^0}\leq u\big]\text{ is nondecreasing}
\qe
where~$\bar\bbE,\overline{\P}$ will denote that expectations and probabilities are conditional on~$\{\bar \imath_1,\ldots, \bar \imath_K,\lambda_{K+1}\}$ and under the hypothesis that~$\mu^0=X^\top X\beta^0$. Note that one can integrate in~$\lambda_{K+1}$ to get the statement of Theorem \ref{thm:independentTest}.

\medskip

\noindent
 {$\circ$ \bf Step 1: }We start by giving the joint law of the LARS knots under the alternative in the orthogonal design case. Lemma~\ref{z:l:1} and~\eqref{zazcosh} show that, conditional on~$\{\bar \imath_1,\ldots, \bar \imath_K,\lambda_{K+1}\}$,~$(\lambda_1,\ldots,\lambda_K)$ is distributed on the set~$\lambda_1\geq \lambda_2\geq \cdots\geq \lambda_K\geq \lambda_{K+1}$ and it has a coordinate-wise nondecreasing density. 
 Now we can assume without loss of generality  that~$\sigma^2 =1$, in addition because of orthogonality 
~$\rho_k^2=1$ implying that~$F_k=\Phi(\lambda_k)$~$\cP_{i,j}=\Phi_i\circ\Phi_j^{-1}=\mathrm{Id}$.  We deduce that, conditional on~$\{\bar \imath_1,\ldots, \bar \imath_K,F_{K+1}\}$,~$(F_1,\ldots,F_K)$ is distributed on the set
\[
\big\{
(f_1,\ldots,f_K)\in\bbR^K\ :\ 
1\geq f_1\geq f_2\geq \cdots\geq f_K\geq F_{K+1})
\big\}\,,
\]
it has an {\bf explicit density} given by~\eqref{zazcosh}, and we denote it by~$h_{\mu^0}$. By the change of variables~$G_k:=\frac{ F_k-F_{K+1}}{ F_{k-1}-F_{K+1}}$ one  obtains that the distribution of~$(G_1,\ldots,G_K)$ is supported on~$[0,1]^K$. More precisely, define
\begin{align*}
\psi(f_1,\ldots,f_K)& :=(g_1,\ldots,g_K)   := (\frac{ f_1-F_{K+1}}{1-F_{K+1}},\ldots,\frac{ f_K-F_{K+1}}{f_{K-1}-F_{K+1}}) \\
\psi^{-1}(g_1,\ldots,g_K)&:=\Big((1-F_{K+1})g_1+F_{K+1},\ldots,(1-F_{K+1})g_1g_2\ldots g_K+F_{K+1}\Big),
\end{align*}
whose inverse Jacobian determinant is 
\[
\det\Big[\frac{\partial \psi}{\partial f_1} \cdots\frac{\partial \psi}{\partial f_K}\Big]^{-1}=\prod_{k=1}^K(f_{k-1}-F_{K+1})=(1-F_{K+1})^K\prod_{k=1}^{K}g_k^{K-k}\,.
\] 
We deduce that the density  of~$(G_1,\ldots,G_K)|\{\bar \imath_1,\ldots,\bar \imath_K,F_{K+1}\}$  at point~$g$ with respect to Lebesgue measure is 
\begin{align}
\label{eq:density}
\bp(g):=
\mathrm{(const)}\mathds 1_{g\in(0,1)^K}\prod_{k=1}^{K}g_k^{K-k}&\cosh\big[{\mu^0_{\bar\imath_k}}((1-F_{K+1})\prod_{\ell=1}^kg_\ell+F_{K+1})\big]\,,
\end{align}
where we have used~\eqref{zazcosh}. From~\eqref{zaza5} and~\eqref{eq:BetaF}, one has 
\eq
\label{eq:yoyoFDR}
\hat p_k =1-\bF_{\beta(1,K-k+1)} \bigg( \frac{ F_k-F_{K+1}}{ F_{k-1}-F_{K+1}} \bigg)
=1-\bF_{\beta(1,K-k+1)}(G_k)
\qe
 where~$\bF_\beta$ is the cumulative distribution  of the Beta distribution in reference. We deduce that for any~$v\in(0,1)$ and for any~$\ell\in[K]$, 
\begin{align*}
\hat p_k &= v
\Leftrightarrow
(G_1,\ldots,G_K)\in 
{[0,1]^K\cap\{\bF_{\beta(1,K-k+1)}^{-1}(1-v)= G_k\}
}
\,,
\end{align*}
so that 
\eq
\label{yoyo:Fg}
\bar\bbP\big[(\,\hat p_1,\ldots,\hat p_K)\in\Gamma\big|\hat p_{k^0}\leq u\big]=\bar\bbP\big[(G_1,\ldots,G_K)\in\bar\Gamma\big|G_{k^0}\geq \bF_{\beta(1,K-\ell+1)}^{-1}(1-u)\big]\,,
\qe
where~$\bar\Gamma$ can be proved to be a {\bf nonincreasing Borel set} from~\eqref{eq:yoyoFDR}.

\medskip

\noindent
 {$\circ$ \bf Step 2: }
 Let~$0<x<y<1$ and denote by~$\mu_x$ the following conditional law  
 \[
 \mu_x:=\mathrm{law}\big[(G_1,\ldots,G_K)|\{\bar \imath_1,\ldots,\bar \imath_K,F_{K+1},G_{k^0}\geq x\}\big]\,.
 \] 
 Remark that if there exists a measurable~$T:[0,1]^K\mapsto [0,1]^K$ such that 
 \begin{itemize}
 \item~$T$ is nondecreasing, meaning that for any~$g\in[0,1]^K$,~$T(g)\geq g$;
 \item~$T$ is such that push-forward of~$\mu_x$ by~$T$ gives~$\mu_y$, namely~$T_\#\mu_x=\mu_y$;
 \end{itemize}
then it holds
 \begin{itemize}
 \item~$\mathds 1_{\{T(g)\in\bar\Gamma\}}\leq\mathds 1_{\{g\in\bar\Gamma\}}$;
 \item~$\mathrm{law}\big[T(G)|\{\bar \imath_1,\ldots,\bar \imath_K,F_{K+1},G_{k^0}\geq x\}\big]=\mathrm{law}\big[G|\{\bar \imath_1,\ldots,\bar \imath_K,F_{K+1},G_{k^0}\geq y\}\big]$ where~$G=(G_1,\ldots,G_K)$.
 \end{itemize}
In this case, we deduce that 
\[
\bar\bbP\big[G\in\bar\Gamma\big|G_{k^0}\geq x\big]
\geq
\bar\bbP\big[T(G)\in\bar\Gamma\big|G_{k^0}\geq x\big]
=
\bar\bbP\big[G\in\bar\Gamma\big|G_{k^0}\geq y\big]\,.
\]
 If one can prove that such function~$T$ exists for any~$0<x<y<1$, it proves that 
 \[
 x\mapsto \bar\bbP\big[G\in\bar\Gamma\big|G_{k^0}\geq x\big]\text{ is nonincreasing}\,,
 \]
and, in view of~\eqref{yoyo:Fg}, it proves~\eqref{yoyonondecresaing}. Proving that such function~$T$ exists is done in the next step.

\medskip

\noindent
 {$\circ$ \bf Step 3: }
 Let~$0<x<y<1$. Consider the {\bf Knothe-Rosenblatt transport map}~$T$ of~$\mu_x$ toward~$\mu_y$ following the order
 \[
 k^0\to k^0+1\to\cdots\to K\to k^0-1\to k^0-2\to\cdots\to 1\,.
 \]
 It is based on a sequence of conditional quantile transforms defined following the ordering above. Its construction is presented for instance in~\cite[Sec.2.3, P.67]{santambrogio2015optimal} or~\cite[P.20]{villani2008optimal}. The transport~$T$ is defined as follows. Given~$z,z'\in[0,1]^K$ such that~$z'=T(z)$ it holds
 \begin{align*}
 z'_{k^0}&=T^{(k^0)}(z_{k^0});\\
 z'_{k^0+1}&=T^{(k^0+1)}(z_{k^0+1},z'_{k^0});\\
 &\vdots\\
 z'_K&=T^{(K)}(z_K,z'_{K-1},\ldots,z'_{k^0});\\
 z'_{k^0-1}&=T^{(k^0-1)}(z_{k^0-1},z'_K,\ldots,z'_{k^0});\\
 &\vdots\\
 z'_1&=T^{(1)}(z_1,z'_2,\ldots,z'_{k^0-1},z'_K,\ldots,z'_{k^0});
 \end{align*}
where~$T^{(k^0)}, T^{(k^0+1)},\ldots,T^{(K)},T^{(k^0-1)}\ldots,T^{(1)}$ will be build in the sequel, in which we will drop their dependencies in the~$z'_k$'s to ease notations. It remains to prove that 
 \begin{itemize}
 \item~$T$ is nondecreasing, meaning that for any~$g\in[0,1]^K$,~$T(g)\geq g$;
 \item~$T$ is such that push-forward of~$\mu_x$ by~$T$ gives~$\mu_y$, namely~$T_\#\mu_x=\mu_y$;
 \end{itemize}
 to conclude. The last point is a property of the Knothe-Rosenblatt transport map. Proving the first point will be done in the rest of  the proof.
 
 \medskip
 
\noindent
 {$\circ$ \it Step 3.1: } We start by the first transport map~$T^{(k^0)}:[0,1]\mapsto[0,1]$. Denote~$\mu_x^{(k^0)}$ the following conditional law  
 \[
 \mu_x^{(k^0)}:=\mathrm{law}\big[G_{k^0}|\{\bar \imath_1,\ldots,\bar \imath_K,F_{K+1},G_{k^0}\geq x\}\big]\,,
 \]
 and~$\bbF_x^{(k^0)}$ its cdf. Note that the Knothe-Rosenblatt construction gives~$T^{(k^0)}=(\bbF_y^{(k^0)})^{-1}\circ\bbF_x^{(k^0)}$. We would like to prove that~$T^{(k^0)}(t)\geq t$ for all~$z\in(0,1)$. This is equivalent to prove that it holds~$ \bbF_x^{(k^0)}\geq \bbF_y^{(k^0)}$. For~$t\leq y$,~$\bbF_y^{(k^0)}(t)=0$ and it implies that~$ \bbF_x^{(k^0)}(t)\geq \bbF_y^{(k^0)}(t)$. Let~$t>y$, using the conditional density~$\bp$ defined in~\eqref{eq:density}, note that
 \begin{align*}
 \bbF_x^{(k^0)}(t)\geq \bbF_y^{(k^0)}(t)\quad&\Leftrightarrow\quad 
 \frac  {\int_x^t\bp}  {\int_x^1\bp}    \geq 
 \frac  {\int_y^t\bp}  {\int_y^1\bp} \\
 &\Leftrightarrow\quad \int_x^t \int_y^1 \bp\otimes\bp
 \geq \int_y^t   \int_x^1\bp\otimes\bp,
 \end{align*}
where, for example
\[
\int_x^t \mbox{ means the integral over the  hyper rectangle }  [x,t]
:=\Big\{(g_1,\ldots,g_K)\in[0,1]^K\ :\!\ x\leq g_{k^0}\leq t\Big\}.
\]
A simple calculation (see also Figure~\ref{fig:yoyo}) gives that
\[
\int_x^t \int_y^1 \bp\otimes\bp = \int_y^t   \int_x^1\bp\otimes\bp +
\int_{[x,y] \times [t,1]}
\bp\otimes\bp\,,
\]
and it proves that~$ \bbF_x^{(k^0)}\geq \bbF_y^{(k^0)}$.

\begin{figure}
\includegraphics[width=0.4\textwidth]{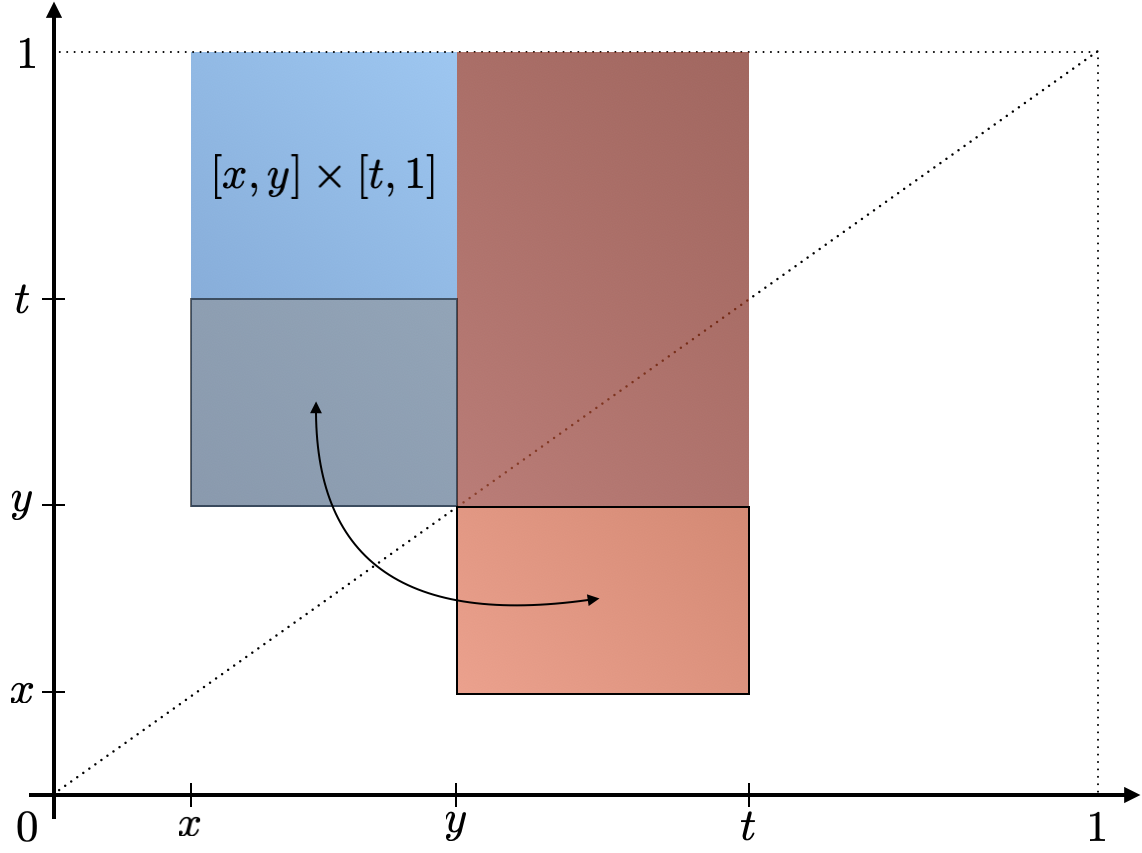}
\caption{Note that, by symmetry  the two boxed regions  have same~$ \bp\otimes\bp$ measure. The blue region is~$[x,t]\times [y,1]~$, its measure is the measure of the red region (namely~$cD_y(t)\times \cD_x(1)$) more the bluest upper left corner (namely~$ [x,y]\times[t,1]$). }
\label{fig:yoyo}
\end{figure}
 
 \medskip
 
 \noindent
 {$\circ$ \it Step 3.2: } We continue with the second transport map in Knothe-Rosenblatt construction. Let~$z_{k^0}\in(x,1)$ and denote~$\mu_{z_{k^0}}^{(k^0+1)}$ the following conditional law  
 \[
 \mu_{z_{k^0}}^{(k^0+1)}:=\mathrm{law}\big[G_{k^0+1}|\{\bar \imath_1,\ldots,\bar \imath_K,F_{K+1}, G_{k^0}=z_{k^0}\}\big]\,,
 \]
 and~$\bbF_{z_{k^0}}^{(k^0+1)}$ its cdf. Let~$z'_{k^0}:=T^{(k^0)}(z'_{k^0})$ and denote~$\mu_{z'_{k^0}}^{(k^0+1)}$ the following conditional law  
 \[
 \mu_{z'_{k^0}}^{(k^0+1)}:=\mathrm{law}\big[G_{k^0+1}|\{\bar \imath_1,\ldots,\bar \imath_K,F_{K+1}, G_{k^0}=z'_{k^0}\}\big]\,,
 \]
 and~$\bbF_{z'_{k^0}}^{(k^0+1)}$ its cdf. Note that~$x<z_{k^0}\leq z'_{k^0}=T^{(k^0)}(z_{k^0})\leq 1$. Again, we would like to prove that~$\bbF_{z_{k^0}}^{(k^0+1)}\geq \bbF_{z'_{k^0}}^{(k^0+1)}$ which implies that the transport map~$T^{(k^0+1)}:=\big(\bbF_{z'_{k^0}}^{(k^0+1)}\big)^{-1}\circ \bbF_{z_{k^0}}^{(k^0+1)}$ satisfies~$T^{(k^0+1)}(u)\geq u$ for all~$u\in(0,1)$.
 
 Recall that the conditional density~$\bp$ of~$G|\{\bar \imath_1,\ldots,\bar \imath_K,F_{K+1}\}$ is given by~\eqref{eq:density} and recall that~$k^0\in I_0$. Observe that~$\mu^0_{k^0}=0$, so that the conditional density of~$G|\{\bar \imath_1,\ldots,\bar \imath_K,F_{K+1}, G_{k^0}=\bz\}$ is
 \begin{align}
 \mathrm{(const)}\,\mathds 1_{g\in(0,1)^K}\mathds 1_{g_{k^0}=\bz}&\prod_{k< k^0}g_k^{K-k}\cosh\big[{\mu^0_{\bar\imath_k}}((1-F_{K+1})\prod_{\ell=1}^kg_\ell+F_{K+1})\big]\label{yoyoz1}
 \\&\times \prod_{k> k^0}g_k^{K-k}\cosh\big[{\mu^0_{\bar\imath_k}}((1-F_{K+1})\,\,\bz\!\!\prod_{1\leq \ell\neq k^0\leq k}\!\!g_\ell+F_{K+1})\big]\,.\notag
 \end{align}
 Set~$\tau:=z'_{k^0}/z_{k^0}\geq 1$ and~$G'_{k^0+1}=\tau G_{k^0+1}$ so that 
 \[
 z_{k^0} G_{k^0+1}=z'_{k^0}G'_{k^0+1}\,.
 \]
 Denote~$G':=(G_1,\ldots,G_{k^0-1},G'_{k^0+1},G_{k^0+2},\ldots,G_K)\in(0,1)^{k^0-1}\times(0,\tau)\times(0,1)^{K-k^0-1}$ and note that the conditional density of~$G'|\{\bar \imath_1,\ldots,\bar \imath_K,F_{K+1}, G_{k^0}=\tau \bz\}$ is
 \begin{align*}
 \mathrm{(const)}\,\mathds 1_{g\in(0,1)^{k^0-1}\times(0,\tau)\times(0,1)^{K-k^0-1}}
 &\prod_{k< k^0}g_k^{K-k}\cosh\big[{\mu^0_{\bar\imath_k}}((1-F_{K+1})\prod_{\ell=1}^kg_\ell+F_{K+1})\big]
 \\&\times \prod_{k> k^0}g_k^{K-k}\cosh\big[{\mu^0_{\bar\imath_k}}((1-F_{K+1})\,\,\bz\!\!\prod_{1\leq \ell\neq k^0\leq k}\!\!g_\ell+F_{K+1})\big]\,,
 \end{align*}
 which, up to some normalising constant, is the same as~\eqref{yoyoz1} up to the following change of support
 \[
 \mathds 1_{g\in(0,1)^K}\leftrightarrow \mathds 1_{g'\in(0,1)^{k^0-1}\times(0,\tau)\times(0,1)^{K-k^0-1}}\,.
 \]
 By an abuse of notation, we denote by~$\bp$ this function, namely
 \begin{align*}
\bp(g)=
 &\prod_{k< k^0}g_k^{K-k}\cosh\big[{\mu^0_{\bar\imath_k}}((1-F_{K+1})\prod_{\ell=1}^kg_\ell+F_{K+1})\big]
 \\&\times \prod_{k> k^0}g_k^{K-k}\cosh\big[{\mu^0_{\bar\imath_k}}((1-F_{K+1})\,\,\bz\!\!\prod_{1\leq \ell\neq k^0\leq k}\!\!g_\ell+F_{K+1})\big]\,.
 \end{align*}
We deduce that
\begin{align}
 \bbF_{z_{k^0}}^{(k^0+1)}(t)\geq \bbF_{z'_{k^0}}^{(k^0+1)}(t)\quad
 &\Leftrightarrow\quad \bar\bbP(G_{k^0+1}\leq t|G_{k^0}=z_{k^0})\geq
 \bar\bbP(G_{k^0+1}\leq t|G_{k^0}=z'_{k^0})\notag\\
 &\Leftrightarrow\quad \bar\bbP(G_{k^0+1}\leq t|G_{k^0}=z_{k^0})\geq
 \bar\bbP(G'_{k^0+1}\leq \tau t|G_{k^0}=\tau z_{k^0})\notag\\
 &\Leftrightarrow\quad \frac{\int_{\cD(t)}\bp}{\int_{\cD(1)}\bp} \geq \frac{\int_{\cD(\tau t)}\bp}{\int_{\cD(\tau)}\bp}\notag\\
&\Leftrightarrow\quad \int_{\cD(t)\times\cD(\tau) }\bp\otimes\bp\geq\int_{\cD(\tau t)\times \cD(1)}\bp\otimes\bp\,,\label{eq:yoyo123}
 \end{align}
where
\[
 \cD(s):=\Big\{(g_1,\ldots,g_{k^0-1},g_{k^0+1}\ldots,g_K)\in(0,1)^{K-1}\ :\ 0< g_{k^0+1}\leq s\Big\}\,.
\]
\begin{figure}
\includegraphics[width=0.4\textwidth]{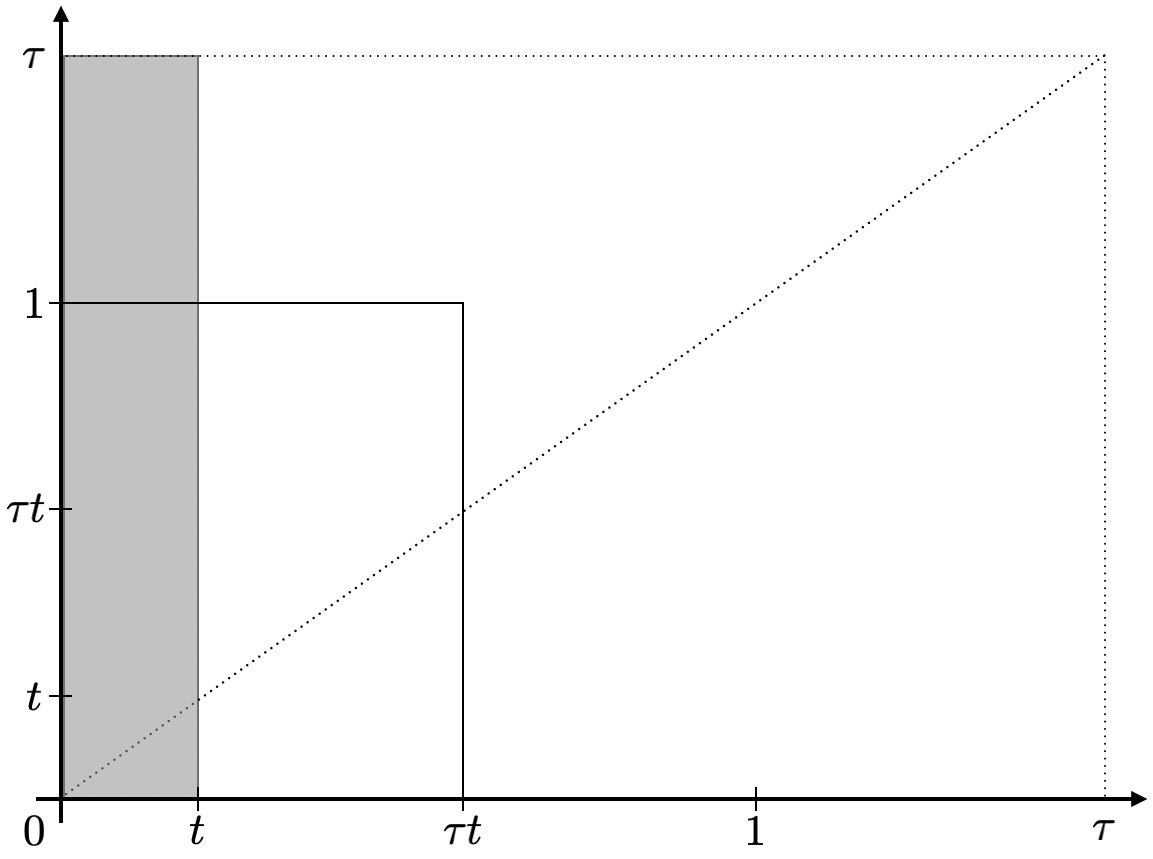}
\caption{The two boxed rectangles have Lebesgue measure, namely~$\tau t$. The~$\bp\otimes\bp$  measure of the grey box is greater than the~$\bp\otimes\bp$ measure of the white box. }\label{fig:yoyo2}
\end{figure}

We now present an inequality on the  to conclude. Observe that we are integrating on domains depicted in Figure~\ref{fig:yoyo2}. The two boxes have same area for the uniform measure and we would like to compare their respective measure for the~$\bp\otimes\bp$ measure. We start by the next lemma whose proof is omitted. 
\begin{lemma}
\label{lem:coshyoyo}
Let~$a,b\geq0$. The function 
\[
z\mapsto\cosh(a\, z+b)\times \cosh(a / z+b)
\]
is non-decreasing on the domain~$[1,\infty)$.
\end{lemma}

\noindent
Now, let~$(g_1,\ldots,g_{k^0-1},g_{k^0+2}\ldots,g_K)\in(0,1)^{K-1}$ be fixed in the integrals~\eqref{eq:yoyo123}. We are the looking at the weights of the domains~$(h_1,h_2)\in(0,t)\times (0,\tau)$ and~$(h_3,h_4)\in(0,\tau t)\times (0,1)$ for the weight function~$w$ given by
\begin{align*}
w(h_1,h_2)=&\mathrm{C}_1h_1^{K-k^0-1}\cosh\big[{\mu^0_{\bar\imath_k}}((1-F_{K+1})\,\,z_{k^0}\!\!\prod_{1\leq \ell <k^0}\!\!g_\ell\times h_1+F_{K+1})\big]\\
&\times \prod_{k> k^0+1}\cosh\big[{\mu^0_{\bar\imath_k}}((1-F_{K+1})\,\,z_{k^0}\!\!\prod_{1\leq \ell\neq k^0,k^0+1\leq k}\!\!g_\ell\times h_1+F_{K+1})\big]\\
&\times h_2^{K-k^0-1}\cosh\big[{\mu^0_{\bar\imath_k}}((1-F_{K+1})\,\,z_{k^0}\!\!\prod_{1\leq \ell <k^0}\!\!g_\ell\times h_2+F_{K+1})\big]\\
&\times \prod_{k> k^0+1}\cosh\big[{\mu^0_{\bar\imath_k}}((1-F_{K+1})\,\,z_{k^0}\!\!\prod_{1\leq \ell\neq k^0,k^0+1\leq k}\!\!g_\ell\times h_2+F_{K+1})\big]\,,
 \end{align*}
 where the constant~$\mathrm{C}_1$ depends on~$(g_1,\ldots,g_{k^0-1},g_{k^0+2}\ldots,g_K)\in(0,1)^{K-1}$. By the change of variables~$h'_1=h_3/t$ and~$h'_2=t h_4$, the right hand term of~\eqref{eq:yoyo123} is given by the integration on the domain~$(h'_1,h'_2)\in(0,t)\times (0,\tau)$ of the weight function~$w'$ given by
\begin{align*}
w'(h_1',h_2')=
 &\mathrm{C}_1{h'_1}^{K-k^0-1}\cosh\big[{\mu^0_{\bar\imath_k}}((1-F_{K+1})\,\,z_{k^0}\!\!\prod_{1\leq \ell <k^0}\!\!g_\ell\times t\times h'_1+F_{K+1})\big]\\
&\times \prod_{k> k^0+1}\cosh\big[{\mu^0_{\bar\imath_k}}((1-F_{K+1})\,\,z_{k^0}\!\!\prod_{1\leq \ell\neq k^0,k^0+1\leq k}\!\!g_\ell\times t\times h'_1+F_{K+1})\big]\\
 &\times {h'_2}^{K-k^0-1}\cosh\big[{\mu^0_{\bar\imath_k}}((1-F_{K+1})\,\,z_{k^0}\!\!\prod_{1\leq \ell <k^0}\!\!g_\ell\times h'_2/t+F_{K+1})\big]\\
&\times \prod_{k> k^0+1}\cosh\big[{\mu^0_{\bar\imath_k}}((1-F_{K+1})\,\,z_{k^0}\!\!\prod_{1\leq \ell\neq k^0,k^0+1\leq k}\!\!g_\ell\times h'_2/t+F_{K+1})\big]\,.
 \end{align*}
Now, invoke Lemma~\ref{lem:coshyoyo} with 
\begin{align*}
a&={\mu^0_{\bar\imath_k}}(1-F_{K+1})\,\,z_{k^0}\!\!\prod_{1\leq \ell <k^0}\!\!g_\ell\times  h\\
b&={\mu^0_{\bar\imath_k}}F_{K+1}\\
z&=t\geq 1\,,
\end{align*}
where~$h=h_1$ or~$h_2$, to get that~$w'\geq w$ and so 
\[
\int_{\cD(t)\times\cD(\tau) }\bp\otimes\bp\geq\int_{\cD(\tau t)\times \cD(1)}\bp\otimes\bp\,,
\]
which concludes this part of the proof.

 \medskip
 
 \noindent
 {$\circ$ \it Step 3.3: } We continue by induction with the other transport maps in Knothe-Rosenblatt's construction. Assume that we have built~$z':=(z'_{k},\ldots,z'_{k^0})$ and~$z:=(z_{k},\ldots,z_{k^0})$ for some~$k>k^0$. Denote~$\mu_{z}^{(k+1)}$ the following conditional law  
 \[
 \mu_{z}^{(k+1)}:=\mathrm{law}\big[G_{k+1}|\{\bar \imath_1,\ldots,\bar \imath_K,F_{K+1}, \underbrace{G_{k}=z_{k},\ldots, G_{k^0}=z_{k^0}}_{\mathrm{denoted\ }G^{[k,k^0]}=z}\}\big]\,,
 \]
 and~$\bbF_{z}^{(k+1)}$ its cdf. Denote~$\mu_{z'}^{(k+1)}$ the following conditional law  
 \[
 \mu_{z'}^{(k+1)}:=\mathrm{law}\big[G_{k+1}|\{\bar \imath_1,\ldots,\bar \imath_K,F_{K+1}, \underbrace{G_{k}=z'_{k},\ldots, G_{k^0}=z'_{k^0}}_{G^{[k,k^0]}=z'}\}\big]\,,
 \]
 and~$\bbF_{z'}^{(k+1)}$ its cdf. Note that~$z\leq z'=T^{(k)}(z)\leq 1$. Again, we would prove that~$\bbF_{z}^{(k+1)}\geq \bbF_{z'}^{(k+1)}$ which implies that the transport map~$T^{(k+1)}:=\big(\bbF_{z'}^{(k+1)}\big)^{-1}\circ \bbF_{z}^{(k+1)}$ satisfies~$T^{(k+1)}(u)\geq u$ for all~$u\in(0,1)$.
 
 For~$\bz\in(0,1)^{k-k^0}\times(x,1)$, the conditional density of~$G|\{\bar \imath_1,\ldots,\bar \imath_K,F_{K+1}, G^{[k,k^0]}=\bz\}$ is
 \begin{align}
 \mathrm{(const)}\,\mathds 1_{g\in(0,1)^K}&\mathds 1_{g^{[k,k^0]}=\bz}
 \prod_{m< k^0}g_m^{K-m}\cosh\big[{\mu^0_{\bar\imath_m}}((1-F_{K+1})\prod_{\ell=1}^mg_\ell+F_{K+1})\big]\notag]\\
 &\times \prod_{k^0\leq m\leq k}\bz_m^{K-m}\cosh\big[{\mu^0_{\bar\imath_m}}((1-F_{K+1})\prod_{1\leq \ell<k^0}\!\!g_\ell\,\,\prod_{n=k^0}^m\bz_n+F_{K+1})\big]\notag\\
&\times \prod_{k<m}g_m^{K-m}\cosh\big[{\mu^0_{\bar\imath_m}}((1-F_{K+1})\prod_{1\leq \ell<k^0}\!\!g_\ell\,\,\prod_{n=k^0}^k\bz_n\prod_{k< \ell\leq m}\!\!g_\ell+F_{K+1})\big]\,.\notag
 \end{align}
Set~$\tau:=\prod_{n=k^0}^kz'_n/\prod_{n=k^0}^kz_n\geq 1$ and~$G'_{k}=\tau G_{k^0+1}$ so that 
 \[
\Big[\prod_{n=k^0}^kz'_n\Big] G_{k+1}=\Big[\prod_{n=k^0}^kz_n\Big] G'_{k+1}\,.
 \]
Then the proof follows the same idea as in {\it Step 3.2} and we will not detail it here.
 \medskip
 
 \noindent
 {$\circ$ \it Step 3.4: } This is the last step of the proof. Assume that we have built~$z':=(z'_{K},\ldots,z'_{k^0})$ and~$z:=(z_{K},\ldots,z_{k^0})$. Denote~$\mu_{z}^{(k^0-1)}$ the following conditional law  
 \[
 \mu_{z}^{(k^0-1)}:=\mathrm{law}\big[G_{k^0-1}|\{\bar \imath_1,\ldots,\bar \imath_K,F_{K+1}, G^{[K,k^0]}=z\}\big]\,,
 \]
 and~$\bbF_{z}^{(k^0-1)}$ its cdf. Denote~$\mu_{z'}^{(k^0-1)}$ the following conditional law  
 \[
 \mu_{z'}^{(k^0-1)}:=\mathrm{law}\big[G_{k^0-1}|\{\bar \imath_1,\ldots,\bar \imath_K,F_{K+1}, G^{[K,k^0]}=z'\}\big]\,,
 \]
 and~$\bbF_{z'}^{(k^01)}$ its cdf. Note that~$z\leq z'=T^{(K)}(z)\leq 1$. Again, we would prove that~$\bbF_{z}^{(k^0-1)}\geq \bbF_{z'}^{(k^0-1)}$ which implies that the transport map~$T^{(k^0-1)}:=\big(\bbF_{z'}^{(k^0-1)}\big)^{-1}\circ \bbF_{z}^{(k^0-1)}$ satisfies~$T^{(k^0-1)}(u)\geq u$ for all~$u\in(0,1)$.

For~$\bz\in(0,1)^{K-k^0}\times(x,1)$, the conditional density of~$G|\{\bar \imath_1,\ldots,\bar \imath_K,F_{K+1}, G^{[K,k^0]}=\bz\}$ is
 \begin{align}
 \mathrm{(const)}\,\mathds 1_{g\in(0,1)^K}&\mathds 1_{g^{[K,k^0]}=\bz}
 \prod_{m< k^0}g_m^{K-m}\cosh\big[{\mu^0_{\bar\imath_m}}((1-F_{K+1})\prod_{\ell=1}^mg_\ell+F_{K+1})\big]\notag\\
 &\times \prod_{k^0\leq m\leq K}\bz_m^{K-m}\cosh\big[{\mu^0_{\bar\imath_m}}((1-F_{K+1})\prod_{1\leq \ell<k^0}\!\!g_\ell\,\,\prod_{n=k^0}^m\bz_n+F_{K+1})\big]\,.\notag
 \end{align}

\noindent
Now, let~$(g_1,\ldots,g_{k^0-2})\in(0,1)^{k^0-2}$ be fixed and denote by
 \begin{align}
\forall g\in(0,1),\quad
w_{z}(g):=&g^{K-k^0+1}\cosh\big[{\mu^0_{\bar\imath_{k^0-1}}}((1-F_{K+1})\prod_{\ell=1}^{k^0-2}g_\ell\times g+F_{K+1})\big]\notag\\
 &\times \prod_{k^0\leq m\leq K}\bz_m^{K-m}\cosh\big[{\mu^0_{\bar\imath_m}}((1-F_{K+1})\prod_{n=k^0}^m\bz_n\prod_{\ell=1}^{k^0-2}g_\ell\times g+F_{K+1})\big]\,.\notag
 \end{align}
 and, substituting~$z$ by~$z'$, define~$w_{z'}$ as well. Let~$t\in(0,1)$. Following the idea of  {\it Step 3.2}, one can check that it is sufficient to prove that
 \[
 \int_0^t\Big(\int_0^1w_{z}(g)w_{z'}(g')\mathrm dg'\Big)\mathrm dg
 \geq
  \int_0^1\Big(\int_0^tw_{z}(g)w_{z'}(g')\mathrm dg'\Big)\mathrm dg\,.
 \]
 Substituting 
 \[
  \int_0^t\Big(\int_0^tw_{z}(g)w_{z'}(g')\mathrm dg'\Big)\mathrm dg
 \]
 on both parts, one is reduced to prove that
  \[
 \int_0^t\Big(\int_t^1w_{z}(g)w_{z'}(g')\mathrm dg'\Big)\mathrm dg
 \geq
\int_0^t\Big(\int_t^1w_{z'}(g)w_{z}(g')\mathrm dg'\Big)\mathrm dg\,.
 \]
 Observe that~$g\leq g'$ in the last two integrals. Now, we have this lemma whose proof is omitted. 
 \begin{lemma}
\label{lem:coshyoyo2}
Let~$0<a\leq a'$ and~$b>0$. The function 
\[
z\mapsto\frac{\cosh(a\, z+b)}{\cosh(a'z+b)}
\]
is non-increasing on the domain~$(0,\infty)$.
\end{lemma}

\noindent
Let~$g\leq g'$. From Lemma~\ref{lem:coshyoyo2}, we deduce that~$\cosh(a g+b)\cosh(a' g'+b)\geq \cosh(a g'+b)\cosh(a' g+b)$, proving that~$w_{z}(g)w_{z'}(g')\geq w_{z'}(g)w_{z}(g')$. It proves that~$T^{(k^0-1)}(u)\geq u$ for all~$u\in(0,1)$.

We then proceed by induction for~$k^0-1\to k^0-2\to\cdots\to 1$. The proof follows the same line as above, {\it Step 3.4}.


\section{A Quasi Monte Carlo (QMC) method: Cubature by lattice rule} 
\label{app:cbc}
Our goal is to  compute  the integral of some function$f$  on  the  hypercube of dimension~$d$, namely
\[
I:= \int_{[0,1]^d}  f(x) dx.
\]
We want to approximate it by a finite sum over~$n$ points 
 \[
I_n:=  \frac 1 n \sum_{i=1}^n f(x^{(i)}).
\] 
A convenient way of constructing the sequence~$x^{(i)}, i=1,\dots,n$ is the so-called {\it lattice rule}: from the first point~$x^{(1)}$ we deduce the others~$x^{(i)}$ by
 \[
x^{(i)} =  \big\{ i . x^{(1)}\big\},
\] 
where  the~$\{\}$ brackets mean that we take the fractional part coordinate by coordinate. In such a case  the error given by
\[
E(f,n,x^{(1)})= I-I_n
\] 
is a function, in particular, of starting point~$x^{(1)}$ . 

The Fast-rank algorithm~\citep{NC}   is a fast algorithm that finds, component  by component and  as a function of the prime~$n$, the sequence of coordinates of~$x^{(1)}$ that minimizes the maximal error when~$f~$ varies in a unit ball~$\mathcal E$ of some {\it RKHS}, namely a tensorial product of {\it Koborov spaces}. In addition it gives an expression  of its minimax error, namely 
\[
 \max_{f\in \mathcal E}(f,n,x^{(1)}).
\]

 In practice, very few properties are known on  the function~$f$, so the result above is not directly applicable. Nevertheless 
 for many functions~$f$, it happens that the convergence  of~$I_n$ to~$I$ is ‘‘{\it fast}'': typically  of the order~$ 1/n$ while the Monte-Carlo method (choosing the~$x^{(i)}$ {\it at random}) converges at rate~$ 1/ \sqrt{n}$. 
 
 A reliable estimate of the error is obtained by adding a {\it Monte-Carlo layer} as in~\cite{genz} for instance. This can be done as follows. Let~$U$ a {\bf unique} uniform variable on~$[0,1]^d$, we define
  \[
x^{(i)}_U:=  \big\{ i. x^{(1)} +U\big\}, \quad I_{n,U}:= \frac1n\sum_{i=1}^n f(x^{(i)}_U).
\] 
Classical computations show that~$I_{n,U}$ is now an unbiased estimator of~$I$.  In a final step,  we perform~$N$ (in practice 15-20) independent  repetitions of the experiment  above an we compute usual asymptotic confidence intervals for independent observations. 

\newpage

 \bibliography{references_all}

\providecommand{\AC}{A.-C}\providecommand{\CA}{C.-A}\providecommand{\CH}{C.-H}\providecommand{\CJ}{C.-J}\providecommand{\JC}{J.-C}\providecommand{\JP}{J.-P}\providecommand{\JB}{J.-B}\providecommand{\JF}{J.-F}\providecommand{\JJ}{J.-J}\providecommand{\JM}{J.-M}\providecommand{\KW}{K.-W}\providecommand{\PL}{P.-L}\providecommand{\RE}{R.-E}\providecommand{\SJ}{S.-J}\providecommand{\XR}{X.-R}\providecommand{\WX}{W.-X}
\begin{thebibliography}{}

\bibitem[Aza{\"\i}s et~al., 2018]{azais2018power}
Aza{\"\i}s, J.-M., De~Castro, Y., and Mourareau, S. (2018).
\newblock Power of the spacing test for least-angle regression.
\newblock {\em Bernoulli}, 24(1):465--492.

\bibitem[Bachoc et~al., 2018]{bachoc2018post}
Bachoc, F., Blanchard, G., Neuvial, P., et~al. (2018).
\newblock On the post selection inference constant under restricted isometry
  properties.
\newblock {\em Electronic Journal of Statistics}, 12(2):3736--3757.

\bibitem[Barber et~al., 2015]{barber2015controlling}
Barber, R.~F., Cand{\`e}s, E.~J., et~al. (2015).
\newblock Controlling the false discovery rate via knockoffs.
\newblock {\em The Annals of Statistics}, 43(5):2055--2085.

\bibitem[Bellec et~al., 2018]{bellec2018slope}
Bellec, P.~C., Lecu{\'e}, G., Tsybakov, A.~B., et~al. (2018).
\newblock Slope meets lasso: improved oracle bounds and optimality.
\newblock {\em The Annals of Statistics}, 46(6B):3603--3642.

\bibitem[Benjamini and Hochberg, 1995]{benjamini1995controlling}
Benjamini, Y. and Hochberg, Y. (1995).
\newblock Controlling the false discovery rate: a practical and powerful
  approach to multiple testing.
\newblock {\em Journal of the Royal statistical society: series B
  (Methodological)}, 57(1):289--300.

\bibitem[Berk et~al., 2013]{berk2013valid}
Berk, R., Brown, L., Buja, A., Zhang, K., Zhao, L., et~al. (2013).
\newblock Valid post-selection inference.
\newblock {\em The Annals of Statistics}, 41(2):802--837.

\bibitem[Bickel et~al., 2009]{bickel2009simultaneous}
Bickel, P.~J., Ritov, Y., Tsybakov, A.~B., et~al. (2009).
\newblock Simultaneous analysis of lasso and dantzig selector.
\newblock {\em The Annals of Statistics}, 37(4):1705--1732.

\bibitem[Blanchard et~al., 2017]{blanchard2017post}
Blanchard, G., Neuvial, P., and Roquain, E. (2017).
\newblock Post hoc inference via joint family-wise error rate control.
\newblock {\em arXiv preprint arXiv:1703.02307}.

\bibitem[Blanchard et~al., 2008]{blanchard2008two}
Blanchard, G., Roquain, E., et~al. (2008).
\newblock Two simple sufficient conditions for fdr control.
\newblock {\em Electronic journal of Statistics}, 2:963--992.

\bibitem[Bogdan et~al., 2015]{bogdan2015slope}
Bogdan, M., Van Den~Berg, E., Sabatti, C., Su, W., and Cand{\`e}s, E.~J.
  (2015).
\newblock Slope---adaptive variable selection via convex optimization.
\newblock {\em The annals of applied statistics}, 9(3):1103.

\bibitem[B{\"u}hlmann and {van de Geer}, 2011]{Buhlmann_vandeGeer11}
B{\"u}hlmann, P. and {van de Geer}, S. (2011).
\newblock {\em {Statistics for high-dimensional data}}.
\newblock {Springer Series in Statistics}. Springer, Heidelberg.
\newblock Methods, theory and applications.

\bibitem[Cand{\`e}s et~al., 2006]{Candes_Romberg_Tao06}
Cand{\`e}s, E.~J., Romberg, J., and Tao, T. (2006).
\newblock {Robust uncertainty principles: Exact signal reconstruction from
  highly incomplete frequency information}.
\newblock {\em {{IEEE} Trans. Inf. Theory}}, 52(2):489--509.

\bibitem[Chen et~al., 1998]{Chen_Donoho_Saunders98}
Chen, S.~S., Donoho, D.~L., and Saunders, M.~A. (1998).
\newblock {Atomic decomposition by basis pursuit}.
\newblock {\em SIAM J. Sci. Comput.}, 20(1):33--61 (electronic).

\bibitem[De~Castro, 2021]{yohann_de_castro_2021_5079768}
De~Castro, Y. (2021).
\newblock github:{ydecastro/lar\_testing: GtSt experiments on real and
  simulated data}, doi:10.5281/zenodo.507976.

\bibitem[Efron et~al., 2004]{efron2004least}
Efron, B., Hastie, T., Johnstone, I., Tibshirani, R., et~al. (2004).
\newblock Least angle regression.
\newblock {\em The Annals of statistics}, 32(2):407--499.

\bibitem[Fithian et~al., 2014]{fithian2014optimal}
Fithian, W., Sun, D., and Taylor, J. (2014).
\newblock Optimal inference after model selection.
\newblock {\em arXiv preprint arXiv:1410.2597}.

\bibitem[Genz, 1992]{genz}
Genz, A. (1992).
\newblock Numerical computation of multivariate normal probabilities.
\newblock {\em Journal of computational and graphical statistics},
  1(2):141--149.

\bibitem[Genz and Bretz, 2009]{genz2009computation}
Genz, A. and Bretz, F. (2009).
\newblock {\em Computation of multivariate normal and t probabilities}, volume
  195.
\newblock Springer Science \& Business Media.

\bibitem[Giraud, 2014]{giraud2014introduction}
Giraud, C. (2014).
\newblock {\em Introduction to high-dimensional statistics}.
\newblock Chapman and Hall/CRC.

\bibitem[Javanmard et~al., 2019]{javanmard2019false}
Javanmard, A., Javadi, H., et~al. (2019).
\newblock False discovery rate control via debiased lasso.
\newblock {\em Electronic Journal of Statistics}, 13(1):1212--1253.

\bibitem[Javanmard and Montanari, 2014]{javanmard2014confidence}
Javanmard, A. and Montanari, A. (2014).
\newblock Confidence intervals and hypothesis testing for high-dimensional
  regression.
\newblock {\em The Journal of Machine Learning Research}, 15(1):2869--2909.

\bibitem[Lockhart et~al., 2014]{lockhart2014significance}
Lockhart, R., Taylor, J., Tibshirani, R.~J., and Tibshirani, R. (2014).
\newblock {A significance test for the lasso}.
\newblock {\em Annals of statistics}, 42(2):413.

\bibitem[Nuyens and Cools, 2006]{NC}
Nuyens, D. and Cools, R. (2006).
\newblock Fast algorithms for component-by-component construction of rank-1
  lattice rules in shift-invariant reproducing kernel hilbert spaces.
\newblock {\em Mathematics of Computation}, 75(254):903--920.

\bibitem[Rhee et~al., 2006]{rhee2006genotypic}
Rhee, S.-Y., Taylor, J., Wadhera, G., Ben-Hur, A., Brutlag, D.~L., and Shafer,
  R.~W. (2006).
\newblock Genotypic predictors of human immunodeficiency virus type 1 drug
  resistance.
\newblock {\em Proceedings of the National Academy of Sciences},
  103(46):17355--17360.

\bibitem[Roquain, 2011]{roquain2010type}
Roquain, E. (2011).
\newblock Type i error rate control for testing many hypotheses: a survey with
  proofs.
\newblock {\em Journal de la Soci{\'e}t{\'e} Fran{\c c}aise de Statistique},
  152(2):3--38.

\bibitem[Santambrogio, 2015]{santambrogio2015optimal}
Santambrogio, F. (2015).
\newblock Optimal transport for applied mathematicians.
\newblock {\em Birk{\"a}user, NY}, 55:58--63.

\bibitem[Taylor and Tibshirani, 2015]{taylor2015statistical}
Taylor, J. and Tibshirani, R.~J. (2015).
\newblock Statistical learning and selective inference.
\newblock {\em Proceedings of the National Academy of Sciences},
  112(25):7629--7634.

\bibitem[Tian et~al., 2018]{tian2018selective}
Tian, X., Loftus, J.~R., and Taylor, J.~E. (2018).
\newblock Selective inference with unknown variance via the square-root lasso.
\newblock {\em Biometrika}, 105(4):755--768.

\bibitem[Tibshirani, 1996]{tibshirani1996regression}
Tibshirani, R. (1996).
\newblock Regression shrinkage and selection via the lasso.
\newblock {\em Journal of the Royal Statistical Society: Series B
  (Methodological)}, 58(1):267--288.

\bibitem[Tibshirani et~al., 2015]{hastie2015statistical}
Tibshirani, R., Wainwright, M., and Hastie, T. (2015).
\newblock {\em Statistical Learning with Sparsity: The Lasso and
  Generalizations}.
\newblock Monographs on Statistics \& Applied Probability. Chapman and Hall/CRC
  press.

\bibitem[Tibshirani et~al., 2016]{taylor2014exact}
Tibshirani, R.~J., Taylor, J., Lockhart, R., and Tibshirani, R. (2016).
\newblock Exact post-selection inference for sequential regression procedures.
\newblock {\em Journal of the American Statistical Association},
  111(514):600--620.

\bibitem[van~de Geer, 2016]{van2016estimation}
van~de Geer, S. (2016).
\newblock Estimation and testing under sparsity.
\newblock {\em Lecture Notes in Mathematics}, 2159.

\bibitem[Villani, 2008]{villani2008optimal}
Villani, C. (2008).
\newblock {\em Optimal transport: old and new}, volume 338.
\newblock Springer Science \& Business Media.

\bibitem[Wainwright, 2009]{wainwright2009sharp}
Wainwright, M.~J. (2009).
\newblock Sharp thresholds for high-dimensional and noisy sparsity recovery
  using l1-constrained quadratic programming (lasso).
\newblock {\em IEEE transactions on information theory}, 55(5):2183--2202.

\end{thebibliography}

\end{document}